\documentclass[11pt]{article}

\usepackage{amsmath,amssymb,amsfonts}
\usepackage{pstricks,pst-node}
\usepackage{xr}

\numberwithin{equation}{section}
\newtheorem{thm}{Theorem}[section] 
\newtheorem{prp}[thm]{Proposition}
\newtheorem{lmm}[thm]{Lemma}   
\newtheorem{crl}[thm]{Corollary} 
\newtheorem{dfn}[thm]{Definition}

\renewcommand{\frak}{\mathfrak}
\renewcommand{\Bbb}{\mathbb}
\def\A{\mathcal A}
\def\C{\mathbb C}
\def\cC{\mathcal C}

\def\D{\mathfrak D}
\def\cD{\mathcal D}
\def\E{\mathbb E}

\def\F{\mathfrak F}

\def\H{\mathcal H}
\def\I{\mathfrak i}
\def\cI{\mathcal I}
\def\L{\mathfrak L}
\def\M{\mathfrak M}
\def\cM{\mathcal M}
\def\R{\mathbb R}
\def\O{\mathcal O}
\def\Q{\mathbb Q}
\def\T{\mathcal T}
\def\U{\mathcal U}
\def\V{\mathcal V}
\def\X{\mathfrak X}
\def\Z{\mathbb Z}

\def\e_ref#1{(\ref{#1})}
\def\under#1{\underline{#1}}
\def\ov#1{\overline{#1}}
\def\wt#1{\widetilde{#1}}
\def\ti#1{\tilde{#1}}

\def\lra{\longrightarrow}
\def\Lra{\Longrightarrow}

\def\Llra{\Longleftrightarrow}

\def\lan{\langle}
\def\ran{\rangle}
\def\lr#1{\lan{#1}\ran}
\def\blr#1{\big\lan{#1}\big\ran}
\def\llrr#1{\lan\!\lan{#1}\ran\!\ran}

\def\al{\alpha}
\def\be{\beta}
\def\de{\delta}
\def\ep{\epsilon}
\def\ga{\gamma}
\def\io{\iota}

\def\na{\nabla}
\def\om{\omega}

\def\th{\theta}
\def\ups{\upsilon}
\def\ve{\varepsilon}
\def\vph{\varphi}

\def\ze{\zeta}

\def\Ga{\Gamma}
\def\La{\Lambda}
\def\Om{\Omega}
\def\Si{\Sigma}

\def\Aut{\textnormal{Aut}}

\def\const{\textnormal{const}}
\def\ev{\textnormal{ev}}
\def\id{\textnormal{id}}
\def\ind{\textnormal{ind}~\,}

\def\rk{\textnormal{rk}~\,}

\def\st{~\textnormal{s.t.}~}

\def\Hom{\textnormal{Hom}}
\def\Gr{\textnormal{Gr}}

\def\P{\Bbb{P}^n}
\def\bP{\Bbb{P}}
\def\PP{\Bbb{P}^2}
\def\PPP{\Bbb{P}^3}
\def\Pf{\Bbb{P}^4}
\def\i{\infty}
\def\eset{\emptyset}
\def\bpar{\bar\partial}

\begin{document}

\title{On the Structure of Certain Natural Cones\linebreak 
over Moduli Spaces of Genus-One Holomorphic Maps}

\author{Aleksey Zinger\thanks{Partially supported by an NSF Postdoctoral Fellowship}}

\date{\today}
\maketitle

\begin{abstract}
\noindent
We show that certain naturally arising cones over the main component of a moduli space
of $J_0$-holomorphic maps into $\P$ have a well-defined euler class.
We also prove that this is the case if the standard complex structure~$J_0$ on~$\P$
is replaced by a nearby almost complex structure~$J$.
The genus-zero analogue of the cone considered in this paper is a vector bundle.
The genus-zero Gromov-Witten invariant of a projective complete intersection
can be viewed as the euler class of such a vector bundle.
As shown in a separate paper, this is also the case for the ``genus-one part" of 
the genus-one GW-invariant.  
The remaining part is a multiple of the genus-zero GW-invariant.
\end{abstract}

\thispagestyle{empty}

\tableofcontents

\section{Introduction}
\label{intro_sec}

\subsection{Motivation}
\label{motiv_subs}

\noindent
The GW-invariants of symplectic manifolds have been an area of 
much research in the past decade.
These invariants are however often hard to compute.\\

\noindent
If $Y$ is a compact Kahler submanifold of the complex projective space $\P$,
one could try to compute the GW-invariants of~$Y$ by relating them to
the GW-invariants of~$\P$.
For example, suppose $Y$ is a  hypersurface in~$\P$ of degree~$a$. 
In other words, if $\ga\!\lra\!\P$ is the tautological  line bundle and
$\L\!=\!\ga^{*\otimes a}\!\lra\!\P$, then
$$Y=s^{-1}(0),$$
for some $s\!\in\!H^0(\P;\L)$ such that $s$ is transverse to the zero set.
If $g$, $k$, and~$d$ are nonnegative integers, let 
$\ov\M_{g,k}(\P,d)$ and $\ov\M_{g,k}(Y,d)$ denote the moduli spaces of stable $J_0$-holomorphic 
\hbox{degree-$d$} maps from genus-$g$ Riemann surfaces with $k$ marked points
to~$\P$ and~$Y$, respectively.
These moduli spaces determine the genus-$g$ degree-$d$ GW-invariants of $\P$ and~$Y$.\\

\noindent
By definition, the moduli space $\ov\M_{g,k}(Y,d)$ is a subset
of the moduli space $\ov\M_{g,k}(\P,d)$. In fact,
\begin{equation}\label{modulirel_e1}
\ov\M_{g,k}(Y,d)=\big\{[\cC,u]\!\in\!\ov\M_{g,k}(\P,d)\!:
s\!\circ\!u=0\in H^0\big(\cC;u^*\L\big)\big\}.
\end{equation}
Here $[\cC,u]$ denotes the equivalence class of the holomorphic map 
$u\!:\cC\!\lra\!\P$ from a genus-$g$ curve~$\cC$ with $k$ marked points.
The relationship~\e_ref{modulirel_e1} can be restated more globally as follows. 
Let 
$$\pi_{g,k}^d\!\!:{\frak U}_{g,k}(\P,d)\lra\ov\M_{g,k}(\P,d)$$ 
be 
the semi-universal family and let
$$\ev_{g,k}^d\!:{\frak U}_{g,k}(\P,d) \lra \P$$ 
be the natural evaluation~map.
In other words, the fiber of $\pi_{g,k}^d$ over $[\cC,u]$ is the curve~$\cC$ 
with $k$ marked points, while 
$$\ev_{g,k}^d\big([\cC,u;z]\big)=u(z)
\qquad\hbox{if}\quad z\!\in\!\cC.$$
We define the section~$s_{g,k}^d$ of the sheaf 
$\pi_{g,k*}^d\ev_{g,k}^{d*}\L\!\lra\!\ov\M_{g,k}(\P,d)$ by
$$s_{g,k}^d\big([\cC,u]\big)=[s\circ u].$$
By~\e_ref{modulirel_e1}, $\ov\M_{g,k}(Y,d)$ is the zero set of this section.\\

\noindent
The previous paragraph suggests that it should be possible
to relate the genus-$g$ degree-$d$ GW-invariants of the hypersurface~$Y$
to the moduli space $\ov\M_{g,k}(\P,d)$ in general and to the~sheaf 
$$\pi_{g,k*}^d\ev_{g,k}^{d*}\L\lra\ov\M_{g,k}(\P,d)$$ 
in particular. 
In~fact, it can be shown that
\begin{equation}\label{genus0_e}
\hbox{GW}_{0,k}^Y(d;\psi)\!\equiv\!
\blr{\psi,\big[\ov\M_{0,k}(Y,d)\big]^{vir}}
=\blr{\psi\cdot e\big(\pi_{0,k*}^d\ev_{0,k}^{d*}\L\big),
\big[\ov\M_{0,k}(\P,d)\big]}
\end{equation}
for all $\psi\!\in\!H^*(\ov\M_{0,k}(\P,d);\Q)$.
The moduli space $\ov\M_{0,k}(\P,d)$ is a smooth orbivariety  and
\begin{equation}\label{g0sheaf_e}
\pi_{0,k*}^d\ev_{0,k}^{d*}\L\lra\ov\M_{0,k}(\P,d)
\end{equation}
is a locally free sheaf, i.e.~a vector bundle.
Furthermore,
\begin{gather*}
\dim_{\C}\ov\M_{0,k}(\P,d)=d(n\!+\!1)+(n\!-\!3)+k, \qquad
\rk\!\!_{\C}\,\pi_{0,k*}^d\ev_{0,k}^{d*}\L=da+1,\\ 
\hbox{and}\qquad
\dim_{\C}^{vir}\ov\M_{0,k}(Y,d)=d(n\!+\!1-a)+(n\!-\!1-3)+k.
\end{gather*}
Thus, the right-hand side of~\e_ref{genus0_e} is well-defined
and vanishes for dimensional reasons precisely when the left-hand side of
\e_ref{genus0_e} does.
In other cases, the right-hand side of~\e_ref{genus0_e}
can be computed via the classical localization theorem of~\cite{AB},
though the complexity of this computation increases rapidly with the degree~$d$.\\

\noindent
If $g\!>\!0$, the sheaf 
$\pi_{g,k*}^d\ev_{g,k}^{d*}\L\!\lra\!\ov\M_{g,k}(\P,d)$ 
is not locally free and does not define an euler class.
Thus, the right-hand side of~\e_ref{genus0_e} does not even
make sense if $0$ is replaced by $g\!>\!0$.
Instead one might try to generalize~\e_ref{genus0_e} as
\begin{equation}\label{allgenus_e}\begin{split}
\hbox{GW}_{g,k}^Y(d;\psi)
&\equiv\blr{\psi,\big[\ov\M_{g,k}(Y,d)\big]^{vir}}\\
&\stackrel{?}{=}\blr{ \psi\cdot
e\big(R^0\pi_{g,k*}^d\ev_{g,k}^{d*}\L-R^1\pi_{g,k*}^d\ev_{g,k}^{d*}\L\big),
\big[\ov\M_{g,k}(\P,d)\big]^{vir}},
\end{split}\end{equation}
where $R^i\pi_{g,k*}^d\ev_{g,k}^{d*}\L\!\lra\!\ov\M_{g,k}(\P,d)$
is the $i$th direct image sheaf.
The right-hand side of~\e_ref{allgenus_e} can be computed via 
the virtual localization theorem of~\cite{GP1}.
However, 
$$N_1(d)\!\equiv\!\hbox{GW}_{g,0}^Y(d;1)\neq
\blr{e\big(R^0\pi_{1,0*}^d\ev_{1,0}^{d*}\L\!-\!
R^1\pi_{1,0*}^d\ev_{1,0}^{d*}\L\big),
\big[\ov\M_1(\Pf,d)\big]^{vir}},$$
according to a low-degree check of~\cite{GP2} and~\cite{K}
for a quintic threefold $Y\!\subset\!\Pf$.\\

\noindent
As shown in~\cite{LZ}, a $g\!=\!1$ analogue of the role played by the euler class of sheaf~\e_ref{g0sheaf_e} is played by the euler class of the~sheaf
\begin{equation}\label{g1sheaf_e}
\pi_{1,k*}^d\ev_{1,k}^{d*}\L\lra\ov\M_{1,k}^0(\P,d),
\end{equation}
where $\ov\M_{1,k}^0(\P,d)$ is the primary, algebraically irreducible, component
of $\ov\M_{1,k}(\P,d)$.
In other words, $\ov\M_{1,k}^0(\P,d)$ is the closure in $\ov\M_{1,k}(\P,d)$,
either in the stable-map or Zariski topology,
of the subspace
$$\M_{1,k}^0(\P,d)=\big\{[{\cal C},u]\!\in\!\ov\M_{1,k}(\P,d)\!:
{\cal C}~\hbox{is smooth}\big\}.$$
One of the results of this paper is that the euler class of 
the sheaf \e_ref{g1sheaf_e} is in fact well-defined:

\begin{thm}
\label{thm1}
If $n$, $d$, and $a$ are positive integers, $k$ is a nonnegative integer,
$\L\!=\!\ga^{*\otimes a}\!\lra\!\P$,
$$\pi_{1,k}^d\!:{\frak U}_{1,k}(\P,d)\lra\ov\M_{1,k}^0(\P,d)$$ 
is the semi-universal family, and 
$$\ev_{1,k}^d\!:{\frak U}_{1,k}(\P,d)\lra\P$$ 
is the natural evaluation~map, the sheaf
$$\pi_{1,k*}^d\ev_{1,k}^{d*}\L\lra\ov\M_{1,k}^0(\P,d)$$
determines a homology class and a cohomology class on $\ov\M_{1,k}^0(\P,d)$:
\begin{gather*}
\hbox{PD}_{\ov\M_{1,k}^0(\P,d)}\big(e(\pi_{1,k*}^d\ev_{1,k}^{d*}\L)\big)
\in H_{2(d(n+1-a)+k)}\big(\ov\M_{1,k}^0(\P,d);\Q\big)\\
\hbox{and}\qquad
e(\pi_{1,k*}^d\ev_{1,k}^{d*}\L)\in H^{2da}(\ov\M_{1,k}^0(\P,d);\Q\big).
\end{gather*}\\
\end{thm}

\noindent
{\it Remark:}
If $a_1,\ldots,a_m\!\in\!\Bbb{Z}^+$ and 
$$\L=\ga^{*\otimes a_1}\!\oplus\!\ldots\!\oplus\!\ga^{*\otimes a_m}\lra\P,$$
then the sheaf $\pi_{1,k*}^d\ev_{1,k}^{d*}\L$ is the direct sum of the sheaves
corresponding to the line bundles~$\ga^{*\otimes a_l}$.
Thus, Theorem~\ref{thm1} applies to any split vector bundle over~$\P$.\\

\noindent
One way to view the statement of this theorem is that the sheaf~\e_ref{g1sheaf_e}
admits a desingularization, and the euler class of every desingularization
of~\e_ref{g1sheaf_e} is the same, in the appropriate sense.
This is not the point of view taken in this paper.
However, one approach to computing the number 
\begin{equation}\label{thm1_e}
\blr{\psi,\hbox{PD}_{\ov\M_{1,k}^0(\P,d)}\big(e(\pi_{1,k*}^d\ev_{1,k}^{d*}\L)\big)}
=\blr{\psi\cdot e\big(\pi_{1,k*}^d\ev_{1,k}^{d*}\L\big),
\big[\ov\M_{1,k}^0(\P,d)\big]}
\end{equation}
for a natural cohomology class $\psi\!\in\!H^*(\ov\M_{1,k}(\P,d);\Q)$ 
is to apply the localization theorem of \cite{AB} to
a desingularization of~\e_ref{g1sheaf_e}.
In~\cite{VZ}, we construct a desingularization of the space~$\ov\M_{1,k}^0(\P,d)$, 
i.e.~a smooth orbivariety $\wt\M_{1,k}^0(\P,d)$ and a~map
$$\ti\pi\!:\wt\M_{1,k}^0(\P,d)\lra\ov\M_{1,k}^0(\P,d),$$
which is biholomorphic onto $\M_{1,k}^0(\P,d)$.
This desingularization of $\ov\M_{1,k}^0(\P,d)$ comes with a desingularization
of the sheaf~\e_ref{g1sheaf_e},
i.e.~a vector bundle 
$$\ti{\cal V}_{1,k}^d\lra\wt\M_{1,k}^0(\P,d)
\qquad\hbox{s.t.}\qquad
\tilde{\pi}_*\tilde{\cal V}_{1,k}^d=\pi_{1,k*}^d\ev_{1,k}^{d*}\L.$$
In particular,
\begin{equation}\label{design_e}
\blr{\psi\cdot e\big(\pi_{1,k*}^d\ev_{1,k}^{d*}\L\big),
\big[\ov\M_{1,k}^0(\P,d)\big]}=
\blr{\tilde{\pi}^*\psi\cdot e(\tilde{\cal V}_{1,k}^d),
\big[\wt\M_{1,k}^0(\P,d)\big]}.
\end{equation}
Since a group action on $\P$ induces actions on $\wt\M_{1,k}^0(\P,d)$
and on $\tilde{\cal V}_{1,k}^d$,
the localization theorem of~\cite{AB} is directly applicable
to the right-hand side of~\e_ref{design_e}, for a natural cohomology class~$\psi$.\\ 

\noindent
Before the results of~\cite{LZ} were announced, no positive-genus analogue of~\e_ref{genus0_e}
had been even conjectured.
On the other hand, Theorem~\ref{thm1} suggests a natural genus-one analogue of~\e_ref{genus0_e},
which is proved in~\cite{LZ}, and a conjectural extension of~\e_ref{genus0_e} to higher
genera, which is stated in~\cite{LZ}.\\

\noindent
Theorem~\ref{thm1} is the $J\!=\!J_0$ case of Theorem~\ref{main_thm},
which is stated in Subsection~\ref{mainres_subs}.
In the next subsection, we describe the main topological arguments that lie 
behind the proof of Theorems~\ref{thm1} and~\ref{main_thm}.

\subsection{General Approach}
\label{appr_subs}

\noindent
In this paper, we view Theorem~\ref{thm1} as a statement about a certain
(orbi-)cone 
\begin{equation}\label{cone_e}
\pi\!:{\cal V}_{1,k}^d\lra\ov\M_{1,k}^0(\P,d).
\end{equation}
In other words, $\pi$ is a continuous map between topological spaces
such that every fiber of $\pi$ is a vector space, up to a quotient by a finite group,
and the vector space operations are continuous.
In this~case,
$${\cal V}_{1,k}^d\big|_{[{\cal C},u]}\!\equiv\!\pi^{-1}\big([{\cal C},u]\big)
=H^0({\cal C};u^*\L)\big/\hbox{Aut}({\cal C},u).$$
Furthermore, the space $\ov\M_{1,k}^0(\P,d)$ is stratified by smooth orbifolds,
and the restriction of ${\cal V}_{1,k}^d$ to every stratum 
of $\ov\M_{1,k}^0(\P,d)$ is a smooth vector bundle.
We will show that the cone~\e_ref{cone_e} admits a continuous multisection~$\vph$
such~that\\
${}\quad$ (${\cal V}1$) $\vph|_{\M_{1,k}^0(\P,d)}$ is smooth and transverse to
the zero~set in ${\cal V}_{1,k}^d|_{\M_{1,k}^0(\P,d)}$, and\\
${}\quad$ (${\cal V}2$) the intersection of $\vph^{-1}(0)$ with a boundary stratum of 
$\ov\M_{1,k}^0(\P,d)$ is a smooth suborbifold\\
${}\qquad\quad~$ of the stratum of real dimension of at most $2(d(n\!+\!1\!-\!a)\!+\!k)\!-\!2$.\\
These two properties, along with the structure of the space $\ov\M_{1,k}^0(\P,d)$,
imply that $\vph^{-1}(0)$ determines an element of 
$H_{2d(n+1-a)+2k}(\ov\M_{1,k}^0(\P,d);\Bbb{Q})$.
We will also show that for any two continuous sections $\vph_0$ and $\vph_1$
of~\e_ref{cone_e} satisfying~(${\cal V}1$) and~(${\cal V}2$),
there exists a continuous homotopy
$$\Phi\!:[0,1]\!\times\!\ov\M_{1,k}^0(\P,d)\lra{\cal V}_{1,k}^d$$
such that $\Phi|_{\{t\}\times\ov\M_{1,k}^0(\P,d)}\!=\!\vph_t$ for $t\!=\!0,1$,\\
${}\quad$ (${\cal V}1'$) $\Phi|_{[0,1]\times\M_{1,k}^0(\P,d)}$ is smooth and 
transverse to the zero~set in $[0,1]\!\times\!{\cal V}_{1,k}^d|_{\M_{1,k}^0(\P,d)}$, and\\
${}\quad$ (${\cal V}2'$) the intersection of $\Phi^{-1}(0)$
with a stratum of $[0,1]\!\times\!\ov\M_{1,k}^0(\P,d)$ is a smooth suborbifold\\
${}\qquad\quad~~$ of the stratum of real dimension of at most 
$2(d(n\!+\!1\!-\!a)\!+\!k)\!-\!1$.\\
The existence of such a homotopy implies that 
$$[\vph_0^{-1}(0)]=[\vph_1^{-1}(0)]\in H_{2(d(n+1-a)+k)}(\ov\M_{1,k}^0(\P,d);\Bbb{Q}).$$
We call this homology class {\it the Poincare dual of the euler class}
of the cone~\e_ref{cone_e} and of the sheaf~\e_ref{g1sheaf_e}. \\

\noindent
The key fact we use in constructing a section $\vph$ of~\e_ref{cone_e}
satisfying~(${\cal V}1$) and~(${\cal V}2$) is Proposition~\ref{g1conebdstr_prp}.
This proposition implies in particular  that on a neighborhood of each boundary
stratum of $\ov\M_{1,k}^0(\P,d)$ the cone ${\cal V}_{1,k}^d$ contains 
a complex vector bundle of a rank sufficiently high to insure 
that a generic multisection~$\vph$ of this bundle satisfies~(${\cal V}2$).\\ 

\noindent
Any homology class $X\!\in\!H_{2da}(\ov\M_{1,k}^0(\P,d);\Q)$ can be represented by
a {\tt pseudocycle}
$$f_X\!:M_X\lra\ov\M_{1,k}^0(\P,d).$$
Here $M_X$ is a compact topological space which is stratified 
by smooth orbifolds, such that the main stratum~$M_X^0$ of $M_X$ 
is an oriented orbifold of real dimension~$2da$, while
the complement of $M_X^0$ in $M_X$  is a union of orbifolds of
real dimension of at most $2da\!-\!2$, and
$f_X$ is a continuous map such that the restriction 
of $f_X$ to each stratum of $M_X$ is smooth.
In particular, every stratum of $M_X$ is mapped into a stratum of $\ov\M_{1,k}^0(\P,d)$;
see Chapter~7 of~\cite{MS} for a discussion of pseudocycles in the basic manifold case.
If $\vph$ is a section of~\e_ref{cone_e} satisfying~(${\cal V}1$)
and~(${\cal V}2$), 
we can also require that\\
${}\quad$ ($\vph X1$) $f_X(M_X)\!\cap\!\vph^{-1}(0)\subset\M_{1,k}^0(\P,d)$,
$f_X^{-1}\big(\vph^{-1}(0)\big)\subset M_X^0$;\\
${}\quad$ ($\vph X2$) $f_X|_{M_X^0}$ intersects $\vph^{-1}(0)$
transversally in~$\M_{1,k}^0(\P,d)$.\\
These assumptions imply that $\vph^{-1}(0)\!\cap\!f_X(M_X^0)$
is a compact oriented zero-dimensional suborbifold of $\M_{1,k}^0(\P,d)$. 
We then set
\begin{equation}\label{eulerclassdfn_e}
\blr{e(\pi_{1,k*}^d\ev_{1,k}^{d*}\L),X}=
~{^{\pm}}\big|\vph^{-1}(0)\!\cap\!f_X(M_X^0)\big|,
\end{equation}
where $^{\pm}\!|{\cal Z}|$ denotes the cardinality of 
a compact oriented zero-dimensional orbifold ${\cal Z}$,
i.e.~the number of elements in~the finite set ${\cal Z}$ counted 
with the appropriate multiplicities.\\

\noindent
If $f_{X,0}\!: M_{X,0}\!\lra\!\ov\M_{1,k}^0(\P,d)$ and
$f_{X,1}\!:M_{X,1}\!\lra\!\ov\M_{1,k}^0(\P,d)$
are two pseudocycles satisfying~($\vph X1$) and~($\vph X2$),
we can choose a pseudocycle equivalence 
$$F\!: \ti{M}\lra\ov\M_{1,k}^0(\P,d)$$
between $f_{X,0}$ and $f_{X,1}$ such that\\
${}\quad$ ($\vph X1'$) $F(\ti{M})\!\cap\!\vph^{-1}(0)\subset\M_{1,k}^0(\P,d)$,
$F^{-1}\big(\vph^{-1}(0)\big)\subset\ti{M}^0$;\\
${}\quad$ ($\vph X2'$) $F|_{\ti{M}^0}$ intersects $\vph^{-1}(0)$
transversally in~$\M_{1,k}^0(\P,d)$.\\
These two assumptions imply that $\vph^{-1}(0)\!\cap\!F(\ti{M}^0)$
is a compact oriented one-dimensional suborbifold of $\M_{1,k}^0(\P,d)$ and
\begin{gather*}
\partial\big(\vph^{-1}(0)\!\cap\!F(\ti{M}^0)\big)
=\vph^{-1}(0)\!\cap\!f_{X,1}(M_{X,1}^0)-
\vph^{-1}(0)\!\cap\!f_{X,0}(M_{X,0}^0)\\
\Lra\qquad
^{\pm}\!\big|\vph^{-1}(0)\!\cap\!f_{X,0}(M_{X,0}^0)\big|
=~^{\pm}\!\big|\vph^{-1}(0)\!\cap\!f_{X,1}(M_{X,1}^0)\big|.
\end{gather*}
Thus, the number in~\e_ref{eulerclassdfn_e} is independent of 
the choice of pseudocycle representative~$f_X$ 
for $X$ satisfying~($\vph X1$) and~($\vph X2$).\\

\noindent
Similarly, if $\vph_0$ and $\vph_1$ are two multisections satisfying
(${\cal V}1$) and~(${\cal V}2$), let $\Phi$ be a homotopy between
$\vph_0$ and $\vph_1$ satisfying (${\cal V}1'$) and~(${\cal V}2'$).
We can then choose a pseudocycle representative 
$$f_X\!:M_X\lra\ov\M_{1,k}^0(\P,d)$$
for $X$ such that\\
${}\quad$ ($\Phi X1$) $f_X(M_X)\!\cap\!\Phi^{-1}(0)\subset\M_{1,k}^0(\P,d)$,
$f_X^{-1}\big(\Phi^{-1}(0)\big)\subset M_X^0$;\\
${}\quad$ ($\Phi X2$) $f_X|_{M_X^0}$ intersects $\Phi^{-1}(0)$
transversally in~$\M_{1,k}^0(\P,d)$,\\
and $f_X$ satisfies  ($\vph X2$) with $\vph\!=\!\vph_0$ and $\vph\!=\!\vph_1$.
These assumptions imply that $\Phi^{-1}(0)\!\cap\!f_X(M_X^0)$
is a compact oriented one-dimensional suborbifold of $\M_{1,k}^0(\P,d)$ and
\begin{gather*}
\partial\big(\Phi^{-1}(0)\!\cap\!f_X(M_X^0)\big)
=\vph_1^{-1}(0)\!\cap\!f_X(M_X^0)-\vph_0^{-1}(0)\!\cap\!f_X(M_X^0)\\
\Lra\qquad
^{\pm}\!\big|\vph_0^{-1}(0)\!\cap\!f_X(M_X^0)\big|
=~^{\pm}\!\big|\vph_1^{-1}(0)\!\cap\!f_X(M_X^0)\big|.
\end{gather*}
Thus, the number in~\e_ref{eulerclassdfn_e} is independent of 
the choice of section~$\vph$ satisfying~(${\cal V}1$) and~(${\cal V}2$).
We conclude that \e_ref{eulerclassdfn_e} defines an element~of
$$\hbox{Hom}\big(H_{2da}(\ov\M_{1,k}^0(\P,d);\Q);\Q)
=H^{2da}(\ov\M_{1,k}^0(\P,d);\Q).$$
We call this cohomology class {\tt the euler class}
of the cone~\e_ref{cone_e} and of the sheaf~\e_ref{g1sheaf_e}.\\

\noindent
We note that the existence of a continuous section $\vph$ of~\e_ref{cone_e}
satisfying~(${\cal V}1$) and~(${\cal V}2$) 
implies that the euler class of every desingularization of~\e_ref{cone_e},
or of~\e_ref{g1sheaf_e}, is the same, in the appropriate sense,
for the following reason.
If 
\begin{equation}\label{desing_e}
\begin{array}{ccc}
\tilde{\cal V}_{1,k}^d& \psline{->}(-1,0)(1,0)\rput(0,.25){\tilde{\pi}_*}& {\cal V}_{1,k}^d\\
\psline{->}(0,0)(0,-1.5)& \hspace{.6in}& \psline{->}(0,0)(0,-1.5)\\
&&\\ &&\\ &&\\
\wt\M_{1,k}^0(\P,d)&\psline{->}(-.8,0.1)(.8,0.1)\rput(0,.35){\tilde{\pi}}&\ov\M_{1,k}^0(\P,d)
\end{array}
\end{equation}
is a desingularization of the cone~\e_ref{cone_e}, or of the sheaf~\e_ref{g1sheaf_e},
the section $\vph$ induces a section $\ti\vph$ of the vector bundle
$$\ti{\cal V}_{1,k}^d\lra\wt\M_{1,k}^0(\P,d)$$
such that $\ti\vph\!=\!\vph$ on $\M_{1,k}^0(\P,d)$ and 
$\ti\vph^{-1}(0)\!-\!\M_{1,k}^0(\P,d)$ is a finite union
of smooth orbifolds of real dimension of at most $2(d(n\!+\!1\!-\!a)\!+\!k)\!-\!2$.
Suppose $X\!\in\!H_{2da}(\ov\M_{1,k}^0(\P,d);\Q)$
is represented by a pseudocycle 
$$f_X\!:M_X\lra\ov\M_{1,k}^0(\P,d),$$
and 
\begin{equation*}\begin{split}
\psi_X\!\equiv\!\hbox{PD}_{\ov\M_{1,k}^0(\P,d)}X \in
& H^{2(d(n+1-a)+k)}(\ov\M_{1,k}^0(\P,d);\Bbb{Q})\\
&\qquad =\hbox{Hom}\big(H_{2(d(n+1-a)+k)}(\ov\M_{1,k}^0(\P,d);\Q);\Q)
\end{split}\end{equation*}
is the {\tt Poincare dual of $X$}, i.e.~the element constructed 
by intersecting $2(d(n\!+\!1\!-\!a)\!+\!k)$-pseudocycles with~$f_X(M_X)$.
The Poincare dual of the cohomology class $\ti\pi^*\psi_X$ in $\wt\M_{1,k}^0(\P,d)$
can then be represented by a pseudocycle 
\begin{gather*}
f_{\ti{X}}\!:M_{\ti{X}}\lra\wt\M_{1,k}^0(\P,d) \qquad\hbox{s.t.}\\
M_X^0\subset M_{\ti{X}}^0,\quad
f_{\ti{X}}(M_{\ti{X}}\!-\!M_X^0)\subset \ti\pi^{-1}\big(f_X(M_X\!-\!M_X^0)\big),\\
\hbox{and}\qquad
f_X\big|_{M_X^0}\!=\!f_{\ti{X}}\big|_{M_X^0}\!:
M_X^0\lra \M_{1,k}^0(\P,d)\subset\ov\M_{1,k}^0(\P,d),\wt\M_{1,k}^0(\P,d).
\end{gather*}
Our assumptions on $\vph$ and $f_X$ then imply that
all intersections of $f_{\ti{X}}(M_{\ti{X}})$ with $\ti\vph^{-1}(0)$
are contained in $f_{\ti{X}}(M_X^0)\!\cap\!\M_{1,k}^0(\P,d)$, are transverse,
and correspond to the intersections of $f_X(M_X)$ with~$\vph^{-1}(0)$.
Thus,
\begin{equation}\label{eulerdfn_e}\begin{split}
\blr{\tilde{\pi}^*\psi_X\cdot e(\tilde{\cal V}_{1,k}^d),
\big[\wt\M_{1,k}^0(\P,d)\big]}
& =~^{\pm}\!\big|\ti\vph^{-1}(0)\!\cap\!f_{\ti{X}}(M_{\ti{X}})\big|\\
& =~^{\pm}\!\big|\vph^{-1}(0)\!\cap\!f_X(M_X)\big|.
\end{split}\end{equation}
In particular, the left-hand side of~\e_ref{eulerdfn_e}
depends only on the homology class~$X$ used in constructing
the cohomology class~$\psi_X$ and is independent 
of the desingularization~\e_ref{desing_e}.\\

\noindent
The above argument also shows that if the cone~\e_ref{cone_e}
admits a multisection~$\vph$ satisfying~(${\cal V}1$) and~(${\cal V}2$) 
and admits a desingularization as in~\e_ref{desing_e}, then the number
$$\blr{e({\cal V}_{1,k}^d),X}
\equiv \blr{\tilde{\pi}^*\psi_X\cdot e(\tilde{\cal V}_{1,k}^d),
\big[\wt\M_{1,k}^0(\P,d)\big]}$$
is well-defined for every homology class $X$ on $\ov\M_{1,k}^0(\P,d)$.
Thus, the euler class $e({\cal V}_{1,k}^d)$ of the cone~\e_ref{cone_e}
and the sheaf~\e_ref{g1sheaf_e} is also well-defined.
In particular, the existence of homotopies satisfying~(${\cal V}1'$) and~(${\cal V}2'$)
is not absolutely necessary for showing that the euler class of~\e_ref{cone_e}
is well-defined.

\subsection{Main Result}
\label{mainres_subs}

\noindent
While the standard complex structure $J_0$ on $\P$ is ideal for many purposes,
such as computing obstruction bundles in the Gromov-Witten theory
and applying the localization theorems of~\cite{AB} and~\cite{GP1},
it is sometimes more convenient to work with an almost complex structure~$J$
on~$\P$ obtained by perturbing~$J_0$.
For this reason, we generalize Theorem~\ref{thm1} to
almost complex structures~$J$ that are close to~$J_0$.\\

\noindent
We denote by $\X_{g,k}(\P,d)$  the space of equivalence classes
of stable degree-$d$ smooth maps from genus-$g$ Riemann surfaces
with $k$~marked points to~$\P$. 
Let $\X_{g,k}^0(\P,d)$ be the subset of $\X_{g,k}(\P,d)$
consisting of stable maps with smooth domains. 
The spaces $\X_{g,k}(\P,d)$ are topologized using 
$L^p_1$-convergence on compact subsets of smooth points of the domain
and certain convergence requirements near the nodes; 
see Section~3 in~\cite{LT} for more details.
Here and throughout the rest of the paper, $p$~denotes a real number 
greater than~two.
The spaces $\X_{g,k}(\P,d)$ are stratified by 
the smooth infinite-dimensional orbifolds $\X_{\cal T}(\P)$
of stable maps from domains of the same geometric type and with the same degree
distribution between the components.
The closure of the main stratum, $\X_{g,k}^0(\P,d)$, is 
$\X_{g,k}(\P,d)$.\\

\noindent
Using modified Sobolev norms, \cite{LT} also defines a cone
$\Ga_{g,k}(T\P,d)\!\lra\!\X_{g,k}(\P,d)$ such that the fiber
of $\Ga_{g,k}(T\P,d)$ over a point $[b]\!=\![\Si,j;u]$ in 
$\X_{g,k}(\P,d)$ is the Banach space
$$\Ga_{g,k}(T\P,d)\big|_b=\Ga(b;T\P)\big/\hbox{Aut}(b),
\qquad\hbox{where}\qquad \Ga(b;T\P)=L^p_1(\Si;u^*T\P).$$
The topology on $\Ga_{g,k}(T\P,d)$ is defined similarly 
to the convergence topology on~$\X_{g,k}(\P,d)$.
If $\L$ is the line bundle $\ga^{*\otimes a}\!\lra\!\P$, let 
$\Ga_{g,k}(\L,d)\!\lra\!\X_{g,k}(\P,d)$ be the cone 
such that the fiber of $\Ga_{g,k}(\L,d)$ over $[b]\!=\![\Si,j;u]$ in 
$\X_{g,k}(\P,d)$ is the Banach space
$$\Ga_{g,k}(\L,d)\big|_b=\Ga(b;\L)\big/\hbox{Aut}(b),
\qquad\hbox{where}\qquad \Ga(b;\L)=L^p_1(\Si;u^*\L),$$
and the topology on $\Ga_{g,k}(\L,d)$ is defined analogously
to the topology on~$\Ga_{g,k}(\P,d)$.\\

\noindent
Let $\na$ denote the hermitian connection in the line bundle $\L\!\lra\!\P$ 
induced from the standard connection on the tautological line bundle over~$\P$.
If $(\Si,j)$ is a Riemann surface and $u\!:\Si\!\lra\!\P$ is a smooth map, 
let
$$\na^u\!:\Ga(\Si;u^*\L)\lra \Ga\big(\Si;T^*\Si\!\otimes\!u^*\L\big)$$
be the pull-back of $\na$ by $u$.
If $b\!=\!(\Si,j;u)$,
we define the corresponding $\bar{\partial}$-operator by
\begin{equation}\label{vdfn_e}
\bpar_{\na,b}\!:\Ga(\Si;u^*\L)\lra
\Ga\big(\Si;\La_{\I,j}^{0,1}T^*\Si\!\otimes\!u^*\L\big),
\quad
\bpar_{\na,b}\xi=\frac{1}{2}
\big(\na^u\xi+\I\na^u\xi\circ j\big),
\end{equation}
where $\I$ is the complex multiplication in the bundle $u^*\L$ and
$$\La_{\I,j}^{0,1}T^*\Si\!\otimes\!u^*\L=
\big\{\eta\!\in\!\hbox{Hom}(T\Si,u^*\L)\!:\eta\!\circ\!j=-\I\eta\big\}.$$
The kernel of $\bpar_{\na,b}$ is necessarily a finite-dimensional
complex vector space.
If $u\!:\Si\!\lra\!\P$ is a $(J_0,j)$-holomorphic map, then
$$\ker\bpar_{\na,b}=H^0\big((\Si,j);u^*\L\big)$$
is the space of holomorphic sections of the line bundle $u^*\L\!\lra\!(\Si,j)$. 
Let
$${\cal V}_{g,k}^d=\big\{[b,\xi]\!\in\!\Ga_{g,k}(\L,d)\!:
[b]\!\in\!\X_{g,k}(\P,d),~
\xi\!\in\!\ker\bpar_{\na,b}\!\subset\Ga_{g,k}(b;\L)\big\}
\subset\Ga_{g,k}(\L,d).$$
The cone ${\cal V}_{g,k}^d\!\lra\!\X_{g,k}(\P,d)$ inherits 
its topology from~$\Ga_{g,k}(\L,d)$.\\

\noindent
If $J$ is an almost complex structure on $\P$, 
let $\ov\M_{g,k}(\P,d;J)$ denote the moduli spaces of stable $J$-holomorphic 
\hbox{degree-$d$} maps from genus-$g$ Riemann surfaces with $k$ marked points
to~$\P$.  
Let
$$\M_{g,k}^0(\P,d;J)=\big\{[{\cal C},u]\!\in\!\ov\M_{g,k}(\P,d;J)\!:
{\cal C}~\hbox{is smooth}\big\}.$$
We denote by $\ov\M_{1,k}^0(\P,d;J)$ the closed subset of $\ov\M_{1,k}(\P,d;J)$ 
containing $\M_{1,k}^0(\P,d;J)$ defined in~\cite{Z5}.
If~$J$ is sufficiently close to~$J_0$, $\ov\M_{1,k}^0(\P,d;J)$ 
is the closure of $\M_{1,k}^0(\P,d;J)$ in~$\ov\M_{1,k}(\P,d;J)$.
We describe the structure of $\ov\M_{1,k}^0(\P,d;J)$ in this case 
in Lemma~\ref{g1bdstr_lmm} below.
Finally, let $\bar\Z^+$ denote the set of nonnegative integers.

\begin{thm}
\label{main_thm}
If $n,d,a\!\in\!\Z^+$ and $k\!\in\!\bar\Z^+$,
there exists $\de_n(d,a)\!\in\!\Bbb{R}^+$ with the following property.
If $J$ is an almost complex structure on~$\P$ such that 
$\|J\!-\!J_0\|_{C^1}\!<\!\de_n(d,a)$,
the moduli space $\ov\M_{1,k}^0(\P,d;J)$ carries a fundamental class
$$\big[\ov\M_{1,k}^0(\P,d;J)\big]\in H_{2(d(n+1)+k)}\big(\ov\M_{1,k}^0(\P,d;J);\Q\big).$$
Furthermore, the cone ${\cal V}_{1,k}^d\!\lra\!\X_{1,k}(\P,d)$
corresponding to the line bundle $\L\!=\!\ga^{*\otimes a}\!\lra\!\P$
determines a homology class and a cohomology class on $\ov\M_{1,k}^0(\P,d;J)$:
\begin{gather*}
\hbox{PD}_{\ov\M_{1,k}^0(\P,d;J)}\big(e({\cal V}_{1,k}^d)\big)
\in H_{2(d(n+1-a)+k)}\big(\ov\M_{1,k}^0(\P,d;J);\Q\big)\\
\hbox{and}\qquad
e({\cal V}_{1,k}^d)\in H^{2da}(\ov\M_{1,k}^0(\P,d;J);\Q\big).
\end{gather*}
Finally, if ${\cal W}\!\lra\!\X_{1,k}(\P,d)$ is a vector orbi-bundle
such that the restriction of ${\cal W}$ to each stratum $\X_{\cal T}(\P)$
of $\X_{1,k}(\P)$ is smooth, then
\begin{equation}\label{main_thm_e}
\blr{e({\cal W})\cdot e({\cal V}_{1,k}^d),\big[\ov\M_{1,k}^0(\P,d;J)\big]}
=\blr{e({\cal W})\cdot e({\cal V}_{1,k}^d),\big[\ov\M_{1,k}^0(\P,d)\big]}.
\end{equation}\\
\end{thm}

\noindent
{\it Remark:} This theorem remains valid if the compact Kahler manifold $(\P,\om_0,J_0)$,
positive integer~$d$, the holomorphic line bundle $\L\!=\!\ga^{*\otimes a}\!\lra\!\P$,
and the connection~$\na$ in $\L$ are replaced by a compact almost Kahler manifold $(X,\om,J_0)$,
a homology class $A\!\in\!H_2(X;\Z)$, and a split positive vector bundle with connection
$(\L,\na)\!\lra\!X$ such that the almost complex structure $J_0$ on $X$ 
is genus-one $A$-regular in the sense of Definition~\ref{g1comp-g1reg_dfn} in~\cite{Z5}.\\

\noindent
It is well-known that the standard complex structure is genus-one $d\ell$-regular,
where $\ell\!\in\!H_2(\P;\Z)$ is the homology class of a line.
Thus, if $J$ is an almost complex structure on $\P$ which is close to~$J_0$,
Corollary~\ref{g1comp-str_crl} and Theorem~\ref{g1comp-reg_thm} in~\cite{Z5} imply 
that $\ov\M_{1,k}^0(\P,d;J)$ is the closure of $\M_{1,k}^0(\P,d;J)$
in $\ov\M_{1,k}^0(\P,d;J)$ and is contained in a small neighborhood
of $\ov\M_{1,k}^0(\P,d)$ in~$\X_{1,k}(\P,d)$.
In addition, the stratification structure of the moduli space $\ov\M_{1,k}^0(\P,d;J)$
is the same as that of~$\ov\M_{1,k}^0(\P,d)$; 
see Lemmas~\ref{g1mainstr_lmm} and~\ref{g1bdstr_lmm} below.
Thus, $\ov\M_{1,k}^0(\P,d;J)$ carries a rational fundamental class;
see the paragraph at the end of Subsection~\ref{g1comp-str_subs} in~\cite{Z5}.\\

\noindent
The two remaining claims of Theorem~\ref{main_thm} are
the subject of Proposition~\ref{g1cone_prp}.
The restriction of the cone ${\cal V}_{1,k}^d$ to $\M_{1,k}^0(\P,d;J)$ 
is a complex vector bundle of the expected rank, i.e.~$da$.
The cone ${\cal V}_{1,k}^d$ admits a multisection $\vph$ that satisfies 
the analogues of~(${\cal V}1$) and~(${\cal V}2$) for $\ov\M_{1,k}^0(\P,d;J)$.
As in the previous subsection, the zero set of this section determines 
a homology class in real codimension $2da$ and a cohomology class 
of real dimension~$2da$.
On the other hand, if $\under{J}\!=\!(J_t)_{t\in[0,1]}$ 
is a smooth family of almost complex structures on~$\P$ such that
$J_t$ is close to $J_0$ for all $t\!\in\![0,1]$,
the moduli space
$$\ov\M_{1,k}^0(\P,d;\under{J})\equiv\!
\bigcup_{t\in[0,1]}\!\!\ov\M_{1,k}^0(\P,d;J_t)$$
is compact, by Theorem~\ref{g1comp-comp_thm} in~\cite{Z5}.
We can construct a multisection $\Phi$ of the cone ${\cal V}_{1,k}^d$ over 
$\ov\M_{1,k}^0(\P,d;\under{J})$ with properties analogous to~(${\cal V}1$) and~(${\cal V}2$).
If ${\cal W}\!\lra\!\X_{1,k}(\P,d)$ is a complex vector bundle of rank
$d(n\!+\!1\!-\!a)\!+\!k$ as in Theorem~\ref{main_thm},
we can choose a section $F$ of ${\cal W}$ over $\ov\M_{1,k}^0(\P,d;\under{J})$
such that $\Phi^{-1}(0)\!\cap\!F^{-1}(0)$ is a compact oriented one-dimensional
suborbifold of $\M_{1,k}^0(\P,A;\under{J})$ and
\begin{gather*}
\partial\big(\Phi^{-1}(0)\!\cap\!F^{-1}(0)\big)
=\Phi^{-1}(0)\!\cap\!F^{-1}(0)\!\cap\!\M_{1,k}^0(\P,d;J_1)
-\Phi^{-1}(0)\!\cap\!F^{-1}(0)\!\cap\!\M_{1,k}^0(\P,d;J_0)\\
\Lra\qquad
^{\pm}\big|\Phi^{-1}(0)\!\cap\!F^{-1}(0)\!\cap\!\M_{1,k}^0(\P,d;J_1)\big|
=~^{\pm}\big|\Phi^{-1}(0)\!\cap\!F^{-1}(0)\!\cap\!\M_{1,k}^0(\P,d)\big|.
\end{gather*}
This equality implies~\e_ref{main_thm_e}.\\

\noindent
In the next section we first summarize our detailed notation for stable maps 
and for related objects.
We then describe the structure of the moduli space $\ov\M_{1,k}^0(\P,d;J)$.
In Subsection~\ref{g1conestr_subs} we deduce Proposition~\ref{g1cone_prp}
from the descriptions of the local structure of the cone~${\cal V}_{1,k}^d$
that appear in Subsections~\ref{g1conelocalstr_subs1} and~\ref{g1conelocalstr_subs2}.
The key results of this paper, Proposition~\ref{g1conebdstr_prp} and 
Lemma~\ref{g1conebdstr_lmm}, are proved in Section~\ref{gluing_sec} by extending 
the gluing construction of Section~\ref{g1comp-strthm_sec} in~\cite{Z5}
from stable $J$-holomorphic maps to holomorphic bundle sections.

\section{Preliminaries}
\label{prelim_sec}

\subsection{Notation: Genus-Zero Maps}
\label{notation0_subs}

\noindent
We now summarize our notation for bubble maps from genus-zero Riemann surfaces
with at least one marked point,
for the spaces of such bubble maps that form
the standard stratifications of  moduli spaces of stable maps,
and for important vector bundles over them.
For more details on the notation described below, the reader is referred
to Subsections~\ref{g1comp-notation0_subs} and~\ref{g1comp-notation1_subs} in~\cite{Z5}.\\
 
\noindent
In general, moduli spaces of stable maps can stratified by the dual graph.
However, in the present situation, it is more convenient to make use
of {\it linearly ordered sets}:

\begin{dfn}
\label{index_set_dfn1}
(1) A finite nonempty partially ordered set $I$ is a {\tt linearly ordered set} if 
for all \hbox{$i_1,i_2,h\!\in\!I$} such that $i_1,i_2\!<\!h$, 
either $i_1\!\le\!i_2$ \hbox{or $i_2\!\le\!i_1$.}\\
(2) A linearly ordered set $I$ is a {\tt rooted tree} if
$I$ has a unique minimal element, 
i.e.~there exists \hbox{$\hat{0}\!\in\!I$} such that $\hat{0}\!\le\!i$ 
for {all $i\!\in\!I$}.
\end{dfn}

\noindent
If $I$ is a linearly ordered set, let $\hat{I}$ be 
the subset of the non-minimal elements of~$I$.
For every $h\!\in\!\hat{I}$,  denote by $\io_h\!\in\!I$
the largest element of $I$ which is smaller than~$h$, i.e.
$\io_h\!=\!\max\big\{i\!\in\!I:i\!<\!h\big\}$.\\

\noindent
If $M$ is a finite set,
a {\tt genus-zero $\P$-valued bubble map with $M$-marked points} is a tuple
$$b=\big(M,I;x,(j,y),u\big),$$
where $I$ is a rooted tree, and
\begin{equation}\label{stablemap_e1}
x\!:\hat{I}\!\lra\!\Bbb{C}\!=\!S^2\!-\!\{\i\},\quad  j\!:M\!\lra\!I,\quad
y\!:M\!\lra\!\Bbb{C},        \hbox{~~~and~~~} 
u\!:I\!\lra\!C^{\i}(S^2;\P)
\end{equation}
are maps such that $u_h(\i)\!=\!u_{\io_h}(x_h)$ for all $h\!\in\!\hat{I}$.
Such a tuple describes a Riemann surface $\Si_b$ and 
a continuous map $u_b\!:\Si_b\!\lra\!\P$.
The irreducible components $\Si_{b,i}$ of $\Si_b$ are indexed by the set~$I$ 
and $u_b|_{\Si_{b,i}}\!=\!u_i$.
The Riemann surface $\Si_b$ carries a special marked point,
i.e.~the point 
$$y_0(b)\!\equiv\!(\hat{0},\i)\in\Si_{b,\hat{0}}$$
if $\hat{0}$ is the minimal element of~$I$,
and $|M|$ other marked points, $(j_l,y_l)\!\in\!\Si_{b,j_l}$ with $l\!\in\!M$.\\

\noindent
The general structure of genus-zero bubble maps is described by tuples 
$${\cal T}=(M,I;j,\under{d}),$$
where $\under{d}\!:I\!\lra\!\Bbb{Z}$ is a map specifying the degree 
of $u_b|_{\Si_{b,i}}$, if $b$ is a bubble map of type~${\cal T}$.
We call such tuples {\it bubble types}.
Let ${\cal U}_{\T}(\P;J)$ denote the subset of $\ov\M_{0,\{0\}\sqcup M}(\P,d;J)$ 
consisting of stable maps $[{\cal C};u]$ such that
$$[{\cal C};u]=
\big[(\Si_b,(\hat{0},\i),(j_l,y_l)_{l\in M});u_b\big],$$
for some bubble map $b$ of type $\T$.
We recall that
$${\cal U}_{\T}(\P;J)= {\cal U}_{\T}^{(0)}(\P;J)\big/
\hbox{Aut}(\T)\!\propto\!(S^1)^I,$$
for a certain submanifold ${\cal U}_{\cal T}^{(0)}(\P;J)$ of 
the space ${\cal H}_{\T}(\P;J)$ of $J$-holomorphic maps into~$\P$
of type~${\cal T}$, not of equivalence classes of such maps;
see Subsection~2.5 in~\cite{Z3}.
For $l\!\in\!\{0\}\!\sqcup\!M$, let 
$$\ev_l\!:{\cal U}_{\T}(\P;J),{\cal U}_{\T}^{(0)}(\P;J)\lra\P$$ 
be the evaluation maps corresponding to the marked point $y_l$.

\subsection{Notation: Genus-One Maps}
\label{notation1_subs}

\noindent
We next set up analogous notation for maps from genus-one Riemann surfaces.
In this case, we also need to specify the structure of the principal component.
Thus, we index the strata of the moduli space $\ov\M_{1,M}(\P,d;J)$
by {\it enhanced linearly ordered sets}:

\begin{dfn}
\label{index_set_dfn2}
An {\tt enhanced linearly ordered set} is a pair $(I,\aleph)$,
where $I$ is a linearly ordered set, $\aleph$ is a subset of $I_0\!\times\!I_0$,
and $I_0$ is the subset of minimal elements of~$I$,
such that if $|I_0|\!>\!1$, 
$$\aleph=\big\{(i_1,i_2),(i_2,i_3),\ldots,(i_{n-1},i_n),(i_n,i_1)\big\}$$
for some bijection $i\!:\{1,\ldots,n\}\!\lra\!I_0$.
\end{dfn}

\begin{figure}
\begin{pspicture}(-1.1,-2)(10,1)
\psset{unit=.4cm}
\pscircle*(6,-3){.2}
\pscircle*(4,-1){.2}\psline[linewidth=.06]{->}(4.14,-1.14)(5.86,-2.86)
\pscircle*(8,-1){.2}\psline[linewidth=.06]{->}(7.86,-1.14)(6.14,-2.86)
\pscircle*(2,1){.2}\psline[linewidth=.06]{->}(2.14,.86)(3.86,-.86)
\pscircle*(6,1){.2}\psline[linewidth=.06]{->}(5.86,.86)(4.14,-.86)
\pscircle*(18,-3){.2}\psline[linewidth=.06](17.86,-3.14)(17.5,-3.5)
\psarc(18,-4){.71}{135}{45}\psline[linewidth=.06]{->}(18.5,-3.5)(18.14,-3.14)
\pscircle*(16,-1){.2}\psline[linewidth=.06]{->}(16.14,-1.14)(17.86,-2.86)
\pscircle*(20,-1){.2}\psline[linewidth=.06]{->}(19.86,-1.14)(18.14,-2.86)
\pscircle*(14,1){.2}\psline[linewidth=.06]{->}(14.14,.86)(15.86,-.86)
\pscircle*(18,1){.2}\psline[linewidth=.06]{->}(17.86,.86)(16.14,-.86)
\pscircle*(30,-2){.2}\pscircle*(30,-4){.2}\pscircle*(29,-3){.2}\pscircle*(31,-3){.2}
\psline[linewidth=.06]{->}(29.86,-2.14)(29.14,-2.86)
\psline[linewidth=.06]{->}(29.14,-3.14)(29.86,-3.86)
\psline[linewidth=.06]{->}(30.14,-3.86)(30.86,-3.14)
\psline[linewidth=.06]{->}(30.86,-2.86)(30.14,-2.14)
\pscircle*(27,-1){.2}\psline[linewidth=.06]{->}(27.14,-1.14)(28.86,-2.86)
\pscircle*(33,-1){.2}\psline[linewidth=.06]{->}(32.86,-1.14)(31.14,-2.86)
\pscircle*(25,1){.2}\psline[linewidth=.06]{->}(25.14,.86)(26.86,-.86)
\pscircle*(29,1){.2}\psline[linewidth=.06]{->}(28.86,.86)(27.14,-.86)
\end{pspicture}
\caption{Some Enhanced Linearly Ordered Sets}
\label{index_set_fig}
\end{figure}
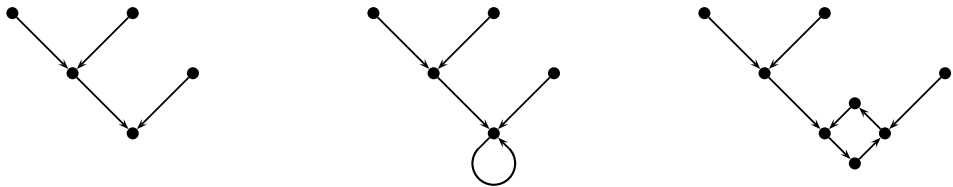

\noindent
An enhanced linearly ordered set can be represented by an oriented connected graph.
In Figure~\ref{index_set_fig}, the dots denote the elements of~$I$.
The arrows outside the loop, if there are any, 
specify the partial ordering of the linearly ordered set~$I$.
In fact, every directed edge outside of the loop
connects a non-minimal element $h$ of $I$ with~$\io_h$.
Inside of the loop, there is a directed edge from $i_1$ to $i_2$
if and only if $(i_1,i_2)\!\in\!\aleph$.\\

\noindent
The subset $\aleph$ of $I_0\!\times\!I_0$ will be used to describe
the structure of the principal curve of the domain of stable maps in 
a stratum of the moduli space~$\ov\M_{1,M}(\P,d;J)$.
If $\aleph\!=\!\eset$, and thus $|I_0|\!=\!1$,
the corresponding principal curve $\Si_{\aleph}$ 
is a smooth torus, with some complex structure.
If $\aleph\!\neq\!\eset$, the principal components form a circle of spheres:
$$\Si_{\aleph}=\Big(\bigsqcup_{i\in I_0}\{i\}\!\times\!S^2\Big)\Big/\sim,
\qquad\hbox{where}\qquad
(i_1,\i)\sim(i_2,0)~~\hbox{if}~~(i_1,i_2)\!\in\!\aleph.$$
A {\tt genus-one $\P$-valued bubble map with $M$-marked points} is a tuple
$$b=\big(M,I,\aleph;S,x,(j,y),u\big),$$
where $S$ is a smooth Riemann surface of genus one if $\aleph\!=\!\eset$
and the circle of spheres $\Si_{\aleph}$ otherwise.
The objects $x$, $j$, $y$, $u$, and $(\Si_b,u_b)$ are as in 
the genus-zero case above, except 
the sphere $\Si_{b,\hat{0}}$ is replaced by the genus-one curve $\Si_{b,\aleph}\!\equiv\!S$.
Furthermore, if $\aleph\!=\!\eset$, and thus $I_0\!=\!\{\hat{0}\}$ is a single-element set,
$u_{\hat{0}}\!\in\!C^{\i}(S;\P)$.
In the genus-one case, the general structure of bubble maps is encoded by
the tuples of the form 
$${\cal T}=(M,I,\aleph;j,\under{d}).$$
Similarly to the genus-zero case, we denote by ${\cal U}_{\cal T}(\P;J)$
the subset of $\ov\M_{1,M}(\P,d;J)$ 
consisting of stable maps $[{\cal C};u]$ such that
$$[{\cal C};u]=\big[(\Si_b,(j_l,y_l)_{l\in M});u_b\big],$$
for some bubble map $b$ of type ${\cal T}$ as above.\\

\noindent
If $\T\!=\!(M,I,\aleph;j,\under{d})$ is a bubble type as above, let
\begin{gather*}
I_1=\big\{h\!\in\!\hat{I}\!:\io_h\!\in\!I_0\big\}, \qquad
M_0=\big\{l\!\in\!M\!:j_l\!\in\!I_0\!\big\},     \qquad\hbox{and}\\
\T_0=\big(M_0\!\sqcup\!I_1,I_0,\aleph;j|_{M_0}\!\sqcup\!\io|_{I_1},
\under{d}|_{I_0}\big),
\end{gather*}
where $I_0$ is the subset of minimal elements of $I$.
For each $h\!\in\!I_1$, we put
\begin{equation*}
I_h=\big\{i\!\in\!I\!:h\!\le\!i\big\}, \qquad
M_h=\big\{l\!\in\!M\!:j_l\!\in\!I_h\!\big\},   \quad\hbox{and}\quad
{\cal T}_h=\big(M_h,I_h;j|_{M_h},\under{d}|_{I_h}\big).
\end{equation*}
The tuple $\T_0$ describes bubble maps from genus-one Riemann surfaces
with the marked points indexed by the set $M_0\!\sqcup I_1$.
By definition, we have a natural isomorphism
\begin{equation}\label{g1gendecomp_e3}\begin{split}
{\cal U}_{\cal T}(\P;J)\approx \Big(\big\{
\big(b_0,(b_h)_{h\in I_1}\big)\!\in\!{\cal U}_{{\cal T}_0}(\P;J)\!\times\!\!
\prod_{h\in I_1}\!\!{\cal U}_{{\cal T}_h}(\P;J)\!:\qquad&\\
\ev_0(b_h)\!=\!\ev_{\io_h}(b_0)~\forall h\!\in\!I_1&\big\}\Big)
\big/\hbox{Aut}^*({\cal T}),
\end{split}\end{equation}
where the group $\hbox{Aut}^*({\cal T})$ is defined by
$$\hbox{Aut}^*({\cal T})=\hbox{Aut}({\cal T})/
\{g\!\in\!\hbox{Aut}({\cal T})\!:g\cdot h\!=\!h~\forall h\!\in\!I_1\}.$$
This decomposition is illustrated in Figure~\ref{g1gendecomp_fig}.
In this figure, we represent an entire stratum of bubble maps
by the domain of the stable maps in that stratum.
We shade the components of the domain on which every (or any) stable map
in ${\cal U}_{\cal T}(\P;J)$ is nonconstant.
The right-hand side of Figure~\ref{g1gendecomp_fig} 
represents the subset of the cartesian product of the three spaces
of bubble maps, corresponding to the three drawings,
on which the appropriate evaluation maps agree pairwise,
as indicated by the dotted lines and defined in~\e_ref{g1gendecomp_e3}.\\

\begin{figure}
\begin{pspicture}(-1.1,-2)(10,1.25)
\psset{unit=.4cm}
\psellipse(8,-1.5)(1.5,2.5)
\psarc[linewidth=.05](6.2,-1.5){2}{-30}{30}\psarc[linewidth=.05](9.8,-1.5){2}{150}{210}
\pscircle[fillstyle=solid,fillcolor=gray](5.5,-1.5){1}\pscircle*(6.5,-1.5){.2}
\pscircle[fillstyle=solid,fillcolor=gray](3.5,-1.5){1}\pscircle*(4.5,-1.5){.2}
\pscircle(10.5,-1.5){1}\pscircle*(9.5,-1.5){.2}
\pscircle[fillstyle=solid,fillcolor=gray](11.91,-.09){1}\pscircle*(11.21,-.79){.2}
\pscircle[fillstyle=solid,fillcolor=gray](11.91,-2.91){1}\pscircle*(11.21,-2.21){.2}
\rput(5.5,0){$h_1$}\rput(3.5,0){$h_2$}\rput(10.3,0){$h_3$}
\rput(13.5,0.1){$h_4$}\rput(13.5,-2.9){$h_5$}
\rput(17.5,-1.5){$\approx$}
\psellipse(23,-1.5)(2.5,1.5)
\psarc[linewidth=.05](23,.3){2}{240}{300}\psarc[linewidth=.05](23,-3.3){2}{60}{120}
\psline[linewidth=.06,linestyle=dotted](28,1)(23,0)
\psline[linewidth=.06,linestyle=dotted](28,-3)(23,-3)
\pscircle*(23,0){.2}\rput(23,.7){$h_1$}
\pscircle*(23,-3){.2}\rput(23,-3.7){$h_3$}
\pscircle[fillstyle=solid,fillcolor=gray](29,1){1}\pscircle*(28,1){.2}
\pscircle[fillstyle=solid,fillcolor=gray](31,1){1}\pscircle*(30,1){.2}
\pscircle(29,-3){1}\pscircle*(28,-3){.2}
\pscircle[fillstyle=solid,fillcolor=gray](30.41,-1.59){1}\pscircle*(29.71,-2.29){.2}
\pscircle[fillstyle=solid,fillcolor=gray](30.41,-4.41){1}\pscircle*(29.71,-3.71){.2}
\rput(29,2.5){$h_1$}\rput(31,2.5){$h_2$}\rput(28.8,-1.5){$h_3$}
\rput(32,-1.4){$h_4$}\rput(32,-4.1){$h_5$}
\end{pspicture}
\caption{An Example of the Decomposition~\e_ref{g1gendecomp_e3}}
\label{g1gendecomp_fig}
\end{figure}

\noindent
Let ${\cal FT}\!\lra\!{\cal U}_{\cal T}(\P;J)$
be the bundle of gluing parameters, or of smoothings at the nodes.
This orbi-bundle has the form 
$${\cal FT}=\Big(\!\bigoplus_{(h,i)\in\aleph}\!\!\!L_{h,0}\!\otimes\!L_{i,1}\oplus
\bigoplus_{h\in\hat{I}}L_{h,0}\!\otimes\!L_{h,1}\Big)\big/\hbox{Aut}({\cal T}),$$
for certain line orbi-bundles $L_{h,0}$ and $L_{h,1}$.
Similarly to the genus-zero case,
\begin{gather}
\label{g1gendecomp_e4a}
{\cal U}_{\cal T}(\P;J)={\cal U}_{\cal T}^{(0)}(\P;J)\big/
\hbox{Aut}({\cal T})\!\propto\!(S^1)^{\hat{I}},\qquad\hbox{where}\\
\label{g1gendecomp_e4b}
{\cal U}_{\cal T}^{(0)}(\P;J)=
\big\{\big(b_0,(b_h)_{h\in I_1}\big)\!\in\!{\cal U}_{{\cal T}_0}(\P;J)\!\times\!\!
\prod_{h\in I_1}\!\!{\cal U}_{{\cal T}_h}^{(0)}(\P;J)\!:
\ev_0(b_h)\!=\!\ev_{\io_h}(b_0)~\forall h\!\in\!I_1\big\}.
\end{gather}
The line bundles $L_{h,0}$ and $L_{h,1}$ arise from the quotient~\e_ref{g1gendecomp_e4a},
and 
$${\cal FT}=\tilde{\cal F}{\cal T}\big/\hbox{Aut}({\cal T})\!\propto\!(S^1)^{\hat{I}},
\qquad\hbox{where}\qquad
\tilde{\cal F}{\cal T}=\tilde{\cal F}_{\aleph}{\cal T}\oplus
\bigoplus_{h\in\hat{I}}\tilde{\cal F}_h{\cal T},$$
$\ti{\cal F}_{\aleph}\T\!\lra\!{\cal U}_{\cal T}^{(0)}(\P;J)$
is the bundle of smoothings for the $|\aleph|$ nodes of the circle of spheres $\Si_{\aleph}$
and $\ti{\cal F}_h\T\!\lra\!{\cal U}_{\T}^{(0)}(\P;J)$ is the line bundle
of smoothings of the attaching node of the bubble indexed by~$h$.
We denote by ${\cal FT}^{\eset}$ and $\tilde{\cal F}{\cal T}^{\eset}$ the subsets of 
${\cal FT}$ and $\tilde{\cal F}{\cal T}$, respectively,
consisting of the elements with all components nonzero.\\

\noindent
Suppose ${\cal T}\!=\!(M,I,\aleph;j,\under{d})$ is a bubble type such that
$d_i\!=\!0$ for all $i\!\in\!I_0$, 
i.e.~every element in ${\cal U}_{\cal T}(\P;J)$ is constant on the principal components.
In this case, the decomposition \e_ref{g1gendecomp_e3} is equivalent to
\begin{equation}\label{g1decomp_e1}\begin{split}
{\cal U}_{\cal T}(\P;J)&\approx
\Big({\cal U}_{{\cal T}_0}(pt)\times
{\cal U}_{\bar{\cal T}}(\P;J)\Big)\big/\hbox{Aut}^*({\cal T})\\
&\subset\Big(\ov\cM_{1,k_0}\times
{\cal U}_{\bar{\cal T}}(\P;J)\Big)\big/\hbox{Aut}^*({\cal T}),
\end{split}\end{equation}
where $k_0\!=\!|I_1|\!+\!|M_0|$,  $\ov\cM_{1,k_0}$ is the moduli space of 
genus-one curves with $k_0$ marked points, and
$${\cal U}_{\bar{\cal T}}(\P;J)=
\big\{(b_h)_{h\in I_1}\!\in\!\prod_{h\in I_1}\!{\cal U}_{{\cal T}_h}(\P;J)\!:
\ev_0(b_{h_1})\!=\!\ev_0(b_{h_2})~\forall h_1,h_2\!\in\!I_1\big\}.$$
Similarly, \e_ref{g1gendecomp_e4a} and~\e_ref{g1gendecomp_e4b} are equivalent~to
\begin{gather}\label{g1decomp_e3a}
{\cal U}_{\cal T}^{(0)}(\P;J)\approx
{\cal U}_{{\cal T}_0}(pt)\times {\cal U}_{\bar{\cal T}}^{(0)}(\P;J)
\subset\ov\cM_{1,k_0}\times{\cal U}_{\bar{\cal T}}^{(0)}(\P;J),
\qquad\hbox{where}\\
\label{g1decomp_e3b}
{\cal U}_{\bar{\cal T}}^{(0)}(\P;J)=
\big\{(b_h)_{h\in I_1}\!\in\!\prod_{h\in I_1}\!{\cal U}_{{\cal T}_h}^{(0)}(\P;J)\!:
\ev_0(b_{h_1})\!=\!\ev_0(b_{h_2})~\forall h_1,h_2\!\in\!I_1\big\}.
\end{gather}
We denote by
$$\pi_P\!:{\cal U}_{\T}(\P;J),{\cal U}_{\T}^{(0)}(\P;J)\lra\ov\cM_{1,k_0}$$
the projections onto the first component in 
the decompositions~\e_ref{g1decomp_e1}  and~\e_ref{g1decomp_e3a}.
Let
$$\ev_P\!:{\cal U}_{\T}(\P;J),{\cal U}_{\T}^{(0)}(\P;J)\lra\P$$
be the maps sending every element $b\!=\!(\Si_b,u_b)$ of ${\cal U}_{\T}(\P;J)$
and ${\cal U}_{\T}^{(0)}(\P;J)$
to the image of the principal component $\Si_{b;P}$  of $\Si_b$ under~$u_b$.\\

\noindent
If $\T\!=\!(M,I,\aleph;j,\under{d})$ is as in the previous paragraph, let
$$\chi(\T)=\big\{i\!\in\!\hat{I}\!:d_i\!\neq\!0;~d_h\!=\!0~\forall h\!<\!i\big\}.$$
The subset $\chi(\T)$ of $I$ indexes the
first-level effective bubbles of every element of~${\cal U}_{\T}^{(0)}(\P;J)$.
For each element $b\!=\!(\Si_b,u_b)$ of ${\cal U}_{\T}^{(0)}(\P;J)$
and $i\!\in\!\chi(\T)$, let
$$\cD_ib=\big\{du_b|_{\Si_{b,i}}\big\}\big|_{\i}e_{\i}\in T_{\ev_P(b)}\P,
\qquad\hbox{where}\qquad e_{\i}=(1,0,0)\in T_{\i}S^2.$$
In geometric terms, the complex span of $\cD_ib$ in $T_{\ev_P(b)}\P$
is the line tangent to the rational component~$\Si_{b,i}$
at the node of~$\Si_{b,i}$ closest to a principal component of~$\Si_b$.
If the branch corresponding to~$\Si_{b,i}$ has a cusp at this node,
then $\cD_ib\!=\!0$.
Let
$$\ti\F\T=\!\bigoplus_{i\in\chi(\T)}\!\!\ti{\cal F}_{h(i)}\T
\lra{\cal U}_{\T}^{(0)}(\P;J), \quad\hbox{where}\quad
h(i)\!=\!\min\{h\!\in\!\hat{I}\!:h\!\le\!i\}\in I_1.$$
We define the bundle map
$$\rho\!: \ti{\cal F}\T\lra\ti\F\T$$
over ${\cal U}_{\T}^{(0)}(\P;J)$ by 
\begin{gather*}
\rho(\ups)\!=\!\big(b;(\rho_i(v))_{i\in\chi(\T)}\big)\in\ti\F\T,
\quad\hbox{where}\quad
\rho_i(v)\!=\!\!\prod_{h\in\hat{I},h\le i}\!\!\!\!\!\!v_h\in\ti{\cal F}_{h(i)}\T,
\qquad\hbox{if}\\
\ups\!=\!\big(b;v_{\aleph},(v_h)_{h\in\hat{I}}\big),~~~
b\!\in\!{\cal U}_{\T}^{(0)}(\P;J),~~~
v_{\aleph}\!\in\!\ti{\cal F}_{\aleph}\T|_b,~~~
v_h\!\in\!\ti{\cal F}_h\T|_b \approx 
\begin{cases}
T_{x_h(b)}\Si_{b;P},& \hbox{if}~h\!\in\!I_1,\\
\C,& \hbox{if}~ h\!\in\!\hat{I}\!-\!I_1,
\end{cases}
\end{gather*}
where $x_h(b)\!\in\!\Si_{b;P}$ is the node joining the bubble $\Si_{b,h}$
of $b$ to the principal component $\Si_{b;P}$ of~$\Si_b$.
This definition is illustrated in Figure~\ref{g1bdstr_fig} on page~\pageref{g1bdstr_fig}.\\

\noindent
Let $\E\!\lra\!\ov\cM_{1,k_0}$ be the Hodge line bundle,
i.e.~the line bundle of holomorphic differentials.
For each $i\!\in\!\chi({\cal T})$, we define the bundle map 
$$\cD_{J,i}\!:\ti{\cal F}_{h(i)}\T \lra \pi_P^*\E^*\!\otimes_J\!\ev_P^*T\P$$
over ${\cal U}_{\T}^{(0)}(\P;J)$ by
$$\big\{\cD_{J,i}(b,w_i)\big\}(\psi)=
\psi_{x_{h(i)}(b)}(w_i)\cdot_J\cD_ib \in T_{\ev_P(b)}\P
\quad\hbox{if}\quad b\!\in\!{\cal U}_{\T}^{(0)}(\P;J),~w_i\!\in\!\ti{\cal F}_{h(i)}\T|_b,
~\psi\!\in\!\pi_P^*\E|_b.$$
Let
$$\cD_{\T}\!: \ti\F\T \lra \pi_P^*\E^*\!\otimes\!\ev_P^*T\P$$
be the bundle map over ${\cal U}_{\T}^{(0)}(\P;J)$ given~by
$$\cD_{\T}\big(b,(w_i)_{i\in\chi(\T)}\big)
  =\sum_{i\in\chi(\T)}\!\!\!\cD_{J,i}(b,w_i).$$
It descends to a bundle map
$$\cD_{\T}\!: \F\T\lra \pi_P^*\E^*\!\otimes\!\ev_P^*T\P
\big/\hbox{Aut}^*(\T)$$
over ${\cal U}_{\T}(\P;J)$, for a bundle $\F\T\!\lra\!{\cal U}_{\T}(\P;J)$.\\

\noindent
Let $\ti{\cal V}_{1,k}^d\!\lra\!{\cal U}_{\T}^{(0)}(\P;J)$ be 
the cone such~that the fiber of $\tilde{\cal V}_{1,k}^d$
over a point $b\!=\!(\Si_b,u_b)$ in ${\cal U}_{\T}^{(0)}(\P;J)$
is $\ker\bpar_{\na,b}$; see Subsection~\ref{mainres_subs}.
If $b\!=\!(\Si_b,u_b)\!\in\!{\cal U}_{\T}^{(0)}(\P;J)$,
$\xi\!=\!(\xi_h)_{h\in I}\!\in\!\Ga(b;\L)$, and $i\!\in\!\chi(\T)$, let
$$\D_{\T,i}\xi=\na_{e_{\i}}\xi_i\in\L_{\ev_P(b)}.$$
The element $\na_{e_{\i}}\xi_i$ of $u_{b,i}^*\L|_{\i}$ is
the covariant derivative of the section $\xi_i\!\in\!\Ga(\Si_{b,i};u_{b,i}^*\L)$
at~$\i\!\in\!\Si_{b,i}$ with respect to the connection~$\na$ 
in $\L$ along~$e_{\i}$; see Subsection~\ref{mainres_subs}.
Note that if $\xi\!\in\!\ker\bpar_{\na,b}$, then
\begin{equation}\label{diffdfn_e}
\na_{c\cdot e_{\i}}\xi_i=c\cdot\D_{\T,i}\xi 
\qquad\forall\, c\!\in\!\C.
\end{equation}
We next define the bundle~map
$$\D_{\T}\!:  \ti\F\T \lra
\Hom\big(\ti{\cal V}_{1,k}^d,\pi_P^*\E^*\!\otimes\!\ev_P^*\L\big)$$
over ${\cal U}_{\T}^{(0)}(\P;J)$ by
\begin{gather*}
\big\{\D_{\T}(b,\xi\!\otimes\!w)\big\}(\psi)
=\sum_{i\in\chi(\T)}\!\psi_{x_{h(i)}(b)}(w_i)\cdot 
\D_{\T,i}\xi \in\L_{\ev_P(b)} \qquad\hbox{if}\\
\xi\in\ti{\cal V}_{1,k}^d|_b\subset\Ga(b;\L), \qquad
w\!=\!(w_i)_{i\in\chi(\T)}\in\ti\F\T|_b, \quad\hbox{and}\quad \psi\in\E_{\pi_P(b)}.
\end{gather*}
By \e_ref{diffdfn_e}, the bundle map $\D_{\T}$  induces a linear bundle~map 
$$\F\T\lra \Hom({\cal V}_{1,k}^d,\pi_P^*\E^*\otimes\ev_P^*\L\big/\hbox{Aut}^*(\T)\big)$$
over~${\cal U}_{\T}(\P;J)$.\\

\noindent
Finally, all vector orbi-bundles we encounter will be assumed to be normed.
Some will come with natural norms; for others, we implicitly choose a norm
once and for~all.
If \hbox{$\pi_{\F}\!:\F\!\lra\!\X$} is a normed vector bundle
and $\de\!:\X\!\lra\!\R$ is any function, possibly constant,
let
$$\F_{\de}=\big\{\ups\!\in\!\F\!: |\ups|\!<\!\de(\pi_{\F}(\ups))\big\}.$$
If $\Om$ is any subset of $\F$, we take  $\Om_{\de}\!=\!\Om\cap\F_{\de}$.

\subsection{The Structure of the Moduli Space $\ov\M_{1,k}^0(\P,d;J)$}
\label{g1str_subs}

\noindent
We now describe the structure of the moduli space $\ov\M_{1,k}^0(\P,d;J)$
near each of its strata.
The first part of Theorem~\ref{main_thm} follows
from the first claims of Lemmas~\ref{g1mainstr_lmm} and~\ref{g1bdstr_lmm} below.
If $k\!\in\!\Z$, we denote by $[k]$ the set of positive integers
that do not exceed~$k$.

\begin{lmm}
\label{g1mainstr_lmm}
If $n$, $k$, and $d$ are as in Theorem~\ref{main_thm},
there exists $\de_n(d)\!\in\!\Bbb{R}^+$ with the following property.
If $J$ is an almost complex structure on~$\P$, 
such that $\|J\!-\!J_0\|_{C^1}\!<\!\de_n(d)$,
and 
$${\cal T}=([k],I,\aleph;j,\under{d})$$
 is a bubble type such that
$\sum_{i\in I}\!d_i\!=\!d$ and $d_i\!\neq\!0$ for some minimal element $i$ of $I$,
then ${\cal U}_{\cal T}(\P;J)$ is a smooth orbifold,
$$\dim{\cal U}_{\cal T}(\P;J)=2\big(d(n\!+\!1)+k-|\aleph|-|\hat{I}|\big),
\qquad\hbox{and}\qquad
{\cal U}_{\cal T}(\P;J)\subset\ov\M_{1,k}^0(\P,d;J).$$
Furthermore, there exist $\de\!\in\!C({\cal U}_{\cal T}(\P;J);\Bbb{R}^+)$,
an open neighborhood $U_{\cal T}$ of ${\cal U}_{\cal T}(\P;J)$ in $\X_{1,k}(\P,d)$, 
and an orientation-preserving homeomorphism
$$\phi_{\cal T}\!:{\cal FT}_{\de}\lra \ov\M_{1,k}^0(\P,d;J)\cap U_{\cal T},$$
which restricts to a diffeomorphism 
${\cal FT}_{\de}^{\eset}\!\lra\!\M_{1,k}^0(\P,d;J)\cap U_{\cal T}$.
\end{lmm}

\noindent
By Theorem~\ref{g1comp-reg_thm} in~\cite{Z5}, there exists 
$\de_n(d)\!\in\!\Bbb{R}^+$ with the following property.
If $J$ is an almost complex structure on~$\P$, 
such that $\|J\!-\!J_0\|_{C^1}\!<\!\de_n(d)$,
$\Si$ is a genus-one prestable Riemann surface, and
$u\!:\Si\!\lra\!\P$ is a $J$-holomorphic map, 
such that the restriction of~$u$ to the principal component(s) of~$\Si$ is not constant,
then the linearization $D_{J,u}$ of the $\bar{\partial}_J$-operator at~$u$ is surjective.
From standard arguments, such as in Chapter~3 of~\cite{MS}, 
it then follows that the stratum ${\cal U}_{\cal T}(\P;J)$ of $\ov\M_{1,k}^0(\P,d;J)$,
where ${\cal T}$ is a bubble type as in Lemma~\ref{g1mainstr_lmm},
is a smooth orbifold of the expected dimension.
Furthermore, there is no obstruction to gluing the maps in~${\cal U}_{\cal T}(\P;J)$,
in the sense of the following paragraph.\\

\noindent
We fix a metric $g_n$ and a connection $\na^n$ on $(T\P,J)$. 
For each sufficiently small element $\ups\!=\!(b,v)$ of $\tilde{\cal F}{\cal T}^{\eset}$
and $b\!=\!(\Si_b,u_b)\!\in\!{\cal U}_{\cal T}^{(0)}(\P;J)$,
let 
$$q_{\ups}\!: \Si_{\ups}\lra\Si_b$$
be the basic gluing map constructed in Subsection~\ref{g1comp-reg1_subs1} of~\cite{Z5}.
In this case, $\Si_{\ups}$ is a smooth elliptic curve.
Let
$$b(\ups)=\big(\Si_{\ups},j_{\ups},u_{\ups}\big),
\qquad\hbox{where}\qquad u_{\ups}=u_b\circ q_{\ups},$$
be the corresponding approximately $J$-holomorphic stable map.
By the previous paragraph,
the linearization~$D_{J,b}$ of the $\bar{\partial}_J$-operator at~$b$ is surjective.
Thus, if $\ups$ is sufficiently small, the linearization
$$D_{J,\ups}\!:\Ga(\ups)\!\equiv\!L^p_1(\Si_{\ups};u_{\ups}^*T\P)\lra
\Ga^{0,1}(\ups)\!\equiv\!
L^p(\Si_{\ups};\La^{0,1}_{J,j}T^*\Si_{\ups}\!\otimes\!u_{\ups}^*T\P),$$
of the $\bar{\partial}_J$-operator at~$b(\ups)$, defined via~$\na^n$,
is also surjective.
In particular, we can obtain an orthogonal decomposition
\begin{equation}\label{gadecomp_e1}
\Ga(\ups)=\Ga_-(\ups)\oplus\Ga_+(\ups)
\end{equation}
such that the linear operator 
$D_{J,\ups}\!:\Ga_+(\ups)\!\lra\!\Ga^{0,1}(\ups)$ is an isomorphism,
while $\Ga_-(\ups)$ is close to~$\ker D_{J,b}$.
The $L^2$-inner product on $\Ga(\ups)$ used in the orthogonal decomposition
is defined via the metric~$g_n$ on~$\P$ and 
the metric~$g_{\ups}$ on~$\Si_{\ups}$ induced by the pregluing construction.
The Banach spaces $\Ga(\ups)$ and $\Ga^{0,1}(\ups)$ carry the norms 
$\|\cdot\|_{\ups,p,1}$ and $\|\cdot\|_{\ups,p}$, respectively,
which are also defined by the pregluing construction.
These norms are equivalent to the ones used in Section~3 of~\cite{LT}.
In particular, the norms of $D_{J,\ups}$ and of the inverse of its restriction
to $\Ga_+(\ups)$ have fiberwise uniform upper bounds, 
i.e.~dependent only on $[b]\!\in\!{\cal U}_{\T}(\P;J)$, and
not on $\ups\!\in\!\ti{\cal F}\T^{\eset}$.
It then follows that the equation 
$$\bar{\partial}_J\exp_{u_{\ups}}\!\ze=0  \qquad\Llra\qquad
[\Si_{\ups},\exp_{u_{\ups}}\!\ze]\in\M_{1,k}^0(\P,d;J)$$
has  a unique small solution $\ze_{\ups}\!\in\!\Ga_+(\ups)$.
Furthermore, 
$$\|\ze_{\ups}\|_{\ups,p,1}\le C(b)|\ups|^{1/p},$$
for some $C\!\in\!C({\cal U}_{\T}(\P;J);\R^+)$.
The diffeomorphism on ${\cal FT}_{\de}^{\eset}$ is given~by
$$\phi_{\cal T}\!:{\cal FT}_{\de}^{\eset}\lra\M_{1,k}^0(\P,d;J),
\qquad \phi_{\cal T}([\ups])=\big[\ti{b}(\ups)\big],
\qquad\hbox{where}\qquad
\tilde{b}(\ups)=(\Si_{\ups},\exp_{u_{\ups}}\!\ze_{\ups});$$
see the paragraph following Lemma~\ref{g1comp-reg0_lmm1} in~\cite{Z5}.
This map extends to a homeomorphism
$$\phi_{\cal T}\!:{\cal FT}_{\de}\lra\ov\M_{1,k}^0(\P,d;J),$$
as can be seen by an argument similar to Subsections~3.9 and~4.1 in~\cite{Z3}.\\

\noindent
We denote by $\M_{1,k}^{\{0\}}(\P,d;J)$ the union of the strata 
$\U_{\T}(\P;J)$ with $\T$ as in Lemma~\ref{g1mainstr_lmm}. In other words,
$$\M_{1,k}^{\{0\}}(\P,d;J)= \big\{ [{\cal C},u]\!\in\!\ov\M_{1,k}(\P,d;J)\!:
u|_{\cC_P}~\hbox{is not constant}\big\},$$
where $\cC_P$ is the principal component of the domain $\cC$ of~$u$.

\begin{lmm}
\label{g1bdstr_lmm}
If $n$, $k$, and $d$ are as in Theorem~\ref{main_thm},
there exists $\de_n(d)\!\in\!\R^+$ with the following property.
If $J$ is an almost complex structure on~$\P$, such that $\|J\!-\!J_0\|_{C^1}\!<\!\de_n(d)$,
and 
$${\cal T}=([k],I,\aleph;j,\under{d})$$
is a bubble type such that
$\sum_{i\in I}\!d_i\!=\!d$ and $d_i\!=\!0$ for all minimal elements $i$ of~$I$, then
${\cal U}_{\cal T}(\P;J)$ is a smooth orbifold,
\begin{gather*}
\dim{\cal U}_{\cal T}(\P;J)=2\big(d(n\!+\!1)\!+\!k-|\aleph|\!-\!|\hat{I}|+n\big),
~~\hbox{and}~~
\ov\M_{1,k}^0(\P,d;J)\cap{\cal U}_{\cal T}(\P;J)=
{\cal U}_{{\cal T};1}(\P;J),\\
\hbox{where}\qquad
{\cal U}_{{\cal T};1}(\P;J)=
\big\{[b]\!\in\!{\cal U}_{\cal T}(\P;J)\!:\dim_{\Bbb{C}}\hbox{Span}_{(\Bbb{C},J)}
\{{\cal D}_ib\!:i\!\in\!\chi({\cal T})\}\!<\!|\chi({\cal T})|\big\}.
\end{gather*}
The space ${\cal U}_{{\cal T};1}(\P;J)$ admits a stratification
by smooth suborbifolds of ${\cal U}_{\cal T}(\P;J)$:
\begin{gather*}
{\cal U}_{\T;1}(\P;J)= \bigsqcup_{m=\max(|\chi(\T)|-n,1)}^{m=|\chi(\T)|}
\!\!\!\!\!\!\!{\cal U}_{\T;1}^m(\P;J)
\qquad\hbox{such that}\\
\begin{split}
\dim{\cal U}_{\T;1}^m(\P;J)
&= 2\big(d(n\!+\!1)\!+\!k-|\aleph|\!-\!|\hat{I}| +n+(|\chi(\T)|\!-\!n\!-\!m)m\big)\\
&\le \dim\M_{1,k}^0(\P,d;J)-2.
\end{split}
\end{gather*}
Furthermore, the space
$${\cal F}^1\T^{\eset}\equiv \big\{[b,\ups]\!\in\!{\cal FT}^{\eset}\!:
\cD_{\T}\big(\rho(\ups)\big)\!=\!0\big\}$$
is a smooth oriented suborbifold of ${\cal FT}$.
Finally, there exist $\de\!\in\!C({\cal U}_{\T}(\P;J);\R^+)$,
an open neighborhood $U_{\T}$ of ${\cal U}_{\T}(\P;J)$ in $\X_{1,k}(\P,d)$, 
and an orientation-preserving diffeomorphism
$$\phi_{\T}\!:
{\cal F}^1{\T}^{\eset}_{\de}\lra \M_{1,k}^0(\P,d;J)\cap U_{\T},$$
which extends to a homeomorphism 
$$\phi_{\T}\!:{\cal F}^1{\T}_{\de}\lra\ov\M_{1,k}^0(\P,d;J)\cap U_{\T},$$
where ${\cal F}^1{\T}$ is the closure of ${\cal F}^1\T^{\eset}$ in ${\cal FT}$.
\end{lmm}

\begin{figure}
\begin{pspicture}(-1.1,-1.8)(10,1.25)
\psset{unit=.4cm}
\psellipse(5,-1.5)(1.5,2.5)
\psarc[linewidth=.05](3.2,-1.5){2}{-30}{30}\psarc[linewidth=.05](6.8,-1.5){2}{150}{210}
\pscircle[fillstyle=solid,fillcolor=gray](2.5,-1.5){1}\pscircle*(3.5,-1.5){.2}
\pscircle[fillstyle=solid,fillcolor=gray](.5,-1.5){1}\pscircle*(1.5,-1.5){.2}
\pscircle(7.5,-1.5){1}\pscircle*(6.5,-1.5){.2}
\pscircle[fillstyle=solid,fillcolor=gray](8.91,-.09){1}\pscircle*(8.21,-.79){.2}
\pscircle[fillstyle=solid,fillcolor=gray](8.91,-2.91){1}\pscircle*(8.21,-2.21){.2}
\rput(2.5,0){$h_1$}\rput(.5,0){$h_2$}\rput(7.3,0){$h_3$}
\rput(10.5,0.1){$h_4$}\rput(10.5,-2.9){$h_5$}
\rput(5,-5){\small ``tacnode"}
\pnode(5,-5){A1}\pnode(3.5,-1.5){B1}
\ncarc[nodesep=.35,arcangleA=-25,arcangleB=-15,ncurv=1]{->}{A1}{B1}
\pnode(5,-4.65){A2}\pnode(7.3,-1.5){B2}
\ncarc[nodesep=0,arcangleA=40,arcangleB=30,ncurv=1]{-}{A2}{B2}
\pnode(8,-.95){B2a}\pnode(8.02,-2.02){B2b}
\ncarc[nodesep=0,arcangleA=0,arcangleB=10,ncurv=1]{->}{B2}{B2a}
\ncarc[nodesep=0,arcangleA=0,arcangleB=10,ncurv=1]{->}{B2}{B2b}
\rput(25,-1.5){\begin{small}\begin{tabular}{l}
$\chi(\T)\!=\!\{h_1,h_4,h_5\}$,~
$\rho(\ups)\!=\!(\ups_{h_1},\ups_{h_3}\ups_{h_4},\ups_{h_3}\ups_{h_5})$\\
\\
${\cal F}^1\T^{\eset}=\big\{[b;v_{h_1},v_{h_2},v_{h_3},v_{h_4},v_{h_5}]\!:
v_{h_2},v_{h_4},v_{h_5}\!\in\!\C^*$;\\
${}\qquad\qquad\qquad v_{h_1}\!\in\!T_{x_{h_1}}\Si_{\aleph}\!-\!\{0\},~
v_{h_3}\!\in\!T_{x_{h_3}}\Si_{\aleph}\!-\!\{0\}$;\\
${}\qquad\qquad\qquad v_{h_1}\cD_{J,h_1}b\!+\!v_{h_3}v_{h_4}\cD_{J,h_4}b\!+\!
v_{h_3}v_{h_5}\cD_{J,h_5}b\!=\!0\big\}$
\end{tabular}\end{small}}
\end{pspicture}
\caption{An Illustration of Lemma~\ref{g1bdstr_lmm}} 
\label{g1bdstr_fig}
\end{figure}

\noindent
We now clarify the statement of Lemma~\ref{g1bdstr_lmm}
and illustrate it using Figure~\ref{g1bdstr_fig}.
As before, the shaded discs represent the components of the domain
on which every stable map~$[b]$ in ${\cal U}_{\cal T}(\P;J)$ is non-constant.
The element $[\Si_b,u_b]$ of ${\cal U}_{\cal T}(\P;J)$ is 
in the stable-map closure of $\M_{1,k}^0(\P,d;J)$ if and only if
the branches of $u_b(\Si_b)$ corresponding to the attaching nodes of 
the first-level effective bubbles of $[\Si_b,u_b]$ form a generalized tacnode.
In the case of Figure~\ref{g1bdstr_fig}, this means that  either\\
${}\quad$ (a) for some $h\!\in\!\{h_1,h_4,h_5\}$, 
the branch of $u_b|_{\Si_{b,h}}$ at the node~$\i$ has a cusp, or\\
${}\quad$ (b) for all $h\!\in\!\{h_1,h_4,h_5\}$,
the branch of $u_b|_{\Si_{b,h}}$ at the node~$\i$ is smooth, but the dimension\\
${}\qquad~~~$ of the span of the three lines tangent to these
branches is less than three.\\

\noindent
The last statement of Lemma~\ref{g1bdstr_lmm} identifies 
a normal neighborhood of ${\cal U}_{{\cal T};1}(\P;J)$ 
in $\ov\M_{1,k}^0(\P,d;J)$ with a small neighborhood of
${\cal U}_{{\cal T};1}(\P;J)$ in the bundle ${\cal F}^1{\cal T}$ 
over~${\cal U}_{{\cal T};1}(\P;J)$.
Each fiber of the projection map ${\cal F}^1{\cal T}\!\lra\!{\cal U}_{{\cal T};1}(\P;J)$ 
is an algebraic variety. See Figure~\ref{g1bdstr_fig} for an example.\\

\noindent
The first three claims of Lemma~\ref{g1bdstr_lmm} follow immediately from 
Theorems~\ref{g1comp-reg_thm} and~\ref{g1comp-str_thm} in~\cite{Z5}
and the decomposition~\e_ref{g1gendecomp_e3}.
The last two statements of Lemma~\ref{g1bdstr_lmm} are a special case of 
the last two statements of the latter theorem.\\

\noindent
If $\T$ is a bubble type as in Lemma~\ref{g1bdstr_lmm} 
and $m$ is a positive integer, let 
$${\cal U}_{\T;1}^m(\P;J)=\big\{
[b]\!\in\!{\cal U}_{\T}(\P;J)\!:\dim_{\C}\hbox{Span}_{(\C,J)}
\{\cD_ib\!:i\!\in\!\chi(\T)\}=|\chi(\T)|\!-\!m\big\}
\subset{\cal U}_{\T;1}(\P;J).$$
By definition, the subspaces ${\cal U}_{\T;1}^m(\P;J)$ of 
${\cal U}_{\T}(\P;J)$ partition~${\cal U}_{\T;1}(\P;J)$.
On the other hand, 
$${\cal U}_{\T;1}^m(\P;J)\neq\eset \qquad\Lra\qquad
\max\big(|\chi(\T)|\!-\!n,1\big)\le m\le |\chi(\T)|.$$
In order to show that the space ${\cal U}_{\T;1}^m(\P;J)$ is
a smooth suborbifold of ${\cal U}_{\T}(\P;J)$ of the claimed dimension,
below we describe ${\cal U}_{\T;1}^m(\P;J)$ in a different way.\\

\noindent 
For each $i\!\in\!\hat{I}$, let
$$M_i=\big\{l\!\in\!M\!:j_l\!=\!i\}\!\sqcup\!\big\{h\!\in\!\hat{I}\!:\io_h\!=\!i\big\}.$$
We denote by
$$\pi_i\!: {\cal U}_{\T}(\P;J) \lra \M_{0,\{0\}\sqcup M_i}^0(\P,d_i;J)$$
the map sending each bubble map $(\Si_b,u_b)$ to its restriction to the component
$\Si_{b,i}\!\subset\Si$.
Let
$$L_0\lra  \M_{0,\{0\}\sqcup M_i}^0(\P,d_i;J)\subset\ov\M_{0,\{0\}\sqcup M_i}(\P,d_i;J)$$
be the universal tangent line bundle for the special point labeled by~$0$,
i.e.~$(i,\i)$ in the notation of Subsection~\ref{notation0_subs}.
We put 
$${\cal F}=\!\bigoplus_{i\in\chi(\T)}\!\!\!\pi_i^*L_0\lra{\cal U}_{\T}(\P;J).$$
While each line bundle $\pi_i^*L_0$ may not be well-defined,
the orbibundle~${\cal F}$ is always well-defined.
We denote~by 
$$\pi_m\!: \Gr_m{\cal F} \lra {\cal U}_{\T}(\P;J) \qquad\hbox{and}\qquad 
\ga_m\lra\hbox{Gr}_m{\cal F}$$
the Grasmannian bundle of $m$-dimensional linear subspaces and
the tautological $m$-plane bundle, respectively.
Let 
\begin{gather*}
{\cal S}_m={\cal D}_m^{-1}(0)\subset\Gr_m{\cal F},\qquad\hbox{where}\\
{\cal D}_m\in\Ga(\Gr_m{\cal F};\ga_m^*\!\otimes\!\pi_m^*\ev_P^*T\P),\qquad
{\cal D}_m([\ups])=\!\sum_{i\in\chi({\cal T})}\!\!\!\!{\cal D}_{J,i}\ups_i\in\ev_P^*T\P
~~~\hbox{if}~~~[\ups]=\big[(\ups_i)_{i\in\chi(\T)}\big].
\end{gather*}
By Theorem~\ref{g1comp-reg_thm} in~\cite{Z5}, the section~${\cal D}_m$
is transverse to the zero set if $\de_n(d)$ is sufficiently small.
Thus, ${\cal S}_m$ is a smooth suborbifold of~$\Gr_m{\cal F}$ of
dimension
\begin{equation*}\begin{split}
\dim{\cal S}_m
&=\dim\Gr_m{\cal F}-2\,\rk\!\!\ga_m^*\!\otimes\!\pi_m^*\ev_P^*T\P\\
&=2\big(d(n\!+\!1)\!+\!k\!-\!|\aleph|\!-\!|\hat{I}|\!+\!n\big)+
2m\big(|\chi({\cal T})|\!-\!m\big)-2\cdot n\cdot m\\
&=2\big(d(n\!+\!1)\!+\!k\!-\!|\aleph|\!-\!|\hat{I}|\!+\!n+
m\big(|\chi({\cal T})|\!-\!n\!-\!m\big)\big).
\end{split}\end{equation*}
The image of ${\cal S}_m$ under the bundle projection map $\pi_m$
is the union of the spaces ${\cal U}_{{\cal T};1}^{m'}(\P;J)$ with $m'\!\ge\!m$.
The map $\pi_m|_{{\cal S}_m}$ is an immersion at $[\ups]\!\in\!{\cal S}_m$ 
if $\pi_m^{-1}(\pi_m([\ups]))\!=\![\ups]$.
The latter is the case if and only $\pi_m([\ups])\!\in\!{\cal U}_{{\cal T};1}^m(\P;J)$.
Thus, the subspace ${\cal U}_{{\cal T};1}^m(\P;J)$ is a smooth suborbifold of
${\cal U}_{\T}(\P;J)$ of~$\dim{\cal S}_m$.

\section{Proof of Theorem~\ref{main_thm}}
\label{main_sec}

\subsection{The Global Structure of the Cone ${\cal V}_{1,k}^d\!\lra\!\ov\M_{1,k}^0(\P,d;J)$}
\label{g1conestr_subs}

\noindent 
In this section, we prove Proposition~\ref{g1cone_prp},
which contains the last two statements of Theorem~\ref{main_thm}.
The key nontrivial ingredient in the proof of Proposition~\ref{g1cone_prp}
is Proposition~\ref{g1conebdstr_prp}, which is proved in Section~\ref{gluing_sec}.

\begin{prp}
\label{g1cone_prp}
If $n$, $k$, $d$, $a$, $\L$, and ${\cal V}_{1,k}^d$ are as in Theorem~\ref{main_thm},
there exists $\de_n(d,a)\!\in\!\Bbb{R}^+$  with the following property.
If $J$ is an almost complex structure on~$\P$, 
such that $\|J\!-\!J_0\|_{C^1}\!<\!\de_n(d,a)$,
the requirements of Lemmas~\ref{g1mainstr_lmm} and~\ref{g1bdstr_lmm} are satisfied.
Furthermore, ${\cal V}_{1,k}^d\!\lra\!\M_{1,k}^0(\P,d;J)$ 
is a smooth complex vector orbibundle of rank~$da$.
In addition, there exists a continuous multisection 
$\vph\!:\ov\M_{1,k}^0(\P,d;J)\!\lra\!{\cal V}_{1,k}^d$
such~that\\ 
${}\quad$ (${\cal V}1$) $\vph|_{\M_{1,k}^0(\P,d;J)}$ is smooth and transverse to
the zero~set in ${\cal V}_{1,k}^d|_{\M_{1,k}^0(\P,d;J)}$;\\
${}\quad$ (${\cal V}2$) the intersection of $\vph^{-1}(0)$ with each boundary stratum 
${\cal U}_{\cal T}(\P;J)$ and ${\cal U}_{{\cal T};1}^m(\P;J)$\\
${}\qquad\quad~$ of $\ov\M_{1,k}^0(\P,d;J)$ 
is a smooth suborbifold of the stratum of real dimension of\\
${}\qquad\quad~$ at most $2(d(n\!+\!1\!-\!a)\!+\!k)\!-\!2$.\\
If $\vph_0$ and $\vph_1$ are two such multisections, 
there exists a continuous homotopy
$$\Phi\!:[0,1]\!\times\!\ov\M_{1,k}^0(\P,d;J)\lra[0,1]\!\times\!{\cal V}_{1,k}^d$$
such that $\Phi|_{\{t\}\times\ov\M_{1,k}^0(\P,d;J)}\!=\!\vph_t$ for $t\!=\!0,1$, and\\
${}\quad$ (${\cal V}1'$)  $\Phi|_{[0,1]\times\M_{1,k}^0(\P,d;J)}$ is smooth and transverse to
the zero~set in $[0,1]\!\times\!{\cal V}_{1,k}^d|_{\M_{1,k}^0(\P,d;J)}$;\\
${}\quad$ (${\cal V}2'$) the intersection of $\Phi^{-1}(0)$ with 
each boundary stratum $[0,1]\!\times\!{\cal U}_{\cal T}(\P;J)$ and\\ 
${}\qquad\quad~$ $[0,1]\!\times\!{\cal U}_{{\cal T};1}^m(\P;J)$
of $[0,1]\!\times\!\ov\M_{1,k}^0(\P,d;J)$ is a smooth suborbifold of the stratum\\ 
${}\qquad\quad~$ of real dimension of at most $2(d(n\!+\!1\!-\!a)\!+\!k)\!-\!1$.\\
Thus,  the cone ${\cal V}_{1,k}^d$
determines a homology class and a cohomology class on $\ov\M_{1,k}^0(\P,d;J)$:
\begin{gather*}
\hbox{PD}_{\ov\M_{1,k}^0(\P,d;J)}\big(e({\cal V}_{1,k}^d)\big)
\in H_{2(d(n+1-a)+k)}\big(\ov\M_{1,k}^0(\P,d;J);\Q\big)\\
\hbox{and}\qquad
e({\cal V}_{1,k}^d)\in H^{2da}(\ov\M_{1,k}^0(\P,d;J);\Q\big).
\end{gather*}
Finally, if ${\cal W}\!\lra\!\X_{1,k}(\P,d)$ is a vector orbi-bundle
such that the restriction of ${\cal W}$ to each stratum $\X_{\T}(\P)$
of $\X_{1,k}(\P,d)$ is smooth, then
$$\blr{e({\cal W})\cdot e({\cal V}_{1,k}^d),\big[\ov\M_{1,k}^0(\P,d;J)\big]}
=\blr{e({\cal W})\cdot e({\cal V}_{1,k}^d),\big[\ov\M_{1,k}^0(\P,d)\big]}.$$\\
\end{prp}

\noindent
The second statement of this proposition is a special case of Lemma~\ref{g1conemainstr_lmm}.
We use  Lemma~\ref{g1conemainstr_lmm} and Proposition~\ref{g1conebdstr_prp} 
to construct a multisection~$\vph$ satisfying~(${\cal V}1$) and~(${\cal V}2$),
starting from the lowest-dimensional strata of $\ov\M_{1,k}^0(\P,d;J)$.
Suppose ${\cal T}$ and $m$ are as in Lemma~\ref{g1bdstr_lmm}
and we have constructed a neighborhood $U$~of
$$\partial\bar{\cal U}_{\T;1}^m(\P;J)\equiv
\bar{\cal U}_{\T;1}^m(\P;J)-{\cal U}_{\T;1}^m(\P;J)$$
in $\ov\M_{1,k}^0(\P,d;J)$ and a continuous multisection 
$\vph$ of the cone ${\cal V}_{1,k}^d$ over $U$ such~that
for all ${\cal T}'$ and $m'$ as in Lemma~\ref{g1bdstr_lmm}
the restriction of $\vph$ to ${\cal U}_{\T';1}^{m'}(\P;J)\!\cap\!U$
is a smooth multisection of the vector bundle ${\cal V}_{1,k;\T'}^{d;m'}$
of Proposition~\ref{g1conebdstr_prp}
which is transverse to the zero set in~${\cal V}_{1,k;{\cal T}'}^{d;m'}$.
We then extend the restriction of $\vph$ to
${\cal U}_{\T;1}^m(\P;J)\!\cap\!U$ to  a smooth section 
of ${\cal V}_{1,k;\T}^{d;m}$ over ${\cal U}_{\T;1}^m(\P;J)$ and
to a continuous section $\vph_{\T}^m$ of ${\cal V}_{1,k;{\T}}^{d;m}$  over
$\ov\M_{1,k}^0(\P,d;J)\!\cap\!U_{\T}^m$, using 
the bundle isomorphism~$\ti\phi_{\T}^m$ of Proposition~\ref{g1conebdstr_prp}.
By the definition of the bundles ${\cal V}_{1,k;\T}^{k;m}$
in Subsection~\ref{g1conelocalstr_subs2},
the restriction of $\vph_{\T}^m$ to each space  
${\cal U}_{\T';1}^{m'}(\P;J)\!\cap\!U_{\T}^m$
is a section of~${\cal V}_{1,k;\T'}^{k;m'}$,
for all $\T'$ and $m'$ as in Lemma~\ref{g1bdstr_lmm}.
We can also insure that 
the restriction of $\vph_{\cal T}^m$ to ${\cal U}_{{\cal T}';1}^{m'}(\P;J)\cap U_{\cal T}^m$
is smooth and transverse to the zero set in~${\cal V}_{1,k;{\cal T}'}^{d;m'}$.
Finally, by using a partition of unity and the newly constructed section~$\vph_{\cal T}^m$,
we can extend the section $\vph$ 
to a neighborhood of $\bar{\cal U}_{{\cal T};1}^m(\P;J)$ in $\ov\M_{1,k}^0(\P,d;J)$,
without changing it on $\bar{\cal U}_{{\cal T};1}^m(\P;J)$ or
on a neighborhood of $\partial\bar{\cal U}_{{\cal T};1}^m(\P;J)$ in $\ov\M_{1,k}^0(\P,d;J)$.
After finitely many steps, we end up with a neighborhood $U$~of
$$\ov\M_{1,k}^0(\P,d;J)-\M_{1,k}^{\{0\}}(\P,d;J)$$
in  $\ov\M_{1,k}^0(\P,d;J)$ and a continuous multisection $\vph$ 
of the cone ${\cal V}_{1,k}^d$ over $U$ such~that
for all ${\cal T}'$ and $m'$ as in Lemma~\ref{g1bdstr_lmm}
the restriction of $\vph$ to ${\cal U}_{{\cal T}';1}^{m'}(\P;J)$
is a smooth multisection of the vector bundle ${\cal V}_{1,k;{\cal T}'}^{d;m'}$
which is transverse to the zero set in~${\cal V}_{1,k;{\cal T}'}^{d;m'}$.
We then extend $\vph$ in the same stratum-by-stratum way to a section over all
of $\ov\M_{1,k}^0(\P,d;J)$,  using Lemma~\ref{g1conemainstr_lmm}.
Since the real dimension of a boundary stratum ${\cal U}_{\cal T}(\P;J)$
of $\ov\M_{1,k}^0(\P,d;J)$, with ${\cal T}$ as in Lemma~\ref{g1mainstr_lmm},
is at least two less than the dimension of~$\M_{1,k}^0(\P,d;J)$,
the transversality of $\vph|_{{\cal U}_{\cal T}(\P;J)}$ 
to the zero set in ${\cal V}_{1,k}^d$ implies (${\cal V}2$) for  
this stratum.
Similarly, the transversality of $\vph|_{{\cal U}_{{\cal T};1}^m(\P;J)}$ 
to the zero set in ${\cal V}_{1,k;{\cal T}}^{d;m}$ 
and the rank statement of Proposition~\ref{g1conebdstr_prp}
imply (${\cal V}2$) for  each stratum ${\cal U}_{{\cal T};1}^m(\P;J)$
of $\ov\M_{1,k}^0(\P,d;J)$, with ${\cal T}$  and $m$ as in Lemma~\ref{g1bdstr_lmm}.
The homotopy statement of Proposition~\ref{g1cone_prp} is proved 
by a nearly identical construction.\\

\noindent
The second-to-last statement of Proposition~\ref{g1cone_prp} follows from 
the preceding claims by the same argument as in Subsection~\ref{appr_subs}.
The final statement of Proposition~\ref{g1cone_prp} follows from 
the proof of the first part of Proposition~\ref{g1cone_prp} and
from the last statement of Theorem~\ref{g1comp-reg_thm} in~\cite{Z5}.
The latter states that there exists $\de_n(d)\in\!\Bbb{R}^+$ 
with the following property.
If $\under{J}\!=\!(J_t)_{t\in[0,1]}$ 
is a $C^1$-smooth family of almost complex structures on~$\P$ such that
$\|J_t\!-\!J_0\|_{C^1}\!\le\!\de_n(d)$ for all $t\!\in[0,1]$,
then the compact moduli space
$$\ov\M_{1,k}^0(\P,d;\under{J})
\equiv\!\bigcup_{t\in[0,1]}\!\!\!\ov\M_{1,k}^0(\P,d;J_t)
~\subset{\frak X}_{1,k}(\P,d)$$
has the general topological structure of a unidimensional variety with boundary. 
It is stratified by the smooth orbifolds with boundary,
$$\U_{\T}(\P;\under{J})
\equiv\!\bigcup_{t\in[0,1]}\!\!\!{\cal U}_{\cal T}(\P;J_t)
\quad\hbox{and}\quad
{\cal U}_{{\cal T};1}^m(\P;\under{J})
\equiv\!\bigcup_{t\in[0,1]}\!\!\!{\cal U}_{{\cal T};1}^m(\P;J_t),$$
each of dimension one greater than the corresponding dimension
given~by Lemmas~\ref{g1mainstr_lmm} or~\ref{g1bdstr_lmm}.
By the same argument as above, we can construct a multisection
$\Phi$ of the cone ${\cal V}_{1,k}^d$ over 
$\ov\M_{1,k}^0(\P,d;\under{J})$ such that\\
${}\quad$ (${\cal V}1''$) $\Phi|_{\M_{1,k}^0(\P,d;\under{J})}$ is smooth and transverse to
the zero~set in ${\cal V}_{1,k}^d|_{\M_{1,k}^0(\P,d;\under{J})}$;\\
${}\quad$ (${\cal V}2''$) the intersection of $\Phi^{-1}(0)$ with each boundary stratum 
${\cal U}_{\T}(\P;\under{J})$ and ${\cal U}_{{\cal T};1}^m(\P;\under{J})$\\
${}\qquad\quad~$ of $\ov\M_{1,k}^0(\P,d;\under{J})$ is a smooth suborbifold of the stratum
of real dimension of at\\
${}\qquad\quad~$ most $2(d(n\!+\!1\!-\!a)\!+\!k)\!-\!1$,\\
and the restrictions $\vph_0\!\equiv\!\Phi|_{\ov\M_{1,k}^0(\P,d;J_0)}$
and $\vph_1\!\equiv\!\Phi|_{\ov\M_{1,k}^0(\P,d;J_1)}$
satisfy conditions (${\cal V}1$) and~(${\cal V}2$).
If \hbox{${\cal W}\!\lra\!\X_{1,k}(\P,d)$} is a complex vector bundle 
of rank $d(n\!+\!1\!-\!a)\!+\!k$ as in Proposition~\ref{g1cone_prp},
we can then choose a continuous section $F$ of 
${\cal W}$ over $\ov\M_{1,k}^0(\P,d;\under{J})$ such~that\\
${}\quad$ ($\Phi{\cal W}1$) 
$\Phi^{-1}(0)\!\cap\!F^{-1}(0)\!\subset\!\M_{1,k}^0(\P,d;\under{J})$;\\
${}\quad$ ($\Phi{\cal W}2$) 
$F|_{\M_{1,k}^0(\P,d;\under{J})}$ is smooth and transverse to the zero set in
${\cal W}|_{\M_{1,k}^0(\P,d;\under{J})}$;\\
${}\quad$ ($\Phi{\cal W}3$) $F^{-1}(0)$ intersects $\Phi^{-1}(0)$
transversely in~$\M_{1,k}^0(\P,d;\under{J})$,\\
${}\quad$ ($\Phi{\cal W}4$) $f_t^{-1}(0)$ intersects $\vph_t^{-1}(0)$ 
transversely in~$\M_{1,k}^0(\P,d;J_t)$ for $t\!=\!0,1$,\\
${}\qquad\qquad$ where $f_t\!\equiv\!F|_{\M_{1,k}^0(\P,d;J_t)}$.\\
It follows that $\Phi^{-1}(0)\!\cap\!F^{-1}(0)$
is a compact oriented one-dimensional suborbifold of $\M_{1,k}^0(\P,d;\under{J})$
and
\begin{gather*}
\partial\big(\Phi^{-1}(0)\!\cap\!F^{-1}(0)\big)
=\vph_1^{-1}(0)\!\cap\!f_1^{-1}(0)-\vph_0^{-1}(0)\!\cap\!f_0^{-1}(0)\\
\Lra\qquad
^{\pm}\big|\vph_1^{-1}(0)\!\cap\!f_1^{-1}(0)\big|
=~^{\pm}\big|\vph_0^{-1}(0)\!\cap\!f_0^{-1}(0)\big|.
\end{gather*}
This equality implies the last claim of Proposition~\ref{g1cone_prp}.

\subsection{The Local Structure of the Cone $\V_{1,k}^d\!\lra\!\ov\M_{1,k}^0(\P,d;J)$, I}
\label{g1conelocalstr_subs1}

\noindent
In this subsection, we describe the structure of 
the cone ${\cal V}_{1,k}^d\!\lra\!\ov\M_{1,k}^0(\P,d;J)$ over
a neighborhood of each stratum ${\cal U}_{\cal T}(\P;J)$
of Lemma~\ref{g1mainstr_lmm}.

\begin{lmm}
\label{g1conemainstr_lmm}
If $n$, $k$, $d$, $a$, $\L$, and ${\cal V}_{1,k}^d$ are as in Theorem~\ref{main_thm},
there exists $\de_n(d)\!\in\!\R^+$  with the following property.
If $J$ is an almost complex structure on~$\P$, 
such that $\|J\!-\!J_0\|_{C^1}\!<\!\de_n(d)$,
and 
$$\T=([k],I,\aleph;j,\under{d})$$
is a bubble type such that $\sum_{i\in I}\!d_i\!=\!d$ and 
$d_i\!\neq\!0$ for some minimal element $i$ of $I$,
then the requirements of Lemma~\ref{g1mainstr_lmm} are satisfied.
Furthermore, the restriction ${\cal V}_{1,k}^d\!\lra\!{\cal U}_{\cal T}(\P;J)$
is a smooth complex vector orbibundle of rank~$da$.
Finally, there exists a continuous vector-bundle isomorphism
$$\tilde{\phi}_{\cal T}\!:
\pi_{{\cal FT}_{\de}}^*\big({\cal V}_{1,k}^d\big|_{{\cal U}_{\cal T}(\P;J)}\big)
\lra {\cal V}_{1,k}^d\big|_{\ov\M_{1,k}^0(\P,d;J)\cap U_{\cal T}},$$
covering the homeomorphism $\phi_{\cal T}$ of Lemma~\ref{g1mainstr_lmm},
such that $\tilde{\phi}_{\cal T}$ is the identity over 
${\cal U}_{\cal T}(\P;J)$ and is smooth over~${\cal FT}_{\de}^{\eset}$.
\end{lmm}

\noindent
The restriction ${\cal V}_{1,k}^d\!\lra\!{\cal U}_{\cal T}(\P;J)$
is the quotient of the cone $\tilde{\cal V}_{1,k}^d\!\lra\!{\cal U}_{\cal T}^{(0)}(\P;J)$
by the group \hbox{$\hbox{Aut}({\cal T})\!\propto\!(S^1)^{\hat{I}}$};
see Subsection~\ref{notation1_subs} for notation.
The fiber of $\tilde{\cal V}_{1,k}^d$ at a point $b\!=\!(\Si_b,u_b)$ 
of ${\cal U}_{\cal T}^{(0)}(\P;J)$ is 
the Dolbeault cohomology group $H_{\bar{\partial}}^0(\Si_b;u_b^*\L)$,
for a holomorphic structure in the bundle~$u_b^*\L$.
Since $d_i\!\neq\!0$ for some minimal element $i\!\in\!I$,
the degree of the restriction of $u_b^*\L$ to the principal curve of~$\Si_b$ is positive.
Thus, by  an argument similar to Subsections~6.2 and~6.3 in~\cite{Z2},
$$H^1_{\bar{\partial}}(\Si_b;u_b^*\L)=\{0\}\qquad\Lra\qquad
\dim\tilde{\cal V}_{1,k}^d|_b
=\dim  H^0_{\bar{\partial}}(\Si_b;u_b^*\L)
=\ind\!\bpar_{\na,b}=da.$$
Since the holomorphic structure in the line bundles $u_b^*\L$
varies smoothly with $b\!\in\!{\cal U}_{\cal T}^{(0)}(\P;J)$,
it follows that $\tilde{\cal V}_{1,k}^d\!\lra\!{\cal U}_{\cal T}^{(0)}(\P;J)$
is a smooth complex vector bundle of rank~$da$ and 
${\cal V}_{1,k}^d\!\lra\!{\cal U}_{\cal T}(\P;J)$ 
is a smooth complex vector orbibundle of rank~$da$.\\

\noindent
We construct a lift $\ti\phi_{\T}$ of $\phi_{\T}$
to the cone  $\V_{1,k}^d\!\lra\!{\cal U}_{\T}(\P;J)$ as follows.
For each sufficiently small element $\ups\!=\!(b,v)$ of $\ti{\cal F}\T^{\eset}$,
we define the~maps
\begin{alignat*}{2}
&R_{\ups}\!:\Ga(b;\L)\!\equiv\!L^p_1\big(\Si_b;u_b^*\L\big)
\lra\Ga(\ups;\L)\!\equiv\!L^p_1\big(\Si_{\ups};u_{\ups}^*\L\big)
&\quad\hbox{by}\quad
&\big\{R_{\ups}\xi\big\}(z)=\xi\big(q_{\ups}(z)\big)
\quad\hbox{if} ~z\!\in\!\Si_{\ups},\\
&\Pi_{\ups}\!:\Ga\big(\ups;\L\big)\lra
\ti\Ga(\ups;\L)\!\equiv\!L^p_1\big(\Si_{\ups};\tilde{u}_{\ups}^*\L\big)
&\quad\hbox{by}\quad
&\big\{\Pi_{\ups}\xi\big\}(z)=\Pi_{\ze_{\ups}(z)}\xi(z)
\quad\hbox{if} ~z\!\in\!\Si_{\ups},
\end{alignat*}
where $\Pi_{\ze_{\ups}(z)}\xi(z)$ is the $\na$-parallel transport of 
$\xi(z)$ along the $g_n$-geodesic 
$$\ga_{\ze_{\ups}(z)}\!:[0,1]\lra\P,\qquad
\tau\lra\exp_{u_{\ups}(z)}\tau\ze_{\ups}(z),$$
and $\ze_{\ups}\!\in\!\Ga(\ups)$ is as in Subsection~\ref{g1str_subs}.
As in Subsection~\ref{mainres_subs}, we use the \cite{LT}-modified 
$L^p_1-$ and $L^p-$Sobolev norms, defined in the present setting as in 
Subsection~3.3 of~\cite{Z3}.
By a direct computation, for some $C\!\in\!C({\cal U}_{\cal T}(\P;J);\Bbb{R}^+)$,
\begin{gather}\label{coneiso_e1}
\big\|\bpar_{\na,b(\ups)}R_{\ups}\xi\big\|_{\ups,p}
\le C(b)|\ups|^{1/p}\|\xi\|_{b,p,1}
\quad\forall\, \xi\!\in\!\Ga_-(b;\L)\!\equiv\!\ker\bpar_{\na,b}
\qquad\hbox{and}\\
\label{coneiso_e2}\begin{split}
\big\|\Pi_{\ups}^{-1}\!\circ\!\bpar_{\na,\tilde{b}(\ups)}\!\circ\!\Pi_{\ups}\xi
-\bpar_{\na,b(\ups)}\xi\big\|_{\ups,p}
&\le C'(b)\|\ze\|_{\ups,p,1}^2\|\xi\|_{\ups,p,1}\\
&\le C(b)|\ups|^{2/p}\|\xi\|_{b,p,1} \qquad\forall~\xi\!\in\!\Ga(\ups;\L);
\end{split}\end{gather}
see the proof of Corollary~2.3 in~\cite{Z1} for the first inequality in~\e_ref{coneiso_e2}.
We denote by $\Ga_-(\ups;\L)$ the image of $\Ga_-(b;\L)$ under the map~$R_{\ups}$
and by $\Ga_+(\ups;\L)$ its $L^2$-orthogonal complement in $\Ga(\ups;\L)$.
Since the operator
$$\bpar_{\na,b}\!:\Ga(b;\L)\lra\Ga^{0,1}(b;\L)
\equiv\!L^p(\Si_b;\La^{0,1}_{\I,j_b}T^*\Si_b\!\otimes\!u_b^*\L)$$
is surjective for all $b\!\in\!{\cal U}_{\cal T}^{(0)}(\P;J)$,
similarly to Subsection~\ref{g1str_subs} the operator
$$\bpar_{\na,b(\ups)}\!:\Ga_+(\ups;\L)\lra
\Ga^{0,1}(\ups;\L)\!\equiv\!L^p(\Si_{\ups};
\La^{0,1}_{\I,j_{\ups}}T^*\Si_{\ups}\!\otimes\!u_{\ups}^*\L)$$
is an isomorphism if $\ups$ is sufficiently small.
Its norm and the norm of its inverse depend only on 
$[b]\!\in\!{\cal U}_{\cal T}(\P;J)$.
Thus, by~\e_ref{coneiso_e1} and~\e_ref{coneiso_e2}, 
for every $\xi\!\in\Ga_-(b;\L)$ there exists 
a unique $\xi_+(\ups)\!\in\!\Ga_+(\ups;\L)$ such~that
$$\Pi_{\ups}^{-1}\!\circ\bpar_{\na,\tilde{b}(\ups)}\circ\Pi_{\ups}
\big(R_{\ups}\xi\!+\!\xi_+(\ups)\big)=0
\qquad\Llra\qquad
\Pi_{\ups}\big(R_{\ups}\xi\!+\!\xi_+(\ups)\big)\in\ker\bpar_{\na,\tilde{b}(\ups)}.$$
Furthermore, 
$$\|\xi_+(\ups)\|_{\ups,p,1}\le C(b)|\ups|^{2/p}\|\xi\|_{b,p,1},$$
for some $C\!\in\!C({\cal U}_{\cal T}(\P;J);\Bbb{R}^+)$.
We can thus define a smooth lift $\tilde{\phi}_{\cal T}$ 
of the diffeomorphism on $\phi_{\cal T}|_{{\cal FT}_{\de}^{\eset}}$ by
\begin{gather*}
\ti\phi_{\T}\!:
\pi_{{\cal FT}_{\de}^{\eset}}^*\big({\cal V}_{1,k}^d\big|_{{\cal U}_{\cal T}(\P;J)}\big)
\lra {\cal V}_{1,k}^d\big|_{\M_{1,k}^0(\P,d;J)\cap U_{\cal T}},
\qquad 
\ti\phi_{\T}([\ups;\xi])=[\tilde{R}_{\ups}\xi],\\
\hbox{where}\qquad
\tilde{R}_{\ups}\xi=\Pi_{\ups}(R_{\ups}\xi\!+\!\xi_+(\ups)\big).
\end{gather*}
This map extends to a continuous bundle homomorphism
$$\tilde{\phi}_{\cal T}\!:
\pi_{{\cal FT}_{\de}}^*\big({\cal V}_{1,k}^d\big|_{{\cal U}_{\cal T}(\P;J)}\big)
\lra {\cal V}_{1,k}^d\big|_{\ov\M_{1,k}^0(\P,d;J)\cap U_{\cal T}},$$
as can be seen by an argument similar to Subsections~3.9 and~4.1 in~\cite{Z3}.

\subsection{The Local Structure of the Cone ${\cal V}_{1,k}^d\!\lra\!\ov\M_{1,k}^0(\P,d;J)$, II}
\label{g1conelocalstr_subs2}

\noindent
In this subsection, we state the central results of the paper:
Proposition~\ref{g1conebdstr_prp} and Lemma~\ref{g1conebdstr_lmm}.
The former is the analogue of Lemma~\ref{g1bdstr_lmm} for the cone~${\cal V}_{1,k}^d$.
The latter can be viewed as a condensed version of Proposition~\ref{g1conebdstr_prp}. 
The proof of these two results takes up the next two subsections.

\begin{prp}
\label{g1conebdstr_prp}
If $n$, $k$, $d$, $a$, $\L$, and ${\cal V}_{1,k}^d$ are as in Theorem~\ref{main_thm},
there exists $\de_n(d)\!\in\!\R^+$  with the following property.
If $J$ is an almost complex structure on~$\P$
such that $\|J\!-\!J_0\|_{C^1}\!<\!\de_n(d)$,
then the requirements of Lemma~\ref{g1bdstr_lmm} 
and of Lemma~\ref{g1conemainstr_lmm} are satisfied for all appropriate bubble types.
Furthermore, if
$${\cal T}=([k],I,\aleph;j,\under{d})$$
is a bubble type such that $\sum_{i\in I}\!d_i\!=\!d$ and 
$d_i\!=\!0$ for all minimal elements $i$ of $I$, then
the restriction \hbox{${\cal V}_{1,k}^d\!\lra\!{\cal U}_{\cal T}(\P;J)$}
is a smooth complex vector orbibundle of rank~$da\!+\!1$.
In addition, for every integer
$$m\in \big(\max(|\chi({\cal T})|\!-\!n,1),|\chi({\cal T}|\big),$$
there exist a neighborhood $U_{\cal T}^m$ of ${\cal U}_{{\cal T};1}^m(\P;J)$ in 
$\X_{1,k}(\P,d)$ and a topological vector orbibundle 
$${\cal V}_{1,k;{\cal T}}^{d;m} \lra \ov\M_{1,k}^0(\P,d;J)\cap U_{\cal T}^m$$
such~that
${\cal V}_{1,k;{\cal T}}^{d;m}\!\lra\!\M_{1,k}^0(\P,d;J)\cap U_{\cal T}^m$
is a smooth vector orbibundle,
$${\cal V}_{1,k;{\cal T}}^{d;m}\subset{\cal V}_{1,k}^d, \qquad\hbox{and}\qquad
\rk\!{\cal V}_{1,k;{\cal T}}^{d;m}=
da+1-m>\frac{1}{2}\dim{\cal U}_{{\cal T};1}^m(\P;J)
-\big(d(n\!+\!1\!-\!a)\!+\!k\big).$$
There also exists a continuous vector-bundle isomorphism
$$\tilde{\phi}_{\cal T}^m\!:
\pi_{{\cal F}^1{\cal T}_{\de}}^*\big({\cal V}_{1,k;{\cal T}}^{d;m}
\big|_{{\cal U}_{{\cal T};1}(\P;J)\cap U_{\cal T}^m}\big)
\lra {\cal V}_{1,k;{\cal T}}^{d;m}\big|_{\ov\M_{1,k}^0(\P,d;J)\cap U_{\cal T}^m},$$
covering the homeomorphism $\phi_{\cal T}$ of Lemma~\ref{g1bdstr_lmm},
such that $\ti\phi_{\T}^m$ is the identity over ${\cal U}_{\T;1}^m(\P;J)$. 
Finally, 
$${\cal U}_{\T';1}^{m'}(\P;J)\!\cap\!\ov{\cal U}_{\T;1}^m(\P;J)\neq\eset
\quad\Lra\quad
\V_{1,k;{\cal T}'}^{d;m'}\big|_{{\cal U}_{\T;1}^m(\P;J)\cap U_{\T'}^{m'}}
\subset \V_{1,k;{\cal T}}^{d;m}\big|_{{\cal U}_{\T;1}^m(\P;J)\cap U_{\T'}^{m'}}.$$\\
\end{prp}

\noindent
The restriction of every element of ${\cal V}_{1,k}^d|_b$
to the domain of the image of~$b$ under the projection onto
the first component in the decomposition~\e_ref{g1decomp_e1} is a constant function.
Thus, every element of ${\cal V}_{1,k}^d|_b$ is determined by its restriction
to the domain of the image of~$b$ under the projection onto
the second component in~\e_ref{g1decomp_e1}.
The statement concerning the restriction ${\cal V}_{1,k}^d\!\lra\!{\cal U}_{\cal T}(\P;J)$
in Proposition~\ref{g1conebdstr_prp} now follows by the same argument as 
for the corresponding statement in Lemma~\ref{g1conemainstr_lmm},
but applied to the second component in the decomposition~\e_ref{g1decomp_e1}.
The index in this case is $da\!+\!1$.\\

\noindent
The bundle ${\cal V}_{1,k;{\cal T}}^{d;m}\!\lra\!\ov\M_{1,k}^0(\P,d;J)\!\cap\!U_{\cal T}^m$
is not unique. 
However, its restriction to ${\cal U}_{{\cal T};1}^m(\P;J)$~is:
\begin{equation*}\begin{split}
{\cal V}_{1,k;{\cal T}}^{d;m}|_{{\cal U}_{{\cal T};1}^m(\P;J)}
\equiv\big\{\xi\!\in\!{\cal V}_{1,k}^d|_b\!\!: b\!\in\!{\cal U}_{{\cal T};1}^m(\P;J);
~\hbox{if}~b_r\!\in\!\M_{1,k}^0(\P,d;J)~\hbox{and}\,
\lim_{r\lra\i}\!\!b_r\!=\!b\in{\cal U}_{{\cal T};1}^m(\P;J),&\\
\hbox{then}~\exists\, \xi_r\!\in\!{\cal V}_{1,k}^d|_{b_r}
~\hbox{s.t.}\,\lim_{r\lra\i}\!\!\xi_r=\xi&\big\}.
\end{split}\end{equation*}
In other words, ${\cal V}_{1,k;{\cal T}}^{k;m}|_{{\cal U}_{{\cal T};1}^m(\P;J)}$
is the largest subspace of ${\cal V}_{1,k}^d|_{{\cal U}_{{\cal T};1}^m(\P;J)}$
with the property that a continuous lift 
$$\ti\phi_{\T}\!:
\pi_{{\cal F}^1\T_{\de}}^*\big({\cal V}_{1,k;\T}^{d;m}
\big|_{{\cal U}_{\T;1}(\P;J)\cap U_{\T}^m}\big)
\lra {\cal V}_{1,k;\T}^{d;m}\big|_{\ov\M_{1,k}^0(\P,d;J)\cap U_{\T}^m}$$
of $\phi_{\T}$ that restricts to the identity over~${\cal U}_{\T;1}^m(\P;J)$
can possibly exist for a vector-bundle extension for the subspace
${\cal V}_{1,k;\T}^{d;m}|_{{\cal U}_{\T;1}^m(\P;J)}$
to a neighborhood of ${\cal U}_{\T;1}^m(\P;J)$ in~${\cal U}_{\T;1}(\P;J)$.
The next lemma describes the subspace 
${\cal V}_{1,k;\T}^{d;m}|_{{\cal U}_{\T;1}^m(\P;J)}$
of ${\cal V}_{1,k}^d|_{{\cal U}_{\T;1}^m(\P;J)}$ explicitly.
Let
$$\ti{\cal F}^1\T
=\big\{\ups\!\in\!\ti{\cal F}\T\!:[\ups]\!\in\!{\cal F}^1\T\big\}.$$

\begin{lmm}
\label{g1conebdstr_lmm}
Suppose $n$, $k$, $d$, $a$, $\L$, ${\cal V}_{1,k}^d$, $J$, and ${\cal T}$ are 
as in the first and second sentences of Proposition~\ref{g1conebdstr_prp}.
If $b\!\in\!{\cal U}_{\T}^{(0)}(\P;J)$, $\xi\!\in\!\ti\V_{1,k}^d|_b$, 
and $\ups_r\!\in\!\ti{\cal F}^1\T^{\eset}$ is a sequence of gluing parameters
such~that 
$$\lim_{r\lra\i}\!\ups_r=b  \quad\hbox{and}\quad
\lim_{r\lra\i}\big[\big(\rho_i(\ups_r)\big)_{i\in\chi({\cal T})}\big]
=[w]\in\bP\ti\F\T|_b$$
then 
$$  \exists~ [\xi_r]\!\in\!\ti\V_{1,k}^k|_{\phi_{\T}([\ups_r])}
\st \lim_{r\lra\i}\![\xi_r]\!=\![\xi]  \qquad\Llra\qquad
\D_{\T}(\xi\!\otimes\!w)=0.$$
Therefore,
\begin{gather}\label{g1conebdstr_lmm2e}
\V_{1,k;{\T}}^{d;m}|_{{\cal U}_{\T;1}^m(\P;J)}
=\big\{\xi\!\in\!\V_{1,k}^d|_{[b]}\!: [b]\!\in\!{\cal U}_{\T;1}^m(\P;J);~
\D_{\T}(\xi\!\otimes\!w)\!=\!0~\forall w\!\in\!\ti\F^1{\T}_b\big\},\\
\hbox{where}\qquad 
\ti\F^1\T_b=\big\{w\!\in\!\ti\F\T_b\!:\cD_{\T}w\!=\!0\big\}.\notag
\end{gather}
Thus, 
${\cal V}_{1,k;{\cal T}}^{d;m}|_{{\cal U}_{{\cal T};1}^m(\P;J)}\!\lra\!
{\cal U}_{{\cal T};1}^m(\P;J)$ is a smooth complex vector orbibundle
of rank \hbox{$da\!+\!1\!-\!m$}.
\end{lmm}

\noindent
The bundle map $\D_{\T}$ constructed at the end of 
Subsection~\ref{notation1_subs} depends on the choice of connection 
in the bundle $\L\!\lra\!\P$.
It may appear that so do the first two statements of Lemma~\ref{g1conebdstr_lmm}.
This is however not the case for the following reason.
Suppose 
$$b\!\equiv\!(\Si_b,u_b)\in{\cal U}_{\T}^{(0)}(\P;J), \qquad
\xi\!\in\!\ti{\cal V}_{1,k}^d|_b, \qquad
\ups_r\!\in\!\ti{\cal F}^1{\T}^{\eset}, \quad\hbox{and}\quad
w\!\equiv\!(w_i)_{i\in\chi(\T)}\in\ti\F\T|_b$$
are as in Lemma~\ref{g1conebdstr_lmm}.
Then, by the definition of $\ti{\cal F}^1\T^{\eset}$ in Lemma~\ref{g1bdstr_lmm},
\begin{equation}\label{coneconn_e1}
\cD_{\T}(b,w)\!\equiv\!
\sum_{i\in\chi(\T)}\!\psi_{x_{h(i)}(b)}(w_i)\cdot_J du_{b,i}|_{\i}e_{\i}=0
\in T\P                       \qquad\forall\,\psi\!\in\!\E_{\pi_P(b)}.
\end{equation}
On the other hand, since the map $u_b$ is constant on every component $\Si_{b,h}$
of the domain $\Si_b$ of $b$ with $h\!<\!i$ for some $i\!\in\!\chi({\cal T})$, 
$\xi$ is a holomorphic function on $\Si_{b,h}$ 
and thus must be a constant $\xi_P\!\in\!\L_{\ev_P(b)}$.
It follows that
\begin{equation}\label{coneconn_e2}
\xi_{i_1}(\i)=\xi_{i_2}(\i)=\xi_P
\qquad\forall~i_1,i_2\in\chi({\cal T}).
\end{equation}
Suppose that $\na'$ is a connection in the line bundle
$\L\!\lra\!\P$ that induces the same $\bar{\partial}$-operator 
in the line bundle $u_b^*\L\!\lra\!\Si_b$ as 
the connection~$\na$; see Subsection~\ref{mainres_subs}.
Then, there exists a complex-valued one-form~$\th$ on~$\P$ such~that
\begin{equation}\label{coneconn_e3}
\na_v\ze-\na_v'\ze=(\th_qv)\cdot\ze(z)
~~~\forall\, q\!\in\!\P,\,v\!\in\!T_q\P,\,\ze\!\in\!\Ga(\P;\L),
\quad\hbox{and}\quad
u_b^*\th\circ j_b=\I\cdot u_b^*\th.
\end{equation}
Thus, if $\D_{\T}$ and $\D_{\T}'$ are the bundle maps corresponding to 
the connections $\na$ and $\na'$ as at the end of Subsection~\ref{notation1_subs},
\begin{equation}\label{coneconn_e4}\begin{split}
\big\{\D_{\T}(\xi\!\otimes\!w)\!-\!\D_{\T}'(\xi\!\otimes\!w)\big\}(\psi)
&=\sum_{i\in\chi(\T)}\! \psi_{x_{h(i)}(b)}(w_i)\cdot
\big(\th_{\ev_P(b)}(du_{b,i}|_{\i}e_{\i})\big)\cdot\xi_i(\i)\\
&=\th_{\ev_P(b)}\Big(\!
\sum_{i\in\chi(\T)}\!\psi_{x_{h(i)}(b)}(w_i)\cdot_J
(du_{b,i}|_{\i}e_{\i})\Big)\cdot\xi_P=0.
\end{split}\end{equation}
The middle equality above follows from~\e_ref{coneconn_e2}, 
the second condition in~\e_ref{coneconn_e3}, and 
the assumption that $u_b$ is a $J$-holomorphic map.
The last equality above is an immediate consequence of~\e_ref{coneconn_e1}.
More generally, the middle equality in~\e_ref{coneconn_e4} implies that the expression 
$\D_{\T}(\xi\!\otimes\!w)$ is intrinsically defined whenever $w\!\in\!\ti\F^1\T$.\\

\noindent
The second statement of Lemma~\ref{g1conebdstr_lmm} follows immediately from 
the definition of ${\cal V}_{1,k;{\cal T}}^{d;m}|_{{\cal U}_{{\cal T};1}^m(\P;J)}$,
the first statement of Lemma~\ref{g1conebdstr_lmm},
and the last statement of Lemma~\ref{g1bdstr_lmm}.
For the final statement of Lemma~\ref{g1conebdstr_lmm}, let
$$\ti{\cal U}_{\T;1}^m(\P;J)=\big\{b\!\in\!{\cal U}_{\T}^{(0)}(\P;J)\!:
[b]\!\in\!{\cal U}_{\T;1}^m(\P;J)\big\}.$$
By the proof of Lemma~\ref{g1bdstr_lmm},
$$\ti\F^1\T \lra \ti{\cal U}_{\T;1}^m(\P;J)$$
is a vector bundle of rank~$m$.
On the other hand, by the same argument as in Subsection~6.2 of~\cite{Z2},
for every $b\!\in\!{\cal U}_{\T}^{(0)}(\P;J)$ and $i\!\in\!\chi(\T)$,
the linear~map
$$\big\{\xi\!=\!(\xi_h)_{h\in I}\!\in\!\ti\V_{1,k}^d|_b\!:
\xi_i(\i)\!=\!0\big\}\lra\L_{\ev_P(b)}, \qquad \xi\lra\na_{e_{\i}}\xi_i,$$
is surjective. It follows that the linear bundle map 
$$\ti\V_{1,k}^d\lra\hbox{Hom}(\F^1\T,\pi_P^*\E^*\!\otimes\!\ev_P^*\L)$$
over $\ti{\cal U}_{\T;1}^m(\P;J)$ induced by $\D_{\T}$
is surjective on every fiber.
Thus, its kernel is a smooth vector bundle of~rank
$$\rk\!\V_{1,k;\T}^{d;m}=
\rk\!\V_{1,k}^d-\rk\!\hbox{Hom}(\F^1\T,\pi_P^*\E^*\!\otimes\!\ev_P^*\L)
=da+1-m,$$
as claimed in the last statement of Lemma~\ref{g1conebdstr_lmm}.\\

\noindent
We prove the remaining claims of Proposition~\ref{g1conebdstr_prp} and
Lemma~\ref{g1conebdstr_lmm} at the end of Section~\ref{gluing_sec},
after reviewing the multi-step genus-one gluing construction
used in~\cite{Z5} and extending it to the cone~$\V_{1,k}^d$.

\section{A Gluing Construction}
\label{gluing_sec}

\subsection{Smoothing Stable Maps}
\label{gluing_subs1}

\noindent
We begin by reviewing the gluing construction 
of Section~\ref{g1comp-strthm_sec} in~\cite{Z5}.
If $b\!=\!(\Si_b,u_b)$ is any genus-one bubble map such that 
$u_b|_{\Si_{b;P}}$ is constant, 
let $\Si_b^0\!\subset\!\Si_b$ be the maximum connected union of 
the irreducible components of $\Si_b$ such that $\Si_{b;P}\!\subset\!\Si_b^0$ 
and $u_b|_{\Si_b^0}$ is constant.
If $u_b|_{\Si_{b;P}}$ is not constant, let $\Si_b^0\!=\!\eset$.
We~put
\begin{alignat*}{1}
\Ga_B(b) &=\big\{\ze\!\in\!\Ga(\Si_b;u_b^*T\P)\!: \ze|_{\Si_b^0}\!=\!0\big\},\\
\Ga_B(b;\L) &=\big\{\xi\!\in\!\Ga(\Si_b;u_b^*\L)\!: \xi|_{\Si_b^0}\!=\!0\big\},
\qquad\hbox{and}\\
\Ga_B^{0,1}(b;\L) &=\big\{\eta\!\in\!\Ga(\Si_b;\La^{0,1}_{\I,j}T^*\Si_b\!\otimes\!u_b^*\L)\!:
\eta|_{\Si_b^0}\!=\!0\big\}.
\end{alignat*}\\

\noindent
Suppose $\T\!=\!([k],I,\aleph;j,\under{d})$ is a bubble type as in 
Proposition~\ref{g1conebdstr_prp}, i.e.~$d_i\!=\!0$ for all $i\!\in\!I_0$, 
where $I_0\!\subset\!I$ is the subset of minimal elements.
We~put
\begin{gather*}
\chi^0(\T)=\big\{h\!\in\!I\!: d_i\!=\!0~\forall~i\!\le\!h\big\}, \qquad
\chi^-(\T)=\!\bigcup_{i\in\chi(\T)}\!\!\!\{h\!\in\!\hat{I}\!:h\!<\!i\}
 \subset\chi^0(\T),\\
\lr{\T}=\max\big\{\big|\{h\!\in\!\hat{I}\!:h\!\le\!i\}\big|\!: i\!\in\!\chi(\T)\big\}\ge 1,
\qquad \cI_{\lr{\T}}^*=\chi(\T), \qquad \cI_{\lr{\T}}=\hat{I}-\chi(\T)-\chi^-(\T)-I_1,
\end{gather*}
where $I_1\!\subset\!I$ is as in Subsection~\ref{notation1_subs}.
For each $s\!\in\!\{0\}\!\cup\![\lr{\T}\!-\!1]$, let 
$$\cI_s=\big\{i\!\in\!\chi(\T)\!\cup\!\chi^-(\T)\!:
\big|\{h\!\in\!\hat{I}\!:h\!<\!i\}\big|\!=\!s\big\}, \qquad
\cI_s^*=\cI_s\cup\bigcup_{t=0}^{s-1}\big(\cI_t\!\cap\!\chi(\T)\big).$$
In the case of Figure~\ref{g1bdstr_fig} on page~\pageref{g1bdstr_fig}, 
$$\lr{\T}=2,\qquad \cI_0=\{h_1,h_3\}, \qquad \cI_1=\{h_4,h_5\}, \qquad \cI_2=\{h_2\}.$$
In general, the set $\cI_{\lr{\T}}$ could be empty, but the sets $\cI_s$ with
$s\!<\!\lr{\T}$ never are.\\

\noindent
If $b$ is a bubble map of type $\T$ as in Subsection~\ref{notation1_subs} and
$s\!\in\![\lr{\T}]$, we~put
$$\Si_b^{(s)}=\!\!\bigcup_{i\in\chi^0(\T)-\chi^-(\T)}\!\!\!\!\!\!\!\!\!\!\!\!\!\Si_{b,i}~~
\cup \bigcup_{h\in\cI_{s-1}^*}\bigcup_{i<h}\Si_{b,i} ~\subset \Si_b.$$
If $h\!\in\!\cI_{s-1}^*$, let
$$\Si_b^h=\bigcup_{h\le i}\!\Si_{b,i}\subset\Si_b,  \qquad
\chi_h(\T)=\big\{i\!\in\!\chi(\T)\!: h\!\le\!i\big\}, \qquad
\ti{\F}_h\T= {\cal U}_{\T}^{(0)}(X;J) \times \C^{\chi_h(\T)}.$$
If in addition $\ups\!=\!(b,v)\!\in\!\wt{\cal F}\T$, we put
\begin{gather*}
\rho_{s;h}(\ups)=\big(b,(\rho_{h;i}(\ups))_{i\in\chi_h(\T)}\big)\in\ti\F_h\T,
\qquad\hbox{where}\quad
\rho_{h;i}(\ups)=\!\!\prod_{h<h'\le i}\!\!\!\!\! v_{h'} \in\C;\\
\cI_{s-1}^0(\ups)=\big\{h\!\in\!\cI_{s-1}^*\!:\rho_{s;h}(\ups)\!=\!0\big\};
\end{gather*}
see Subsection~\ref{notation1_subs} for notation.\\

\noindent
If $\ups\!=\!(b,v)\in\ti{\cal F}\T$, let
$$\ups_0=\big(b,v_{\aleph},(v_h)_{h\in I_1}\big)
\qquad\hbox{if}\quad
\ups=\big(b,v_{\aleph},(v_h)_{h\in\hat{I}}\big).$$
Let $\ups_{\lr{0}}\!=\!\ups$ and $\ups_{\lr{\lr{\T}+1}}\!=\!b$.
If $s\!\in\![\lr{\T}]$, we put
$$\ups_s=\big(b,(v_h)_{h\in\cI_s}\big) \qquad\hbox{and}\qquad
\ups_{\lr{s}}=\big(b,(v_h)_{h\in\cI_t,t\ge s}\big).$$
The component $\ups_{\lr{\T}}$ of $\ups$ consists of smoothings 
at the nodes of $\Si_b$ that do not lie on the principal component~$\Si_{b;P}$ of $\Si_b$
and do not lie  between  $\Si_{b;P}$ and the bubble components indexed by 
the set~$\chi(\T)$.
In Section~\ref{g1comp-strthm_sec} of~\cite{Z5},
these nodes are smoothed out at the first step of the gluing construction,
as specified by~$\ups_{\lr{\T}}$.
After that, the nodes indexed by the set $\cI_{\lr{\T}-1}$ are smoothed out, and so on.
At the last step, the nodes that lie on the principle component are smoothed 
according to~$\ups_0$, provided $\ups\!\in\!\ti{\cal F}^1\T^{\eset}$ 
is sufficiently small.\\ 

\noindent
If $\ups\!\in\!\ti{\cal F}^1\T^{\eset}$ is sufficiently small and 
$s\!\in\!\{0\}\!\cup\![\lr{\T}]$, let
$$q_{\ups_{\lr{s}}}\!:  \Si_{\ups_{\lr{s}}}\lra\Si_b$$
be the basic gluing map constructed in Subsection~2.2 of~\cite{Z3}.
Via the construction of Subsection~3.3 in~\cite{Z3}, the map $q_{\ups_{\lr{s}}}$
induces a metric $g_{\ups_{\lr{s}}}$ and a weight function~$\rho_{\ups_{\lr{s}}}$
that define weighted $L^p_1$-Sobolev norms $\|\cdot\|_{\ups,p,1}$ on the spaces
$\Ga_B(b')$ and $\Ga_B(b';\L)$ and weighted $L^p$-Sobolev norms $\|\cdot\|_{\ups,p}$ on 
the corresponding spaces of differentials, 
for any bubble map $b'\!=\!(\Si_{\ups_{\lr{s}}},u)$ such that
$u$ is constant on $q_{\ups_{\lr{s}}}^{-1}(\Si_b^{(s)})$ if $s\!>\!0$.
In this case, $(\Si_{\ups_{\lr{s}}},g_{\ups_{\lr{s}}})$ is obtained
from $\Si_b$ with its metric~$g_b$ by replacing the nodes of $\Si_b$ indexed
by the sets $\cI_t$ with $t\!\ge\!s$ by thin necks.
The norms  $\|\cdot\|_{\ups,p,1}$ and $\|\cdot\|_{\ups,p}$ are analogous to 
the ones used in Section~3 of~\cite{LT}.
Let 
$$q_{\ups_s;\lr{\T}+1-s}\!: \Si_{\ups_{\lr{s}}}\lra\Si_{\ups_{\lr{s+1}}}$$
be the basic gluing map of Subsection~2.2 in~\cite{Z3} corresponding 
to the gluing parameter $\ups_s$.
We recall that 
$$q_{\ups_{\lr{s}}}=q_{\ups_{\lr{s+1}}} \circ q_{\ups_{s};\lr{\T}+1-s}$$
for all $s\!\in\!\{0\}\!\cup\![\lr{\T}\!-\!1]$.
If $s\!\in\![\lr{\T}]$ and $h\!\in\!\cI_{s-1}^*$, let
$$\Si_{\ups_{\lr{s}}}^h=q_{\ups_{\lr{s}}}^{-1}\big(\Si_b^h\big)
\subset \Si_{\ups_{\lr{s}}}.$$
We note that $\Si_{\ups_{\lr{s}}}^h$ is a union of components of $\Si_{\ups_{\lr{s}}}$. \\

\noindent
For any $\ups\!=\!(b,v)\!\in\!\ti{\cal F}\T$, we put
$$\ti{b}_{\lr{\T}+1}(\ups) \equiv \big( \Si_b,\ti{u}_{\ups,\lr{\T}+1}\big)
=(\Si_b,u_b).$$
In Section~\ref{g1comp-strthm_sec} of~\cite{Z5}, for $J$ sufficiently close to~$J_0$, 
$\de\!\in\!C({\cal U}_{\T}(\P;J);\R^+)$ sufficiently small
and all $\ups\!\in\!\wt{\cal F}^1\T_{\de}^{\eset}$, we construct
$J$-holomorphic bubble maps 
$$\ti{b}_s(\ups)=(\Si_{\ups_{\lr{s}}},\ti{u}_{\ups,s}\big) \qquad
\forall~s=0,\ldots,\lr{\T}$$
such that the following properties are satisfied.
First, for all $s\!\in\![\lr{\T}]$,
\begin{equation}\label{mapgluing_e1}
\Si_{\ti{b}_s(\ups)}^0=q_{\ups_{\lr{s}}}^{-1}\big(\Si_b^{(s)}\big)
\qquad\hbox{and}\qquad
\ti{u}_{\ups,s}\big(\Si_{\ti{b}_s(\ups)}^0\big) =u_b(\Si_b^0) \!\equiv\!\ev_P(b).
\end{equation}
Second, for all $s\!\in\![\lr{\T}]$,
\begin{equation}\label{mapgluing_e2}\begin{split}
&\qquad \ti{u}_{\ups,s}=\exp_{u_{\ups,s}}\!\ze_{\ups,s}\\
\hbox{for some}\qquad 
&\ze_{\ups,s}\in\Ga_B\big(b_s(\ups)\big)
~~\hbox{s.t.}~~
\big\|\ze_{\ups,s}\big\|_{\ups_{\lr{s}},p,1} \le  C(b)|\ups|^{1/p},
\end{split}\end{equation}
where
$$b_s(\ups)=\big(\Si_{\ups_{\lr{s}}},u_{\ups_,s}\big), \qquad
u_{\ups,s}=\ti{u}_{s+1}\circ q_{\ups_s;\lr{\T}+1-s}.$$
Third, 
\begin{equation}\label{mapgluing_e3}\begin{split}
&\qquad \ti{u}_{\ups,0}=\exp_{u_{\ups,0}}\!\ze_{\ups,0}\\
\hbox{for some}\qquad 
&\ze_{\ups,0}\in\Ga_B\big(b_0(\ups)\big)
~~\hbox{s.t.}~~
\big\|\ze_{\ups,0}\big\|_{\ups,p,1} \le  C(b)\big|\rho(\ups)\big|,
\end{split}\end{equation}
where
$$b_0(\ups)=\big(\Si_{\ups},u_{\ups,0}\big), \qquad
u_{\ups,0}=\ti{u}_1\circ \ti{q}_{\ups_s;\lr{\T}+1},$$
and 
$$\ti{q}_{\ups_s;\lr{\T}+1}\!: \Si_{\ups}\lra \Si_{\ups_{\lr{1}}}$$
is the modified gluing map corresponding to the parameter of $\de(b)^{1/2}$
constructed in Subsection~\ref{g1comp-reg1_subs2} of~\cite{Z5}. 
Finally, the maps $\ups\!\lra\!\ze_{\ups,s}$ are smooth over 
$\ti{\cal F}^1\T_{\de}^{\eset}$ and extend continuously over $\ti{\cal F}^1\T_{\de}$.
These extensions satisfy
\begin{gather}\label{mapgluing_e4a}
\ze_{b,s}=0\quad\forall\,b\!\in\!{\cal U}_{\T}^{(0)}(X;J),\,s\!\in\!\{0\}\!\cup\![\lr{\T}]
\qquad\hbox{and}\\
\label{mapgluing_e4b}
\ze_{\ups,s}|_{\Si_{\ups_{\lr{s}}}^h}=0 \quad
\forall~\ups\!\in\!\ti{\cal F}^1\T,~s\!\in\![\lr{\T}],~h\!\in\!\cI_{s-1}^0(\ups).
\end{gather}
The homeomorphism of Lemma~\ref{g1bdstr_lmm} is given by
$$\phi_{\T}\!:{\cal F}^1{\T}_{\de}\lra\ov\M_{1,k}^0(\P,d;J)\cap U_{\T},
\qquad \phi_{\T}\big([\ups]\big)=\big[\ti{b}_0(\ups)\big].$$\\

\noindent
{\it Remark:} The bubble maps $b_s(\ups)$ and $\ti{b}_s(\ups)$ above correspond to 
the bubble maps $b_s(\ti\mu_0(\ups,\ze_{\ups,0}))$ and $\ti{b}_s(\ti\mu_0(\ups,\ze_{\ups,0}))$
in Section~\ref{g1comp-strthm_sec} of~\cite{Z5}, where $\ti\mu_0(\ups,\ze_{\ups,0})$
is the perturbation of~$\ups$ constructed in Subsection~\ref{g1comp-strthm_subs} in~\cite{Z5}.

\subsection{Smoothing Bundle Sections, I}
\label{gluing_subs2}

\noindent
In this subsection we extend all but the last step of the gluing construction 
summarized above to the cone $\V_{1,k}^d$ over~$\ov\M_{1,k}(\P,d;J)$.\\

\noindent
In order to do this, we will use a convenient family of connections
in the line bundles $u^*\L\!\lra\!\Si$, which is chosen in Lemma~\ref{conn_lmm} below. 
First, if $b\!=\!(\Si_b,u_b)$ is a stable genus-one bubble map such that 
$u_b|_{\Si_{b;P}}$ is constant,
$g_b$ is a Hermitian metric  in the line bundle $u_b^*\L\!\lra\!\Si_b$,
and  $\na^b$ is a connection in~$u_b^*\L$, we will call the pair
{\tt $(g,\na)$-admissible}~if\\
${}\quad$ ($g\na1$) $\na^b$ is $g_b$-compatible and $\bpar_{\na,b}$-compatible;\\
${}\quad$ ($g\na2$) $g_b\!=\!g_{u_b}$ and $\na^b\!=\!\na^{u_b}$ on $\Si_b^0$,\\
where $g_{u_b}$ is the Hermitian metric in $u_b^*\L$ induced
from the standard metric in~$\L$.
The second condition in~($g\na1$) means~that
$$\bpar_{\na,b}
\equiv\frac{1}{2}\big(\na^{u_b}+\I\na^{u_b}\circ j\big)
=\frac{1}{2}\big(\na^b+\I\na^b\circ j\big),$$ 
with notation as in~\e_ref{vdfn_e}.
If the pair $(g,\na)$ satisfies ($g\na1$), the connection~$\na^b$ is uniquely determined
by the metric~$g_b$.
The second conditions in~($g\na1$) and in~($g\na2$) imply that
the bundle map $\D_{\T}$ does not change if it is defined using 
the connection~$\na^b$ instead of~$\na^{u_b}$; see Subsection~\ref{notation1_subs}.\\

\noindent
If $b\!\in\!{\cal U}_{\T}^{(0)}(\P;J)$, $\de\!\in\!\R^+$, $i\!\in\!\hat{I}$, let
\begin{alignat*}{1}
A_{b,i}^-(\de)&= \big\{(i,z)\!\in\!\Si_{b,i}\!=\!\{i\}\!\times\!S^2\!: 
|z|\!\ge\!\de^{-1/2}/2\big\}\subset \Si_b; \\
\partial^- A_{b,i}^-(\de)&= \big\{(i,z)\!\in\!\Si_{b,i}\!=\!\{i\}\!\times\!S^2\!: 
|z|\!=\!\de^{-1/2}/2\big\}\subset \Si_b;\\
\Si_b^0(\de)&=\!\bigcup_{i\in\chi(\T)}\!\!\!\!A_{b,i}^-(\de)
\cup\bigcup_{h\in\chi^0(\T)}\!\!\!\!\!\!\Si_{b,h}.
\end{alignat*}
If $\ups\!\in\!\ti{\cal F}\T$ is sufficiently small, we put
$$A_{\ups,i}^-(\de)=q_{\ups}^{-1}(A_{b,i}^-(\de))\subset\Si_{\ups},
\quad
\partial^-A_{\ups,i}^-(\de)=q_{\ups}^{-1}\big(\partial^-A_{b,i}^-(\de)\big),
\quad  \Si_{\ups}^0(\de)=q_{\ups}^{-1}\big(\Si_b^0(\de)\big).$$
If $s\!\in\![\lr{\T}\!+\!1]$ and $h\!\in\!\cI_{s-1}^*$, let
$$\Si_{\ups_{\lr{s}}}^{h;0}(\de)=
\Si_{\ups_{\lr{s}}}^0(\de)\cap\Si_{\ups_{\lr{s}}}^h.$$

\begin{lmm}
\label{conn_lmm}
If $n$, $d$, $k$, $a$, and $\L$ are as in Proposition~\ref{g1conebdstr_prp},  
there exists $\de_n(d)\!\in\!\R^+$ such that for every almost complex structure $J$ on~$\P$,
such that $\|J\!-\!J_0\|_{C^1}\!\le\!\de_n(d)$, and a bubble type $\T$  as above,
there exist $\de,C\!\in\!C({\cal U}_{\T}(\P;J);\R^+)$ with the following property.
For~every 
$$\ups\!\equiv\!(b,v)\in\ti{\cal F}^1\T_{\de}^{\eset}
\qquad\hbox{and}\qquad s\in [\lr{\T}\!+\!1],$$
there exist a metric $g_{\ups,s}$ and 
a connection $\na^{\ups,s}$ in the line bundle
$\ti{u}_{\ups,s}^*\L\!\lra\!\Si_{\ups}$ such~that\\
${}\quad$ (1) all pairs $(g_{\ups,s},\na^{\ups,s})$ are admissible;\\
${}\quad$ (2) the curvature of $\na^{\ups,s}$ vanishes on $\Si_{\ups_{\lr{s}}}^0(2\de(b))$.\\
Furthermore, the maps $\ups\!\lra\!(g_{\ups,s},\na^{\ups,s})$
are $\Aut(\T)\!\propto\!(S^1)^I$-invariant and smooth over~$\ti{\cal F}^1\T_{\de}^{\eset}$.
They extend continuously over~$\ti{\cal F}^1\T_{\de}$.
The extension satisfies (1) and~(2).
In addition, 
\begin{gather}\label{connlmm_e1}
(g_{b,s},\na^{b,s}) = (g_{b,\lr{\T}+1},\na^{b,\lr{\T}+1})
\qquad\forall\,b\!\in\!{\cal U}_{\T}^{(0)}(\P;J),\,s\!\in\![\lr{\T}];\\ 
\label{connlmm_e2}
(g_{\ups,s},\na^{\ups,s})\big|_{\Si_{\ups_{\lr{s}}}^h}
=q_{\ups_s;\lr{\T}+1-s}^*
(g_{\ups,s+1},\na^{\ups,s+1})\big|_{\Si_{\ups_{\lr{s}}}^h}
\qquad\forall\, s\!\in\![\lr{T}],\,h\!\in\!\cI_{s-1}^0(\ups).
\end{gather}\\
\end{lmm}

\noindent
This lemma is an analogue of Lemma~\ref{g1comp-reg0_lmm3} in~\cite{Z5}
for the bundle $\L$ and is proved in a similar way.
Let $\be\!:\R^+\!\lra\![0,1]$ be a smooth function such that
$$\be(t)\in
\begin{cases}
0,&\hbox{if}~t\!\le\!1;\\
1,&\hbox{if}~t\!\ge\!2.\\
\end{cases}$$
If $r\!\in\!\R^+$, let $\be_r(t)\!=\!\be(t/\sqrt{r})$.
We define $\be_b\!\in\!C^{\i}(\Si_b;\R)$ by
\begin{equation}\label{cutoff_e}
\be_b(z)=
\begin{cases}
1,&\hbox{if}~z\!\in\!\Si_{b,i},~i\!\in\!\chi^0(\T);\\
1-\be_{\de(b)}(r(z)/2),&\hbox{if}~z\!\in\!\Si_{b,i},~i\!\in\!\chi(\T);\\
0,&\hbox{otherwise},
\end{cases}\end{equation}
where $r(z)\!=\!|q_S^{-1}(z)|$ if $q_S\!:\C\!\lra\!S^2$ is
the stereographic projection mapping the origin to the south pole of~$S^2$.
In other words, $\be_b\!=\!1$ on $\Si_b^0(2\de(b))$ and vanishes
outside of $\Si_b^0(8\de(b))\!\subset\!\Si_b$.
Let $\be_{\ups}\!=\!\be_b\!\circ\!q_{\ups}$.\\

\noindent
For $s\!\in\![\lr{\T}\!+\!1]$, $h\!\in\!\cI_{s-1}^*$, and 
$\ups\!\in\!\ti{\cal F}^1\T_{\de}^{\eset}$,
we use parallel transport with respect to the connection $\na^{\ti{u}_{\ups,s}}$
along the meridians to the south pole of the sphere $\Si_{\ups_{\lr{s}}}^h$
to identify $\ti{u}_{\ups,s}^*\L$ over $\Si_{\ups_{\lr{s}}}^{h;0}(8\de(b))$
with the trivial holomorphic line bundle 
$$\Si_{\ups_{\lr{s}}}^{h;0}(8\de(b))\times\L_{\ev_P(b)}.$$
A connection $\na^{\ups,s}$ with the desired properties can then 
be found by solving an equation of the~form
\begin{equation}\label{conn_e}
\bar{\partial}\th=\be_{\ups}\Om_{\ups,h}, \qquad \th(\i)=0, \qquad
\th\in C^{\i}(\Si_{\ups_{\lr{s}}}^{h;0}(8\de(b));\C),
\end{equation}
where $\Om_{\ups,h}\!\in\!C^{\i}(\Si_{\ups}^{h;0}(8\de(b));\C)$
is determined by $\ups$ and satisfies
$$\|\Om_{\ups,h}\|_{\ups_{\lr{s}},p}\le C(b)\de(b)^{1/p}.$$
This bound follows immediately from the definition of the set~$\chi(\T)$
and~\e_ref{mapgluing_e2}.
The equation~\e_ref{conn_e} can be viewed as an equation on~$\Si_{\ups_{\lr{s}}}^h$,
which is a two-sphere with the metric~$g_{\ups_{\lr{s}}}$ arising in the pregluing construction.
If $\de(b)\!\in\!\R^+$ is sufficiently small, 
\e_ref{conn_e} has a unique solution $\th_{\ups,h}\!\in\!C^{\i}(\Si_{\ups_{\lr{s}}}^h;\C)$.
The curvature of the connection  
$$\ti\na^{\ups,h}\equiv \na^{\ti{u}_{\ups,s}}+\be_{\ups}\th_{\ups,h}$$
then vanishes on $\Si_{\ups_{\lr{s}}}^{h;0}(2\de(b))$.\\

\noindent
Let $g_{\ups,h}$ be the metric in $\ti{u}_{\ups,s}^*\L|_{\Si_{\ups_{\lr{s}}}^h}$ obtained by 
patching the flat metric in $\ti{u}_{\ups,s}^*\L|_{\Si_{\ups_{\lr{s}}}^{h;0}(8\de(b))}$
induced via parallel transport from $\i\!\in\!\Si_{\ups_{\lr{s}}}^h$ with respect to 
$\ti\na^{\ups,h}$ with  the metric $g_{\ti{u}_{\ups,s}}$ 
over 
$$\Si_{\ups_{\lr{s}}}^{h;0}(8\de(b))-\Si_{\ups_{\lr{s}}}^{h;0}(4\de(b)).$$
We put
$$g_{\ups,s}|_z= \begin{cases}
g_{\ups,h}|_z,&\hbox{if}~z\!\in\!\Si_{\ups_{\lr{s}}}^h,~h\!\in\!\cI_{s-1}^*;\\
g_{\ti{u}_{\ups,s}}|_z,&\hbox{if}~z\!\in\!\Si_{\ti{b}_s(\ups)}^0.
\end{cases}$$
Since $\Si_{\ti{b}_s(\ups)}^0$ is the union of the components of $\Si_{\ups_{\lr{s}}}$
that are not in $\Si_{\ups_{\lr{s}}}^h$ for any $h\!\in\!\cI_{s-1}^*$ by~\e_ref{mapgluing_e1},
the metric $g_{\ups,s}$ on $\ti{u}_{\ups,s}^*\L$ is well-defined.
In particular, the two definitions agree at the node of~$\Si_{\ups_{\lr{s}}}^h$.
Let $\na^{\ups,s}$ be the $\bpar_{\na,\ti{u}_{\ups,s}}$-compatible
and $g_{\ups,s}$-compatible connection.
By construction,
$\na^{\ups,s}\!=\!\ti\na^{\ups,h}$ on~$\Si_{\ups_{\lr{s}}}^{h;0}(2\de(b))$.
Thus, the pair $(g_{\ups,s},\na^{\ups,s})$ satisfies 
the requirements (1) and~(2) of Lemma~\ref{conn_lmm}.
By construction,  the map $\ups\!\lra\!(g_{\ups,s},\na^{\ups,s})$
is $\Aut(\T)\!\propto\!(S^1)^I$-invariant and smooth.
Since the maps $\ups\!\lra\!\ze_{\ups,s}$ extend continuously over $\ti{\cal F}^1\T_{\de}$,
so does the map  $\ups\!\lra\!(g_{\ups,s},\na^{\ups,s})$, as can be seen by an argument
analogous to Subsections~3.9 and~4.1 in~\cite{Z3}.
It is immediate from the construction that~\e_ref{connlmm_e1} is satisfied,
while~\e_ref{connlmm_e2} follows from~\e_ref{mapgluing_e4b}.\\

\noindent
For each $s\!\in\![\lr{\T}]$, we will next choose a family of identifications 
$$\Pi_{\ups,s}\big|_z\!: u_{\ups,s}^*\L|_z\lra \ti{u}_{\ups,s}^*\L_z,
\qquad z\in\Si_{\ups_{\lr{s}}},$$
which is smooth in $\ups$ on $\ti{\cal F}^1\T^{\eset}$ and in~$z$.
If $z\!\in\!\Si_{\ups_{\lr{s}}}^{h;0}(2\de(b))$  for some $h\!\in\!\cI_{s-1}^*$, 
let $\Pi_{\i,z}^{\ups,s}$ and $\ti\Pi_{\i,z}^{\ups,s}$ 
be the parallel transports in the line bundles $u_{\ups,s}^*\L$
and $\ti{u}_{\ups,s}^*\L$, respectively, along a path from $\i$ to $z$ in 
$\Si_{\ups_{\lr{s}}}^{h;0}(2\de(b))$ with respect to 
the connections $q_{\ups_s;\lr{\T}+1-s}^*\na^{\ups,s+1}$ and~$\na^{\ups,s}$.
Due to the requirement~(2) of Lemma~\ref{conn_lmm},
these parallel transports are path-independent.
If $z\!\in\!\Si_{\ups_{\lr{s}}}$ and $\xi\!\in\!u_{\ups,s}^*\L|_z$, 
we require that
\begin{equation}\label{bundiden_e}
\Pi_{\ups,s}|_z\xi= \begin{cases}
\xi,&\hbox{if}~z\!\in\!\Si_{b_s(\ups)}^0;\\
\ti\Pi_{\i,z}^{\ups,s}\big\{\Pi_{\i,z}^{\ups,s}\big\}^{-1}\!\xi, &\hbox{if}~
z\!\in\!\Si_{\ups_{\lr{s}}}^{h;0}(\de(b)),~h\!\in\!\cI_{s-1}^*;\\
\Pi_{\ze_{\ups,s}(z)}\xi, &\hbox{if}~
z\!\not\in\!\Si_{\ups_{\lr{s}}}^{h;0}(2\de(b))~\forall\,h\!\in\!\cI_{s-1}^*.
\end{cases}
\end{equation}
We patch the last two identifications in~\e_ref{bundiden_e} over 
$\Si_{\ups_{\lr{s}}}^{h;0}(2\de(b))\!-\!\Si_{\ups_{\lr{s}}}^{h;0}(\de(b))$,
using a cutoff function constructed from~$\be$.
Let
$$\Pi_{\ups,s}\!:\Ga\big(\Si_{\ups_{\lr{s}}};u_{\ups,s}^*\L\big)
\lra\Ga\big(\Si_{\ups_{\lr{s}}};\ti{u}_{\ups,s}^*\L\big)$$
be the operator induced by the maps $\Pi_{\ups,s}|_z$.
We note that if
$$\xi\in \ti\Ga_-\big(\ups_{\lr{s}};\L)\!\equiv\!
\ker\bpar_{\na,\ti{b}(\ups)},$$
then $\{\ti\Pi_{\i,\cdot}^{\ups,s}\}^{-1}\xi$ is a holomorphic function on
$\Si_{\ups_{\lr{s}}}^{h;0}(2\de(b))$,
since covariant differentiation commutes 
with parallel transport due to~(2) of Lemma~\ref{conn_lmm}.\\

\noindent
For $b\!\in\!{\cal U}_{\T}^{(0)}(X;J)$, $s\!\in\![\lr{\T}]$, and 
$h\!\in\!\cI_{s-1}^*$, let
\begin{gather*}
\Ga_h(b;\L)=\big\{\xi\!\in\!\Ga_B(b;\L)\!: \xi|_{\Si_b-\Si_b^h}\!=\!0\big\},
\qquad \Ga_{h;-}(b;\L)=\Ga_h(b;\L)\cap\Ga_-(b;\L),\\
\Ga_h^{0,1}(b;\L)=\big\{\eta\!\in\!\Ga_B^{0,1}(b;\L)\!:\eta|_{\Si_b-\Si_b^h}\!=\!0\big\}.
\end{gather*}
If $\ups\!\in\!\ti{\cal F}^1\T_{\de}^{\eset}$, $s\!\in\![\lr{\T}\!+\!1]$, and 
$h\!\in\!\cI_{s-1}^*$, we~put
$$\ti\Ga_{h;-}\big(\ups_{\lr{s}};\L\big)=\big\{\xi\!\in\!\ti\Ga_-\big(\ups_{\lr{s}};\L)\!:
\xi|_{\Si_{\ups_{\lr{s}}}-\Si_{\ups_{\lr{s}}}^h}\!=\!0\big\}.$$
For each $m\!\in\!\Z^+$, we define
$$\D_{s;h}^{(m)}\!: \ti\Ga_-\big(\ups_{\lr{s}};\L\big) \lra \L_{\ev_P(b)}
\qquad\hbox{by}\qquad
\D_{s;h}^{(m)}\xi=\frac{d}{dw_h}\{\ti\Pi_{\i,\cdot}^{\ups,s}\}^{-1}\xi\Big|_{w_h=0}
=\big\{\na_{e_{\i}}^{\ups,s}\big\}^m\xi|_{\Si_{\ups_{\lr{s}}}^h},$$
where $w_h$ is the standard holomorphic coordinate around $\i$ in~$\Si_{\ups_{\lr{s}}}^h$.
We will construct isomorphisms 
$$\ti{R}_{\ups,s}\!: \Ga_-(b;\L)\!\equiv\!\ker\,\bpar_{\na,b}\lra
\ti\Ga_-\big(\ups_{\lr{s}};\L\big) \qquad\forall~s\!\in\!\big[\lr{\T}\big]$$
such that the following properties are satisfied.
First, for all $h\!\in\!\cI_{s-1}^*$,
\begin{equation}\label{bundleisom_e0}
\ti{R}_{\ups,s}\xi\in\ti\Ga_{h;-}\big(\ups_{\lr{s}};\L\big)
\qquad\forall\, \xi\!\in\!\Ga_{h;-}(b;\L).
\end{equation}
Second, for all $h\!\in\!\cI_{s-1}^*$,
\begin{equation}\label{bundleisom_e1}
\D_{s;h}^{(1)}\ti{R}_{\ups,s}\xi =
\al_{s;h}\big(\rho_{s;h}(\ups);\xi\big) \equiv
\!\sum_{i\in\chi_h(\T)} \!\!\! \rho_{h;i}(\ups)\,\D_{\T,i}\xi
\qquad\forall\, \xi\!\in\!\Ga_-(b;\L).
\end{equation}
Finally, the maps $\ups\!\lra\!\ti{R}_{\ups,s}$ are smooth over 
$\ti{\cal F}^1\T_{\de}^{\eset}$ and extend continuously over $\ti{\cal F}^1\T_{\de}$.
These extensions satisfy
\begin{equation}\label{bundleisom_e2}
\ti{R}_{\ups,s}|_b\!=\!\id\!:  \Ga_-(b;\L) \lra \Ga_-(b;\L).
\end{equation}\\

\noindent
In order to construct isomorphisms $\ti{R}_{\ups,s}$, we observe that certain
operators are surjective.
If $b\!\in\!{\cal U}_{\T}^{(0)}(X;J)$, $s\!\in\![\lr{\T}]$, $h\!\in\!\cI_{s-1}^*$, 
and $[w_h]\!\in\!\bP\ti\F_h\T|_b$, let 
$$\Ga_{h;-}\big(b;\L;[w_h]\big)=\big\{\xi\!\in\!\Ga_{h;-}(b;\L)\!: 
\al_{s;h}(w_h;\xi)\!=\!0\big\}.$$
We denote the $L^2$-orthogonal complement of $\Ga_{h;-}(b;\L;[w_h])$
in $\Ga_{h;-}(b;\L)$ by  $\Ga_{h;-}^{\perp}(b;\L;[w_h])$.
Since $\Si_b^h$ is a genus-zero Riemann surface and the degree of $u_b^*\L$ over
every component of $\Si_b^h$ is nonnegative,
$$H^1\big(\Si_b^h;\{u_b|_{\Si_b^h}\}^*\L\!\otimes\!\O(-z)\big)=\{0\}
\qquad\forall~z\!\in\!\Si_b^{h*},$$
where $\Si_b^{h*}\!\subset\!\Si_b^h$ is the subset of smooth points.
Thus, the operator
$$\bpar_{\na,b}^h\!: \Ga_h(b;\L)\lra \Ga_h^{0,1}(b;\L)$$
induced by $\bpar_{\na,b}$ is surjective.
Similarly, since the degree of $u_b^*\L|_{\Si_{b,i}}$ is positive for all $i\!\in\!\chi(\T)$,
$$H^1\big(\Si_b^h;\{u_b|_{\Si_b^i}\}^*\L\!\otimes\!\O(-2z)\big)=\{0\}
\qquad\forall~z\!\in\!\Si_b^{i*}\cap\Si_b^h.$$
Thus, for every element $w_h\!\in\!\ti\F_h\T|_b$, the linear map
$$\al_{s;h}(w_h;\cdot)\!: \Ga_{h;-}(b;\L)\lra\L_{\ev_P(b)}$$
is surjective. In particular,
$$\al_{s;h}(w_h;\cdot)\!: \Ga_{h;-}^{\perp}(b;\L;[w_h]) \lra\L_{\ev_P(b)}$$
is an isomorphism and 
\begin{equation}\label{bundleisom_e3}
C(b)^{-1}|w_h||\xi| \le \big|\al_{s;h}(w_h;\xi)\big|
\le C(b)|w_h||\xi| \qquad\forall~\xi\in\Ga_{h;-}^{\perp}(b;\L;[w_h]),
\end{equation}
for some $C\!\in\!C({\cal U}_{\T}(\P;J);\R^+)$.\\

\noindent
If $\ups\!\!\in\!\ti{\cal F}^1\T_{\de}^{\eset}$, $s\!\in\![\lr{\T}]$, and 
$h\!\in\!\cI_{s-1}^*$, 
we denote~by $\Ga_h(\ups_{\lr{s}};\L)$ and $\Ga_h^{0,1}(\ups_{\lr{s}};\L)$
the completions  of the spaces
\begin{gather*}
\big\{\xi\!\in\!\Ga_B(b_s(\ups);\L)\!:
\xi|_{\Si_{\ups_{\lr{s}}}-\Si_{\ups_{\lr{s}}}^h}\!=\!0\big\} 
\qquad\hbox{and}\qquad
\big\{\eta\!\in\!\Ga_B^{0,1}(b_s(\ups);\L)\!:
\eta|_{\Si_{\ups_{\lr{s}}}-\Si_{\ups_{\lr{s}}}^h}\!=\!0\big\}
\end{gather*}
with respect to the norms $\|\cdot\|_{\ups_{\lr{s}},p,1}$ and 
$\|\cdot\|_{\ups_{\lr{s}},p}$, respectively.
Let 
\begin{gather*}
\Ga_{h;-}(\ups_{\lr{s}};\L)=\Ga_-(\ups_{\lr{s}};\L)\cap \Ga_h(\ups_{\lr{s}};\L),
\qquad\hbox{where}\\
\Ga_-(\ups_{\lr{s}};\L)=
\big\{\xi\!\circ\!q_{\ups_s;\lr{\T}+1-s}\!:
\xi\!\in\!\ti\Ga_-(\ups_{\lr{s+1}};\L)\big\}.
\end{gather*}
We denote the $L^2$-orthogonal complement of $\Ga_{h;-}(\ups_{\lr{s}};\L)$
in $\Ga_h(\ups_{\lr{s}};\L)$ by $\Ga_{h;+}(\ups_{\lr{s}};\L)$.
By~\e_ref{mapgluing_e2} and the same argument as in Subsection~3.5 in~\cite{Z3},
\begin{equation}\label{bundleisom_e4}
C(b)^{-1}\|\xi\|_{\ups_{\lr{s}},p,1} \le
 \big\|\bpar_{\na,b_s(\ups)}^h\xi\big\|_{\ups_{\lr{s}},p}
 \le  C(b)\|\xi\|_{\ups_{\lr{s}},p,1}
\qquad\forall~\xi\!\in\!\Ga_{h;+}(\ups_{\lr{s}};\L)
\end{equation}
for some $C\!\in\!C({\cal U}_{\T}(\P;J);\R^+)$, provided 
$\de\!\in\!C({\cal U}_{\T}(\P;J);\R^+)$ is sufficiently small. In particular, the operator
$$\bpar_{\na,b_s(\ups)}^h\!: \Ga_{h;+}(\ups_{\lr{s}};\L)\lra \Ga_h^{0,1}(\ups_{\lr{s}};\L)$$
is an isomorphism.
On the other hand,
by the construction of the map $q_{\ups_s;\lr{\T}+1-s}$ in Subsection~2.2 of~\cite{Z3},
\begin{equation}\label{bundleisom_e5}\begin{split}
&\qquad\big\|\bpar_{\na,b_s(\ups)}\big(\xi\!\circ\!q_{\ups_s;\lr{\T}+1-s}\big)
\big\|_{\ups_{\lr{s}},p} \le C(b)|\ups|^{1/p}\|\xi\|_{\ups_{\lr{s+1}},p,1}\\
&\hbox{and}\qquad  
\bpar_{\na,b_s(\ups)}\big(\xi\!\circ\!q_{\ups_s;\lr{\T}+1-s}\big)
\big|_{\Si_{b_s(\ups)}^0}=0
\qquad\forall~\xi\!\in\!\ti\Ga_-(\ups_{\lr{s+1}};\L).
\end{split}\end{equation}
Thus, by the analogue of~\e_ref{coneiso_e2} for $\ze_{\ups,s}$,
there exist unique linear maps 
$$\ve_{\ups,s;h}\!: \ti\Ga_-(\ups_{\lr{s+1}};\L) \lra \Ga_{h;+}(\ups_{\lr{s}};\L),
\qquad h\!\in\!\cI_{s-1}^*,$$
such that 
$$\ti{R}_{\ups,s}'\xi\equiv 
\Pi_{\ups,s}\Big(\xi\!\circ\!q_{\ups_s;\lr{\T}+1-s}
+\sum_{h\in\cI_{s-1}^*}\!\!\!\!\ve_{\ups,s;h}(\xi)\Big) \in \ti\Ga_-(\ups_{\lr{s}};\L)
\qquad\forall~\xi\!\in\!\ti\Ga_-(\ups_{\lr{s+1}};\L).$$
Furthermore, for all $\xi\!\in\!\ti\Ga_-(\ups_{\lr{s+1}};\L)$ and $h\!\in\!\cI_{s-1}^*$,
\begin{equation}\label{bundleisom_e7}\begin{split}
\big\|\ve_{\ups,s;h}(\xi)\big\|_{\ups_{\lr{s}},p,1}
&\le C(b)\big( \big\|\ze_{\ups,s}|_{\Si_{\ups_{\lr{s}}}^h}
              \big\|_{\ups_{\lr{s}},p,1}^2\|\xi\|_{\ups_{\lr{s+1}},p,1}\\
&\qquad\qquad\qquad +\big\|\bpar_{\na,b_s(\ups)}\big(\xi\!\circ\!q_{\ups_s;\lr{\T}+1-s}\big)
|_{\Si_{\ups_{\lr{s}}}^h}\big\|_{\ups_{\lr{s}},p}\big).
\end{split}\end{equation}
In addition, for all $h,h'\!\in\!\cI_{s-1}^*$ such that $h'\!\neq\!h$,
$$\ve_{\ups,s;h}(\xi)=0 \qquad\forall\,\xi\!\in\!\ti\Ga_{h';-}(\ups_{\lr{s+1}};\L).$$\\

\noindent
The expansion in Lemma~\ref{derivest_lmm} below is 
a key step in constructing the homomorphisms~$\ti{R}_{\ups,s}$ with the desired properties.
For every $h\!\in\!\cI_{s-1}^*\!-\!\chi(\T)$, let
$$\chi_h'(\T)=\big\{h'\!\in\!\hat{I}\!: \io_{h'}\!=\!h\big\}.$$
If $b\!\in\!{\cal U}_{\T}^{(0)}(\P;J)$ and $h'\!\in\!\chi_h'(\T)$, we denote by
$$x_{h'}(b)\in\C=\Si_{b,h}\!-\!\{\i\}$$
the node shared by $\Si_{b,h}$ and $\Si_{b,h'}$.

\begin{lmm}
\label{derivest_lmm}
If $n$, $d$, $k$, $a$, and $\L$ are as in Proposition~\ref{g1conebdstr_prp},  
there exists $\de_n(d)\!\in\!\R^+$ such that for every almost complex structure $J$ on~$\P$,
such that $\|J\!-\!J_0\|_{C^1}\!\le\!\de_n(d)$, and a bubble type $\T$  as above,
there exist $\de,C\!\in\!C({\cal U}_{\T}(\P;J);\R^+)$  such that the requirement of
Lemma~\ref{conn_lmm} is satisfied.
Furthermore, for every $\ups\!=\!(b,v)\!\in\!\ti{\cal F}^1\T_{\de}^{\eset}$ 
and $s\!\in\![\lr{\T}]$, there exists an isomorphism
$$\ti{R}_{\ups,s}'\!: \ti\Ga_-(\ups_{\lr{s+1}};\L)
\lra  \ti\Ga_-(\ups_{\lr{s}};\L)$$
such that
\begin{gather}\label{derivestlmm_e1a}
\big\|\ti{R}_{\ups,s}'\xi-\Pi_{\ze_{\ups,s}}(\xi\!\circ\!q_{\lr{\T}+1-s})
\big\|_{\ups_{\lr{s},p,1}}\le C(b)|\ups|^{1/p}\|\xi\|_{\ups_{\lr{s+1},p,1}}
\quad\forall\, \xi\!\in\!\ti\Ga_-(\ups_{\lr{s+1}};\L);\\
\label{derivestlmm_e1b}
\ti{R}_{\ups,s}'\xi\in\ti\Ga_{h;-}(\ups_{\lr{s}};\L)
\qquad\forall\,\xi\!\in\!\ti\Ga_{h;-}(\ups_{\lr{s+1}};\L).
\end{gather}
In addition, there exist homomorphisms 
$$\ve_{\ups,h;i}\!: \ti\Ga_-(\ups_{\lr{s+1}};\L) \lra\L_{\ev_P(b)},
\qquad h\!\in\!\cI_{s-1}^*\!\cap\!\chi^-(\T),~i\!\in\!\chi_h(\T),$$
such that for all $\xi\!\in\!\ti\Ga_-(\ups_{\lr{s+1}};\L)$, 
$h\!\in\!\cI_{s-1}^*\!\cap\!\chi^-(\T)$, and $i\!\in\!\chi_h(\T)$, 
\begin{equation}\label{derivestlmm_e2}\begin{split}
&\qquad\quad \big|\ve_{\ups,h;i}(\xi)\big|\le 
C(b)|\ups|^{1/p}\|\xi\|_{\ups_{\lr{s+1}},p,1} \qquad\hbox{and}\\
&\D_{s;h}^{(1)}\big\{\ti{R}_{\ups,s}'\xi\big\}=
\sum_{h'\in\chi_h'(\T)}\!\!\!\!v_{h'}\D_{s+1;h'}^{(1)}\xi  
+\sum_{i\in\chi_h(\T)}\!\!\!\!\rho_{h;i}(\ups)\ve_{\ups,h;i}(\xi).
\end{split}\end{equation}
Furthermore, the maps $\ups\!\lra\!\ti{R}_{\ups,s}'$ and $\ups\!\lra\!\ve_{\ups,h;i}$
are $\Aut(\T)\!\propto\!(S^1)^I$-invariant and smooth over~$\ti{\cal F}^1\T_{\de}^{\eset}$.
They extend continuously over~$\ti{\cal F}^1\T_{\de}$.
These extensions satisfy
\begin{equation}\label{derivestlmm_e3}
\ti{R}_{b,s}'=\id \quad\forall\,b\!\in\!{\cal U}_{\T}^{(0)}(\P;J) \qquad\hbox{and}\qquad
\ve_{\ups,h;i}=0 \quad \forall\,h\!\in\!\cI_{s-1}^0(\ups),\,i\!\in\!\chi_h(\T).
\end{equation}\\
\end{lmm}

\noindent
Isomorphisms $\ti{R}_{\ups,s}'$ satisfying \e_ref{derivestlmm_e1a} 
and \e_ref{derivestlmm_e1b} have already been constructed.
The estimate~\e_ref{derivestlmm_e2} is obtained by applying the integration-by-parts
argument in the proof of Theorem~2.8 in~\cite{Z2} to the holomorphic functions
$$\{\ti\Pi_{\i,\cdot}^{\ups,s}\}^{-1}\ti{R}_{\ups,s}'\xi\!:
\Si_{\ups_{\lr{s}}}^{h;0}\big(\de(b)\big)\lra\L_{\ev_P(b)}
\qquad\hbox{and}\qquad
\{\ti\Pi_{\i,\cdot}^{\ups,s+1}\}^{-1}\xi\!:
\Si_{\ups_{\lr{s+1}}}^{h';0}\big(\de(b)\big)\lra\L_{\ev_P(b)}.$$
The homomorphism $\ve_{\ups,h;i}$ is given by
\begin{equation}\label{derivestlmm_e5}
\ve_{\ups,h;i}(\xi)=\frac{1}{2\pi\I}
\oint_{\partial^-A_{\ups_{\lr{s}},i}^-(\de(b))}\!
\big\{\{\Pi_{\i,\cdot}^{\ups,s}\}^{-1}\ve_{\ups,s;h}(\xi)\big\}(w_i)\frac{dw_i}{w_i^2},
\end{equation}
where $w_i$ is the coordinate on a neighborhood of the circle 
$\partial^-A_{\ups_{\lr{s}},i}^-(\de(b))$ induced from the standard
holomorphic coordinate centered at $\i$ in $\Si_{b,i}\!=\!S^2$;
see the proof of Lemma~\ref{g1comp-reg0_lmm4} in~\cite{Z5} for details.
By the continuity of the maps
$$\ups\lra \ze_{\ups,s},\,\na^{\ups,s}$$
over $\ti{\cal F}^1\T_{\de}$ and the same argument as in Subsection~4.1 of~\cite{Z3}, 
the homomorphisms $\ve_{\ups,s;h}$ extend continuously over~$\ti{\cal F}^1\T_{\de}$.
Thus, by~\e_ref{derivestlmm_e5}, the homomorphisms $\ve_{\ups,h;i}$ also extend 
continuously over~$\ti{\cal F}^1\T_{\de}$.
By~\e_ref{mapgluing_e4b}, \e_ref{connlmm_e2}, and \e_ref{bundleisom_e7},
$$\ve_{\ups,s;h}=0 \qquad \forall\,h\!\in\!\cI_{s-1}^0(\ups).$$
This observation, along with~\e_ref{derivestlmm_e5}, implies the second claim
in~\e_ref{derivestlmm_e3}.
The first claim in~\e_ref{derivestlmm_e3} follows from~\e_ref{mapgluing_e4a},
\e_ref{connlmm_e1}, and the construction of~$\ti{R}_{\ups,s}'$.\\

\noindent
Suppose $\ups\!=\!(b,v)\!\in\!\ti{\cal F}^1\T_{\de}^{\eset}$,
$s\!\in\![\lr{\T}]$, and we have constructed an isomorphism 
$$\ti{R}_{\ups,s+1}\!:\Ga_-(b;\L) \lra \ti\Ga_-\big(\ups_{\lr{s+1}};\L\big)$$
that satisfies \e_ref{bundleisom_e0}-\e_ref{bundleisom_e2}.
We note that for every $s\!\in\![\lr{\T}]$ and 
$h\!\in\!\cI_{s-1}^*\!\cap\!\chi^-(\T)$,
$$\rho_{s;h}(\ups)=\big(v_{h'}\rho_{s+1;h'}(\ups)\big)_{\io_{h'}=h}
\qquad\forall~ 
\ups\!\equiv\!\big(b,(v_i)_{i\in\aleph\cup\hat{I}}\big)\in\ti{\cal F}\T.$$
Thus, by \e_ref{bundleisom_e0} and  \e_ref{bundleisom_e1} with $s$ replaced by $s\!+\!1$,
~\e_ref{derivestlmm_e1b}, and~\e_ref{derivestlmm_e2},  there exists a homomorphism
$$\ti\ve_{\ups,s;h}\!: \ti\F_h\T \lra \Hom\big(\Ga_-(b;\L),\L_{\ev_P(b)}\big)$$
such that
\begin{gather}\label{derivestlmm_e8a}
\big|\ti\ve_{\ups,s;h}\big|\le C(b)|\ups|^{1/p},\quad
\ti\ve_{\ups,s;h}(w_h;\xi)=0 ~~~\forall\,w_h\!\in\!\ti\F_h\T,\,
\xi\!\in\!\Ga_{h';-}(b;\L),\, h'\!\in\!\cI_{s-1}^*\!-\!\{h\},\\
\label{derivestlmm_e8b}
\D_{s;h}^{(1)}\big\{\ti{R}_{\ups,s}'\ti{R}_{\ups,s+1}\xi\big\} 
=\al_{s;h}\big(\rho_{s;h}(\ups);\xi\big)+\ti\ve_{\ups,s;h}\big(\rho_{s;h}(\ups);\xi\big)
\qquad\forall\,\xi\!\in\!\Ga_-(b;\L).
\end{gather}
We note that for $h\!\in\!\cI_{s-1}^*\!-\!\chi^-(\T)$, the existence of 
such $\ti\ve_{\ups,s;h}$ is immediate from \e_ref{bundleisom_e1} with $s$ 
replaced by $s\!+\!1$, \e_ref{derivestlmm_e1a}, and~\e_ref{derivestlmm_e1b}.
Let $[\rho_{s;h}(\ups)]$ denote the image of $\rho_{s;h}(\ups)$ under the projection map
$\ti\F_h\T\!-\!\{0\}\!\lra\!\bP\ti\F_h\T$.
Since 
$$\al_{s;h}\big(w_h;\cdot\big)\!: 
\Ga_{h;-}^{\perp}(b;\L;[\rho_{s;h}(\ups)]) \lra\L_{\ev_P(b)}$$
is an isomorphism for each $h\!\in\!\cI_{s-1}^*$, 
by the first bound in~\e_ref{derivestlmm_e8a}, \e_ref{derivestlmm_e8b}, 
and~\e_ref{bundleisom_e3} there exists a unique homomorphism
$$\mu_{\ups,s;h}\!:   \Ga_-(b;\L)\lra \Ga_{h;-}^{\perp}(b;\L;[\rho_{s;h}(\ups)]),$$
such that
\begin{equation}\label{derivestlmm_e9a}
\D_{s;h}^{(1)}\big\{\ti{R}_{\ups,s}'\ti{R}_{\ups,s+1}
                       \big(\xi\!+\!\mu_{\ups,s;h}(\xi)\big)\big\} 
=\al_{s;h}\big(\rho_{s;h}(\ups);\xi\big).
\end{equation}
Furthermore, by~\e_ref{bundleisom_e3} and~\e_ref{derivestlmm_e8a},
\begin{equation}\label{derivestlmm_e9}
\big|\mu_{\ups,s;h}\big|\le  C(b)|\ups|^{1/p}, \qquad
\mu_{\ups,s;h}(\xi)=0 ~~~\forall\, \xi\!\in\!\Ga_{h';-}(b;\L),\,
h'\!\in\!\cI_{s-1}^*\!-\!\{h\}.
\end{equation}
We define
$$\ti{R}_{\ups,s}\!:\Ga_-(b;\L) \lra \ti\Ga_-\big(\ups_{\lr{s}};\L\big)
\qquad\hbox{by}\qquad
\ti{R}_{\ups,s}(\xi)=  \ti{R}_{\ups,s}'\ti{R}_{\ups,s+1}
  \Big(\xi+\sum_{h\in\cI_{s-1}^*}\!\!\!\mu_{\ups,s;h}(\xi)\Big).$$
By \e_ref{bundleisom_e0} with $s$ replaced by $s\!+\!1$, \e_ref{derivestlmm_e1b}, 
and the second statement in~\e_ref{derivestlmm_e9}, $\ti{R}_{\ups,s}$
satisfies~\e_ref{bundleisom_e0}.
Since 
$$\D_{s;h}^{(1)}\xi=0 \qquad\forall\, \xi\!\in\!\ti\Ga_{h';-}\big(\ups_{\lr{s}};\L\big),
\, h'\!\in\!\cI_{s-1}^*\!-\!\{h\},$$
$\ti{R}_{\ups,s}$ satisfies~\e_ref{bundleisom_e1} by~\e_ref{derivestlmm_e9a},
along with \e_ref{bundleisom_e0} with $s$ replaced by $s\!+\!1$, \e_ref{derivestlmm_e1b}, 
and the second statement in~\e_ref{derivestlmm_e9}.\\

\noindent
It remains to show that for every $h\!\in\!\cI_{s-1}^*$ the family of homomorphisms
$$\mu_{\ups,s;h}\!: \Ga_-(b;\L)\lra \Ga_{h;-}^{\perp}(b;\L;[\rho_{s;h}(\ups)])
\subset \Ga_{h;-}(b;\L), \qquad \ups\!\in\!\ti{\cal F}^1\T_{\de}^{\eset},$$
extends continuously over $\ti{\cal F}^1\T_{\de}$.
Each homomorphism~$\ti\ve_{\ups,s;h}$ of the previous paragraph extends 
continuously over $\ti{\cal F}^1\T_{\de}$, as this is case for the homomorphisms
$\ve_{\ups,h;i}$ by Lemma~\ref{derivest_lmm}.
Furthermore,
\begin{equation}\label{derivestlmm_e12}
\ti\ve_{b,s;h}=0 ~~~\forall\,
b\!\in\!{\cal U}_{\T}^{(0)}(\P;J),\, h\!\in\!\cI_{s-1}^*,
\quad\hbox{and}\quad
\ve_{\ups,s;h}=0  ~~~\forall\,
\ups\!\in\!\ti{\cal F}^1\T_{\de},\,h\!\in\!\cI_{s-1}^0(\ups).
\end{equation}
The first claim above follows from 
\e_ref{bundleisom_e2} with $s$ replaced by $s\!+\!1$ and 
first statement in~\e_ref{derivestlmm_e3}.
The second claim in~\e_ref{derivestlmm_e12} follows from the second 
statement in~\e_ref{derivestlmm_e3}.
If $\ups\!\in\!\ti{\cal F}^1\T_{\de}$ and $h\!\in\!\cI_{s-1}^*\!-\!\cI_{s-1}^0(\ups)$,
we define $\mu_{\ups,s;h}$ as in the previous paragraph.
This extension is continuous at~$\ups$ since $\ti\ve_{\ups,s;h}$~is.
If $h\!\in\!\cI_{s-1}^0(\ups)$, we take $\mu_{\ups,s;h}\!=\!0$.
This extension is continuous by the continuity of $\ti\ve_{\ups,s;h}$
{\it and} the second statement in~\e_ref{derivestlmm_e12}.
Finally, $\ti{R}_{\ups,s}$ satisfies~\e_ref{bundleisom_e2} by
the first statement in~\e_ref{derivestlmm_e12}, along with 
\e_ref{bundleisom_e2} with $s$ replaced by $s\!+\!1$ and the first statement
in~\e_ref{derivestlmm_e3}.\\

\noindent
{\it Remark:} The key point in the previous paragraph is the second statement 
in~\e_ref{derivestlmm_e12}, because the lines $\Ga_{h;-}^{\perp}(b;\L;[\rho_{s;h}(\ups)])$
may {\it not} extend continuously over~$\ti{\cal F}^1\T_{\de}$. 

\begin{crl}
\label{gluing0_crl}
If $n$, $d$, $k$, $a$, and $\L$ are as in Proposition~\ref{g1conebdstr_prp},  
there exists $\de_n(d)\!\in\!\R^+$ such that for every almost complex structure $J$ on~$\P$,
such that $\|J\!-\!J_0\|_{C^1}\!\le\!\de_n(d)$, and a bubble type $\T$  as above,
there exist $\de,C\!\in\!C({\cal U}_{\T}(\P;J);\R^+)$  such that the requirement of
Lemma~\ref{conn_lmm} is satisfied.
In addition, for every $\ups\!=\!(b,v)\!\in\!\ti{\cal F}^1\T_{\de}^{\eset}$
there exists an isomorphism
$$\ti{R}_{\ups,1}\!: \Ga_-(b;\L)\lra  \ti\Ga_-(\ups_{\lr{1}};\L)$$
such that for every $\xi\!\in\!\Ga_-(b;\L)$, $h\!\in\!\cI_0^*$,
and $\ep\!\in\!(0,2\de(b))$,
\begin{gather}\label{gluing0crl_e2a}
\big\|\na^{\ups,1}\ti{R}_{\ups,1}\xi
 \big\|_{C^0(A_{\ups_{\lr{1}},h}^-(\de(b))),g_{\ups_{\lr{1}}}}
\le C(b)\big|\rho_{1;h}(\ups)\big|\cdot\|\xi\|_{b,p,1}, \qquad\hbox{and}\\
\label{gluing0crl_e2b}
\oint_{\partial^-A_{\ups_{\lr{1}},h}^-(\ep)}\!
\big\{\{\ti\Pi_{\i,\cdot}^{\ups,1}\}^{-1}\ti{R}_{\ups,1}\xi\big\}(w_h) \frac{dw_h}{w_h^2}
=2\pi\I \sum_{i\in\chi_h(\T)}\!\!\rho_{h;i}(\ups)\D_{\T,i}\xi,
\end{gather}
where $w_h$ is the standard holomorphic on the neighborhood of $\i$ in 
$\Si_{\ups_{\lr{1}},h}\!=\!S^2$.
Finally, the map $\ups\!\lra\!\ti{R}_{\ups,1}$ is 
$\Aut(\T)\!\propto\!(S^1)^I$-invariant and smooth on $\ti{\cal F}^1\T_{\de}^{\eset}$.
It extends continuously over~$\ti{\cal F}^1\T_{\de}$.
This extension satisfies
\begin{equation}\label{gluing0crl_e3}
\ti{R}_{b,1}=\id \quad\forall\,b\!\in\!{\cal U}_{\T}^{(0)}(\P;J).
\end{equation}\\
\end{crl}

\noindent
The homomorphism $\ti{R}_{\ups,1}$ constructed above  satisfies 
the extension requirements of the corollary.
Since $\{\ti\Pi_{\i,\cdot}^{\ups,1}\}^{-1}\ti{R}_{\ups,1}\xi$
is holomorphic on $A_{\ups_{\lr{1}},h}^-(\ep)$,
\e_ref{gluing0crl_e2b} is equivalent to the $s\!=\!1$ case of~\e_ref{bundleisom_e1}.\\

\noindent
It remains to verify \e_ref{gluing0crl_e2a}.
Let 
$$u_{\ups_{\lr{1}}}=u_b\circ q_{\ups_{\lr{1}}}.$$
For each $h\!\in\!\cI_0$ and $z\!\in\!\Si_{\ups_{\lr{1}}}^{h;0}(\de(b))$, 
we denote by $\Pi_{\i,z}^{\ups_{\lr{1}}}$ the parallel transport in the line bundle 
$u_{\ups_{\lr{1}}}^*\L$ along a path from $\i$ to $z$ in 
$\Si_{\ups_{\lr{s}}}^{h;0}(\de(b))$  with respect to 
the connection $q_{\ups_{\lr{1}}}^*\na^{\ups,\lr{\T}+1}$.
By the construction of the homomorphism $\ti{R}_{\ups,1}$ above,
$$\{\ti\Pi_{\i,\cdot}^{\ups,1}\}^{-1}\ti{R}_{\ups,1}\xi
  \big|_{\Si_{\ups_{\lr{1}}}^{h;0}(\de(b))}  =
\{\Pi_{\i,\cdot}^{\ups_{\lr{1}}}\}^{-1}(\xi\!\circ\!q_{\ups_{\lr{1}}})
  \big|_{\Si_{\ups_{\lr{1}}}^{h;0}(\de(b))}
+\ve_{\ups}(\xi)\qquad\forall\,\xi\!\in\!\Ga_-(b;\L),$$
for some homomorphism 
\begin{gather*}
\ve_{\ups}\!:\Ga_-(b;\L)\lra C^{\i}\big(\Si_{\ups_{\lr{1}}}^{h;0}(\de(b));\L_{\ev_P(b)}\big)
\qquad\hbox{s.t.}\\
\big\|\ve_{\ups}(\xi)\big\|_{C^0(\Si_{\ups_{\lr{s}}}^{h;0}(\de(b)))}
\le C(b)|\ups|^{1/p}\|\xi\|_{b,p,1} \qquad\forall\,\xi\!\in\!\Ga_-(b;\L).
\end{gather*}
Thus, by the same integration-by-parts argument as in the proof of Theorem~2.2
in~\cite{Z2}, there exist homomorphisms 
$$\ve_{\ups,h;i}^{(l)}\!: \Ga_-(b;\L)\lra\L_{\ev_P(b)},
\qquad\forall~h\!\in\!\cI_0^*,\,i\!\in\!\chi_h(\T),\,l\!\in\!\Z^+,$$
such that for all $h\!\in\!\cI_0^*$, $i\!\in\!\chi_h(\T)$, $l,m\!\in\!\Z^+$,
and $\xi\!\in\!\Ga_-(b;\L)$
\begin{equation}\label{gluing0crl_e5}\begin{split}
&\qquad\qquad\qquad\qquad
\big|\ve_{\ups,h;i}^{(l)}(\xi)\big|\le C(b)\de(b)^{-l/2}\|\xi\|_{b,p,1},\\
&\D_{s;h}^{(m)}\big\{\ti{R}_{\ups,1}\xi\big\}=
\sum_{l=1}^{l=m}\binom{m\!-\!1}{l\!-\!1} \sum_{i\in\chi_h(\T)}\!\!\!\!
x_i^{m-l}(\ups)\rho_{h;i}^l(\ups)\big(\D_{\lr{\T}+1;i}^{(l)}\xi  
+\ve_{\ups,h;i}^{(l)}(\xi)\big).
\end{split}\end{equation}
The number $x_i(\ups)\!\in\!\C$ is given explicitly in the paragraph preceding
Lemma~\ref{g1comp-reg0_lmm3} in~\cite{Z5}.
It is close~to 
$$x_{h'}(b)\in\C=\Si_{b,h}\!-\!\{\i\}, \qquad\hbox{where}\qquad
h'\!\le\!i,\,\io_{h'}\!=\!h.$$
The estimate~\e_ref{gluing0crl_e2a} is obtained by summing up the derivatives of
$\ti{R}_{\ups,1}\xi|_{\Si_{\ups_{\lr{1}}}^h}$ at~$\i$ with the appropriate coefficients,
using~\e_ref{gluing0crl_e5}; see the proof of Lemma~4.2 in~\cite{Z4}
for a similar argument.

\subsection{Smoothing Bundle Sections, II}
\label{gluing_subs3}

\noindent
In this subsection, we take the inductive construction of the previous subsection
one step further to define a homomorphism $\ti{R}_{\ups}\!\equiv\!\ti{R}_{\ups,0}$.
However, in this case we will encounter an obstruction bundle.
The homomorphism $\ti{R}_{\ups}$ will not extend continuously over~$\ti{\cal F}^1\T$,
but its restriction to a cone contained in $\ti\V_{1,k}^d$~will.\\

\noindent
We first recall certain facts concerning the modified gluing map
$$\ti{q}_{\ups_0;\lr{\T}+1}\!:\Si_{\ups}\lra\Si_{\ups_{\lr{1}}}$$
corresponding to the parameter $\de(b)^{1/2}$, as 
constructed in Subsection~\ref{g1comp-reg1_subs2} of~\cite{Z5}.
Suppose
$$\ups\equiv(b,v_{\aleph},(v_h)_{h\in\hat{I}}\big)\in\ti{\cal F}^1\T_{\de}^{\eset}.$$
The map $\ti{q}_{\ups_0;\lr{\T}+1}$ is biholomorphic outside $|\aleph|$
thin necks $A_{\ups,h}$, with $h\!\in\!\aleph$, of $(\Si_{\ups},g_{\ups})$
and the $|I_1|$ annuli 
$$\ti\A_{b,h} \equiv \ti\A_{b,h}^-\cup\ti\A_{b,h}^+,$$
with $h\!\in\!I_1$, where
$$\ti\A_{b,h}^{\pm} \equiv \ti\A_{b,h}^{\pm}\big(\de(b)\big)
\subset \Si_{b;P} \!\approx\! \Si_{\ups}$$
are annuli independent of~$\ups$.
In addition,
\begin{gather}\label{glmap_e1a}\begin{split}
&\qquad\qquad\qquad\qquad 
\ti{u}_{\ups,1}|_{\ti{q}_{\ups_0;\lr{\T}+1}(A_{\ups,h})}=\const 
~~\forall\,h\!\in\!\aleph, \\
&\ti{u}_{\ups,1}|_{\ti{q}_{\ups_0;\lr{\T}+1}(\ti\A_{b,h})}=\const 
~~\forall\,h\!\in\!I_1\!-\!\cI_0, \qquad
\ti{u}_{\ups,1}|_{\ti{q}_{\ups_0;\lr{\T}+1}(\ti\A_{b,h}^+)}=\const 
~~\forall\,h\!\in\!\cI_0;
\end{split}\\
\label{glmap_e1b}
\ti{q}_{\ups_0;\lr{\T}+1}(\ti\A_{b,h}^-) \subset 
  A_{\ups_{\lr{1}},h}^-\big(|v_h|^2/\de(b)\big) 
  ~~~\hbox{and}~~~
\big\|d\ti{q}_{\ups_0;\lr{\T}+1}\big\|_{C^0(\ti\A_{b,h}^-)}
\le C(b)|v_h| \quad\forall~h\!\in\!\cI_0,
\end{gather}
if the $C^0$-norm of $d\ti{q}_{\ups_0;\lr{\T}+1}$ is computed
 with respect to the metrics $g_{\ups}$ on $\Si_{\ups}$
and $g_{\ups_{\lr{1}}}$ on~$\Si_{\ups_{\lr{1}}}$.
Furthermore, 
\begin{equation}\label{glmap_e2}
\big\|d\ti{q}_{\ups_0;\lr{\T}+1}\big\|_{C^0} \le C(b).
\end{equation}\\

\noindent
We now proceed similarly to the previous subsection.
If $\ups\!\in\!\ti{\cal F}^1\T_{\de}^{\eset}$, we denote 
the completions of the spaces
$$\Ga\big(\Si_{\ups};u_{\ups,0}^*\L\big) \qquad\hbox{and}\qquad
 \Ga\big(\Si_{\ups};\La^{0,1}_{\I,j}T^*\Si_{\ups}\!\otimes\!u_{\ups,0}^*\L\big)$$
with respect to the Sobolev norms $\|\cdot\|_{\ups,p,1}$ and $\|\cdot\|_{\ups,p}$
by $\Ga(\ups;\L)$ and $\Ga^{0,1}(\ups;\L)$.
Let
$$\Ga_-(\ups;\L)=\big\{R_{\ups,0}'\xi\!:\xi\!\in\!\ti\Ga_-(\ups_{\lr{1}};\L)\big\},
\qquad\hbox{where}\quad
R_{\ups,0}'\xi=\xi\!\circ\!\ti{q}_{\ups_0;\lr{\T}+1}.$$
Let $R_{\ups,0}\!=\!R_{\ups,0}'\ti{R}_{\ups,1}$.
By \e_ref{gluing0crl_e2a}, \e_ref{glmap_e1a}, and~\e_ref{glmap_e1b},
\begin{equation}\label{glmap_e4}\begin{split}
\big\|\bpar_{\na,b_0(\ups)}R_{\ups,0}'\xi\big\|_{\ups,p}
&\le C(b)|\rho(\ups)|\|\xi\|_{\ups_{\lr{1}},p,1}\\
&\le C'(b)|\rho(\ups)|\|R_{\ups,0}'\xi\|_{\ups,p,1}
\qquad \forall\,\xi\!\in\!\Ga_-(\ups_{\lr{1}};\L).
\end{split}\end{equation}
Let $\Ga_+(\ups;\L)$ denote the $L^2$-orthogonal complement of 
$\Ga_-(\ups;\L)$ in $\Ga(\ups;\L)$.
Similarly to~\e_ref{bundleisom_e4},
\begin{equation}\label{glmap_e5}
C(b)^{-1}\|\xi\|_{\ups,p,1} \le
 \big\|\bpar_{\na,b_0(\ups)}\xi\big\|_{\ups,p} \le  C(b)\|\xi\|_{\ups,p,1}
\qquad\forall~\xi\!\in\!\Ga_+(\ups;\L)
\end{equation}
for some $C\!\in\!C({\cal U}_{\T}(\P;J);\R^+)$, provided 
$\de\!\in\!C({\cal U}_{\T}(\P;J);\R^+)$ is sufficiently small.
Let $\Ga_+^{0,1}(\ups;\L)$ be the image of $\Ga_+(\ups;\L)$ under~$\bpar_{\na,b_0(\ups)}$.\\

\noindent
In contrast to the previous subsection, the operator $\bpar_{\na,b_0(\ups)}$
is not surjective.
We next describe a complement of $\Ga_+^{0,1}(\ups;\L)$ in~$\Ga^{0,1}(\ups;\L)$. 
Since the operator $\bpar_{\na,b}^B$ is surjective,
the cokernel of $\bpar_{\na,b}$
can be identified with the vector space
$$\Ga^{0,1}_-(b;\L) \equiv \H_{b;P}\!\otimes\!\L_{\ev_P(b)}
\approx\E_{\pi_P(b)}^*\!\otimes\!\L_{\ev_P(b)},$$
where $\H_{b;P}$ is the space of harmonic antilinear differentials on 
the main component $\Si_{b;P}$ of~$\Si_b$.
If $\aleph\!\neq\!\eset$, i.e.~$\Si_{b;P}$ is a circle of spheres,
the elements of $\H_{b;P}$ have simple poles at the nodes of $\Si_{b;P}$
with the residues adding up to zero at each node.
Since the Riemann surfaces $\Si_{\ups}$, with $\ups\!\in\!{\cal F}^1\T_{\de}$,
are deformations of $\Si_b$, with $b\!\in\!{\cal U}_{\T}^{(0)}(\P;J)$, 
there exists a family of isomorphisms
$$R^{0,1}_{\ups;P}\!: \H_{b;P} \lra  \H_{\ups;P}\!\equiv\!\H_{b_0(\ups);P},
\qquad \ups\!=\!(b,v)\in {\cal F}^1\T_{\de},$$
such that the family of induced homomorphisms
$$ \H_{b;P}\lra\Ga^{0,1}(\ups;\C)^*,
\qquad \{R^{0,1}_{\ups;P}\eta\}(\eta')=\llrr{R^{0,1}_{\ups;P}\eta,\eta'}_2
\quad\forall\,\eta\!\in\!\H_{b;P},\,\eta'\!\in\!\Ga^{0,1}(\ups;\C),$$
is $\Aut(\T)\!\propto\!(S^1)^I$-invariant and smooth on $\ti{\cal F}^1\T_{\de}^{\eset}$,
continuous on~$\ti{\cal F}^1\T_{\de}$, and
\begin{equation}\label{glmap_e6}
R^{0,1}_{\ups;P}|_b=\id \qquad\forall\, b\!\in\!{\cal U}_{\T}^{(0)}(\P;J).
\end{equation}\\

\noindent
With notation as in~\e_ref{cutoff_e}, we define $\ti\be_b\!\in\!C^{\i}(\Si_b;\R)$ by
$$\ti\be_b(z)=
\begin{cases}
1,&\hbox{if}~z\!\in\!\Si_{b,i},~i\!\in\!\chi^0(\T);\\
1-\be_{\de(b)}(r(z)),&\hbox{if}~z\!\in\!\Si_{b,i},~i\!\in\!\chi(\T);\\
0,&\hbox{otherwise}.
\end{cases}$$
In other words, $\ti\be_b\!=\!1$ on $\Si_b^0(\de(b)/2)$ and vanishes
outside of $\Si_b^0(2\de(b))\!\subset\!\Si_b$.
Let $\ti\be_{\ups}\!=\!\ti\be_b\!\circ\!q_{\ups}$.
If $z\!\in\!\Si_{\ups}^0(2\de(b))$,
we denote by $\Pi_z^{\ups,0}$ the parallel transport in the line bundle 
$u_{\ups,0}^*\L$ along a path from 
$x\!\in\!\ti{q}_{\ups_0;\lr{\T}+1}^{-1}(\Si_{\ups_{\lr{1}};P})$ 
to $z$ in $\Si_{\ups}^0(2\de(b))$  with respect to 
the connection $\ti{q}_{\ups_0;\lr{\T}+1}^*\na^{\ups,1}$.
For each 
\begin{equation}\label{upseta_e}
\ups=(b,v)\!\in\!\ti{\cal F}^1\T_{\de}^{\eset} \qquad\hbox{and}\qquad
\eta\!\in\!\Ga^{0,1}_-(b;\L),
\end{equation}
let $R_{\ups}^{0,1}\eta\!\in\!\Ga^{0,1}(\ups;\L)$ be given by
$$\{R_{\ups}^{0,1}\eta\}_zw
= \ti\be_{\ups}(z) \, \Pi_z^{\ups,0}\eta_z(w) \in \L_{u_{\ups,0}(z)}
\qquad\,z\!\in\!\Si_{\ups},\,w\!\in\!T_z\Si_{\ups}.$$
The image of $\Ga^{0,1}_-(b;\L)$ in $\Ga^{0,1}(\ups;\L)$
is a complement of $\Ga^{0,1}_+(\ups;\L)$ in~$\Ga^{0,1}(\ups;\L)$,
as can be seen from Lemma~\ref{gluing1_lmm} below.\\

\noindent
If $\eta\!\in\!\Ga_-^{0,1}(b;\L)$,  we put
$$\|\eta\|=\sum_{h\in\cI_0}|\eta|_{x_h(b)},$$
where $|\eta|_{x_h(b)}$ is the norm of $\eta|_{x_h(b)}$ with respect to
the metric $g_{\pi_P(b)}$ on $\Si_{b;P}$.
If $\ups$ and $\eta$ are as in~\e_ref{upseta_e} and $\|\eta\|\!=\!1$,
we define~by
$$\pi^{0,1}_{\ups;-}\!: \Ga^{0,1}(\ups;\L)\lra\Ga_-^{0,1}(b;\L)
\qquad\hbox{by}\qquad
\pi^{0,1}_{\ups;-}(\eta')=\llrr{\eta',R_{\ups}^{0,1}\eta}_2\eta
\quad\forall\,\eta'\!\in\!\Ga^{0,1}(\ups;\L).$$
Since the space $\Ga_-^{0,1}(b;\L)$ is one-dimensional, 
$\pi^{0,1}_{\ups;-}$ is independent of the choice of~$\eta$.
We note that since $p\!>\!2$, by Holder's inequality
\begin{equation}\label{gluing1lmm_e0}
\big\|\pi^{0,1}_{\ups;-}\eta'\big\|\le C(b)\|\eta'\|_{\ups,p}
\qquad\forall\,\eta'\!\in\!\Ga^{0,1}(\ups;\L).
\end{equation}

\begin{lmm}
\label{gluing1_lmm}
If $n$, $d$, $k$, $a$, and $\L$ are as in Proposition~\ref{g1conebdstr_prp},  
there exists $\de_n(d)\!\in\!\R^+$ such that for every almost complex structure $J$ on~$\P$,
such that $\|J\!-\!J_0\|_{C^1}\!\le\!\de_n(d)$, and a bubble type $\T$  as above,
there exist $\de,C\!\in\!C({\cal U}_{\T}(\P;J);\R^+)$  such that the requirements of
of Corollary~\ref{gluing0_crl} are satisfied.
Furthermore, with notation as above,  for all 
$\ups\!=\!(b,v)\!\in\!\ti{\cal F}^1\T_{\de}^{\eset}$,
\begin{gather}
\label{gluing1lmm_e1}
\pi^{0,1}_{\ups;-}\bpar_{\na,b_0(\ups)}R_{\ups,0}\xi
=-2\pi\I\, \D_{\T}\big(\xi\!\otimes\!\rho(\ups)\big)
\qquad\forall\,\xi\!\in\!\Ga_-(b;\L);\\
\label{gluing1lmm_e2}
\big\|\pi^{0,1}_{\ups;-}\bpar_{\na,b_0(\ups)}\xi\big\|
\le C(b)\big|\rho(\ups)\big|\|\xi\|_{\ups,p,1}
\qquad\forall\,\xi\!\in\!\Ga(\ups;\L).
\end{gather}
Finally, the map $\ups\!\lra\!\pi^{0,1}_{\ups;-}$ is $\Aut(\T)\!\propto\!(S^1)^I$-invariant 
and smooth on $\ti{\cal F}^1\T_{\de}^{\eset}$.
It extends continuously over~$\ti{\cal F}^1\T_{\de}$.
\end{lmm}

\noindent
The identity~\e_ref{gluing1lmm_e1} requires the restriction on the homomorphisms
$R_{\ups;P}^{0,1}$ and identification of gluing parameters
described in Subsection~\ref{g1comp-reg1_subs2} of~\cite{Z4}.
It follows from~\e_ref{gluing0crl_e2b} by the same integration-by-parts argument 
as used in the proof of Proposition~4.4 in~\cite{Z2}.
The estimate~\e_ref{gluing1lmm_e2} is obtained by computing 
$\bpar_{\na,b_0(\ups)}^*R_{\ups}^{0,1}\eta$;
see the proof of Lemma~2.2 in~\cite{Z2}.\\

\noindent
With notation as in the two previous subsections, let
$$\Pi_{\ze_{\ups,0}}\!: u_{\ups,0}^*\L \lra \ti{u}_{\ups,0}^*\L$$
be the $\na$-parallel transport along the geodesics 
$\tau\!\lra\!\exp_{u_{\ups,0}(z)}\!\ze_{\ups,0}(z)$, with $\tau\!\in\![0,1]$.
We~put
\begin{gather*}
L_{\ups,0}= \Pi_{\ze_{\ups,0}}^{-1}\circ \bpar_{\na,\ti{b}_0(\ups)}\circ \Pi_{\ze_{\ups,0}}
- \bpar_{\na,b_0(\ups)}: \Ga\big(\ups;\L\big)\lra\Ga^{0,1}\big(\ups;\L\big);\\
\ti\Ga_-'(\ups;\L)=
\big\{\Pi_{\ze_{\ups,0}}^{-1}\xi\!:\xi\!\in\!\ti\Ga_-(\ups;\L)\big\}
\subset\Ga\big(\ups;\L\big).
\end{gather*}
We denote by 
$$\pi_{\ups;-}\!:\Ga(\ups;\L)\lra\Ga_-(\ups;\L) \qquad\hbox{and}\qquad
\ti\pi_{\ups;-}\!:\Ga(\ups;\L)\lra\ti\Ga_-'(\ups;\L)$$
the $L^2$-projection maps.
Let $\Ga_-'(\ups;\L)$ be the image of $\tilde{\Ga}_-'(\ups;\L)$ under~$\pi_{\ups;-}$.
By the analogue of~\e_ref{coneiso_e2} for $\ze_{\ups,0}$ and~\e_ref{mapgluing_e3},
\begin{equation}\label{glmap_e8}
\big\|L_{\ups,0}\xi\big\|_{\ups,p}\le C(b)\big|\rho(\ups)\big|^2\|\xi\|_{\ups,p,1}
\qquad\forall\,\xi\!\in\!\Ga(\ups;\L).
\end{equation}
By \e_ref{glmap_e4}, \e_ref{glmap_e5}, and~\e_ref{glmap_e8},
\begin{equation}\label{g1map_e9a}
\big\|\xi-\pi_{\ups;-}\xi\big\|_{\ups,p,1}
\le C(b)\!\big|\rho(\ups)\big|\|\xi\|_{\ups,p,1}
\qquad\forall~\xi\!\in\!\ti\Ga_-'(\ups;\L).
\end{equation}
By~\e_ref{gluing1lmm_e0}-\e_ref{g1map_e9a},
\begin{equation}\label{g1map_e9b}
\big|\D_{\T}(\xi\!\otimes\!\rho(\ups))\big|
\le C(b)\big|\rho(\ups)\big|^2\|R_{\ups,0}\xi\|_{\ups,p,1}
\qquad\forall~R_{\ups,0}\xi\!\in\!\Ga_-'(\ups;\L).
\end{equation}\\

\noindent
For each $b\!\in\!{\cal U}_{\T}^{(0)}(\P;J)$ and $[w]\!\in\!\bP\ti\F\T|_b$, let
$$\Ga_-\big(b;\L;[w]\big)=\big\{\xi\!\in\!\Ga_-(b;\L);~ 
\D_{\T}(\xi\!\otimes\!w)\!=\!0\big\}.$$
Similarly to the previous subsection, the map $\D_{\T}$ is surjective.
Thus, the $L^2$-orthogonal complement $\Ga_-^{\perp}(b;\L;[w])$
of $\Ga_-(b;\L;[w])$ in~$\Ga_-(b;\L)$ is one-dimensional.
Furthermore, there exists $C\!\in\!C({\cal U}_{\T}(\P;J);\R^+)$ such~that
\begin{equation}\label{gluing1_lmm2e1}
C(b)^{-1}|w|\cdot\|\xi\|_{b,p,1} \le \big|\D_{\T}(\xi\!\otimes\!w)\big|
\le C(b)|w|\cdot\|\xi\|_{b,p,1}
\qquad\forall~\xi\!\in\!\Ga_-^{\perp}\big(b;\L;[w]\big).
\end{equation}
If $\ups\!\in\!\ti{\cal F}^1\T^{\eset}_{\de}$, let 
$$\Ga_-\big(\ups;\L;[w]\big)=
          \big\{R_{\ups,0}\xi\!: \xi\!\in\!\Ga_-(b;\L;[w])\big\}
\subset\Ga_-(\ups;\L).$$
We denote by $\Ga_-^{\perp}(\ups;\L;[w])$ the $L^2$-orthogonal complement
of $\Ga_-(\ups;\L;[w])$ in~$\Ga_-(\ups;\L)$.
Since $R_{\ups,0}$ is close to an isometry on $\Ga_-(b;\L)$
with respect to the $L^2$ and $L^p_1$-norms,
\begin{equation}\label{gluing1_lmm2e3}
\big|\D_{\T}(\xi\!\otimes\!w\big|
\ge C(b)^{-1}|w|\|R_{\ups,0}\xi\|_{\ups,p,1}
\qquad\forall~R_{\ups}\xi\!\in\!\Ga_-^{\perp}\big(\ups;\L;[\rho(\ups)]\big),
\end{equation}
by~\e_ref{gluing1_lmm2e1}.
We note that
$$\dim\Ga_-\big(\ups;\L;[\rho(\ups)]\big)
              =\dim\ti\Ga_-'(\ups;\L)=\dim\Ga_-'(\ups;\L).$$
Thus, by \e_ref{g1map_e9a}, \e_ref{g1map_e9b}, and \e_ref{gluing1_lmm2e3}
applied with $w\!=\!\rho(\ups)$, the~map
$$\ti\pi_{\ups;-}\!:\Ga_-\big(\ups;\L;[\rho(\ups)]\big)\lra\ti\Ga_-'(\ups;\L)$$
is an isomorphism. Furthermore,
\begin{equation}\label{g1map_e10}
\big\|\xi-\ti\pi_{\ups;-}\xi\big\|_{\ups,p,1}
\le C(b)\big|\rho(\ups)\big|\|\xi\|_{\ups,p,1}
\qquad\forall~\xi\!\in\!\Ga_-\big(\ups;\L;[\rho(\ups)]\big).
\end{equation}\\

\noindent
If $b\!\in\!{\cal U}_{\T}^{(0)}(\P;J)$, let
$$\ti\Ga_-(b;\L)=\big\{\xi\!\in\!\Ga_-(b;\L)\!: 
\D_{\T}(\xi\!\otimes\!w)\!=\!0~\forall\,w\!\in\!\ti\F^1\T|_b\big\}.$$

\begin{crl}
\label{gluing1_crl1}
If $n$, $d$, $k$, $a$, and $\L$ are as in Proposition~\ref{g1conebdstr_prp},  
there exists $\de_n(d)\!\in\!\R^+$ such that for every almost complex structure $J$ on~$\P$,
such that $\|J\!-\!J_0\|_{C^1}\!\le\!\de_n(d)$, and a bubble type $\T$  as above,
there exist $\de,C\!\in\!C({\cal U}_{\T}(\P;J);\R^+)$  with the following property.
For every $\ups\!=\!(b,v)\!\in\!\ti{\cal F}^1\T_{\de}$
there exists a homomorphism
$$\ti{R}_{\ups}\!: \Ga_-(b;\L)\lra  \ti\Ga_-(\ups;\L)$$
such that the map $\ups\!\lra\!\ti{R}_{\ups}$ is 
$\Aut(\T)\!\propto\!(S^1)^I$-invariant and smooth on $\ti{\cal F}^1\T_{\de}^{\eset}$.
Furthermore, the map $\ups\!\lra\ti{R}_{\ups}|_{\ti\Ga_-(b;\L)}$ is continuous
on $\ti{\cal F}^1\T_{\de}^{\eset}$ and 
\begin{equation}\label{gluing1crl1_e1}
\ti{R}_b=\id \quad\forall\,b\!\in\!{\cal U}_{\T}^{(0)}(\P;J).
\end{equation}\\
\end{crl}

\noindent
If $\ups\!\in\!\ti{\cal F}^1\T_{\de}^{\eset}$,
the homomorphism $\ti{R}_{\ups}$ is defined~by
$$\ti{R}_{\ups}\xi=\Pi_{\ze_{\ups,0}}\ti\pi_{\ups;-}R_{\ups,0}\xi
\qquad\forall\,\xi\!\in\!\Ga_-(b;\L).$$
Since the maps
$$\ups \lra b_0(\ups),\, \ze_{\ups,0},\, R_{\ups,0},\, \Ga_-\big(\ups;\L;[\rho(\ups)]\big)$$
are continuous over $\ti{\cal F}^1\T_{\de}\!-\!\rho^{-1}(0)$,
this family of homomorphisms extends continuously over  
$\ti{\cal F}^1\T_{\de}\!-\!\rho^{-1}(0)$,
as can be seen by an argument similar to Subsection~3.9 and~4.1 in~\cite{Z3}.
This extension is formally described in the same way as the homomorphisms 
$\ti{R}_{\ups}$ for $\ups\!\in\!\ti{\cal F}^1\T_{\de}^{\eset}$.
On the other hand, if $\rho(\ups)\!=\!0$, we~put
$$\ti{R}_{\ups}\xi=\Pi_{\ze_{\ups,0}}R_{\ups,0}\xi=R_{\ups,0}\xi
\qquad\forall\,\xi\!\in\!\ti\Ga_-(b;\L).$$
The second equality above holds by~\e_ref{mapgluing_e4a}.
By~\e_ref{gluing0crl_e3}, the requirement~\e_ref{gluing1crl1_e1} is satisfied.\\

\noindent
It remains to check that the extension described above is continuous at every
$$\ups^*\!\equiv\!(b^*,v^*) \in \ti{\cal F}^1\T_{\de}\!\cap\!\rho^{-1}(0).$$
We note that by~\e_ref{g1map_e10},
\begin{equation}\label{g1map_e12a}
\ti{R}_{\ups}\xi=\Pi_{\ze_{\ups,0}}\big(R_{\ups,0}\xi\!+\!\ve_{\ups,0}(\xi)\big)
\qquad\forall\,\xi\!\in\!\Ga_-(b;\L),\,
\ups\!\in\!\ti{\cal F}^1\T_{\de}^{\eset},
\end{equation}
for some homomorphism
$$\ve_{\ups,0}\!: \Ga_-(b;\L)\lra \Ga(\ups;\L)$$
such that
\begin{equation}\label{g1map_e12b}
\big\|\ve_{\ups,0}(\xi)\big\|_{\ups,p,1}\le C(b)\big|\rho(\ups)\big|\|\xi\|_{b,p,1}
\qquad\forall\, \xi\!\in\!\Ga_-\big(b;\L;[\rho(\ups)]\big).
\end{equation}
Suppose $\ups_r\!\!\equiv\!(b_r,v_r)\in\!\ti{\cal F}^1\T_{\de}^{\eset}$ 
and $\xi_r\!\in\!\ti\Ga(b_r;\L)$ are sequence such that
$$\lim_{r\lra\i}\ups_r=b^* \qquad\hbox{and}\qquad
\lim_{r\lra\i}\xi_r=\xi^*\in\Ga_-(b^*;\L).$$
Since $\ti\Ga(b_r;\L)\!\subset\!\Ga_-(b_r;\L;[\rho(\ups_r)])$ and 
the maps
$$\ups\lra b_0(\ups),\, \ze_{\ups,0},\, R_{\ups,0}$$
are continuous over $\ti{\cal F}^1\T_{\de}$,
$$\lim_{r\lra\i}\ti{R}_{\ups_r}\xi_r=
\lim_{r\lra\i}\Pi_{\ze_{\ups_r,0}}\big(R_{\ups_r,0}\xi\!+\!\ve_{\ups_r,0}(\xi)\big)
=\ti{R}_{\ups}\xi^*,$$
by~\e_ref{g1map_e12a} and~\e_ref{g1map_e12b}, as needed.\\

\noindent
Corollary~\ref{gluing1_crl1} concludes the proof of Lemma~\ref{g1conebdstr_lmm}.
It remains to finish the proof of Proposition~\ref{g1conebdstr_prp}.
By Corollary~\ref{gluing1_crl1}, $\ti{R}_{\ups}$ induces an injective homomorphism
$$R_{[\ups]}\!:{\cal V}_{1,k;{\T}}^{d;m}\big|_b
\lra{\cal V}_{1,k}^d\big|_{\phi_{\T}([\ups])}$$
for $b\!\in\!{\cal U}_{\T;1}^m(\P;J)$ and
$[\ups]\!=\![b,v]\!\in\!{\cal F}^1{\cal T}_{\de}$.
If $U$ is an open subset of ${\cal U}_{\T}(\P;J)$ and ${\cal W}\!\lra\!U$ 
is a smooth subbundle of ${\cal V}_{1,k}^d|_U$  such~that
$${\cal W}_b\subset{\cal V}_{1,k;\T}^{d;m}\big|_b
\qquad\forall~ b\in U\!\cap\!{\cal U}_{\T;1}^m(\P;J),~
m\in\big(\!\max(|\chi(\T)|\!-\!n,1),|\chi(\T)|\big),$$
then the map $[\ups]\!\lra\!R_{[\ups]}$ induces 
a continuous injective bundle homomorphism
$$\ti\phi_{\cal W}\!: \pi_{{\cal F}^1\T_{\de}|_U}^*{\cal W}
\lra{\cal V}_{1,k}^d$$
that restricts to the identity over $U$ and is smooth over ${\cal F}^1\T_{\de}^{\eset}$.\\

\noindent
Finally, for each $m\!\in\!(\max(|\chi({\cal T})|-\!n,1),|\chi({\cal T})|)$,
let $U_{\cal T}^m\!\subset\!U_{\cal T}$ be a small neighborhood of
${\cal U}_{{\cal T};1}^m(\P;J)$ in $\X_{1,k}(\P,d)$ and~let
$${\cal W}_{1,k;{\cal T}}^{d;m}\lra {\cal U}_{\cal T}(\P;J)\cap U_{\cal T}^m$$
be a subbundle of ${\cal V}_{1,k}^d$ such that
\begin{gather}
\label{bundleextcond_e1}
{\cal W}_{1,k;\T}^{d;m}\big|_{{\cal U}_{{\cal T};1}^m(\P;J)}
={\cal V}_{1,k;\T}^{d;m}\big|_{{\cal U}_{{\cal T};1}^m(\P;J)}
\qquad\hbox{and}\\
\label{bundleextcond_e2}
{\cal W}_{1,k;\T}^{d;m}\big|_b\subset {\cal V}_{1,k;\T}^{d;m'}\big|_b
\qquad \forall~ b\in{\cal U}_{\T;1}^{m'}(\P;J)\!\cap\!U_{\T}^m,~
m'\in\big(\!\max(|\chi({\cal T})|\!-\!n,1),|\chi({\cal T})|\big).
\end{gather}
By the next paragraph, such an extension of 
${\cal V}_{{1,k};{\cal T}}^{d;m}\big|_{{\cal U}_{{\cal T};1}^m(\P;J)}$
to ${\cal U}_{\cal T}(\P;J)\!\cap\!U_{\cal T}^m$ exists if $U_{\cal T}^m$
is sufficiently small.
By the previous paragraph, the bundle homomorphism
$$\ti\phi_{\T}^m\!\equiv\!\ti\phi_{{\cal W}_{1,k;\T}^{d;m}}\!:
\pi_{{\cal F}^1\T_{\de}|_U}^*{\cal W}_{1,k;\T}^{d;m} \lra\V_{1,k}^d$$
is continuous and injective, restricts to the identity over 
${\cal U}_{\T;1}(\P;J)\!\cap\!U_{\T}$, and is smooth over ${\cal F}^1\T_{\de}^{\eset}|_U$.
We define the bundle 
$${\cal V}_{{1,k};{\cal T}}^{d,m}\lra\ov\M_{1,k}^0(\P,d;J)\cap U_{\cal T}^m$$
to be the image of $\tilde{\phi}_{\cal T}^m$.
This bundle has the claimed rank by the last statement of Lemma~\ref{g1conebdstr_lmm}.
The last condition of Proposition~\ref{g1conebdstr_prp} is satisfied
by the definition of the bundles 
${\cal V}_{1,k;{\cal T}}^{d;m}|_{{\cal U}_{{\cal T};1}^m(\P;J)}$
following Proposition~\ref{g1conebdstr_prp}.
The proof of Proposition~\ref{g1conebdstr_prp} is now complete.\\

\noindent
We now prove the extension claim used in the previous paragraph.
By definition,
$$\ti\F^1\T=\big\{w\!\in\!\ti\F\T\!: \cD_{\T}w\!=\!0\big\}.$$
Since $\cD_{\T}$ is a continuous bundle section, if $\ti{U}$ is a sufficiently small
neighborhood of $\ti\U_{\T;1}^m(\P;J)$ in $\U_{\T}^{(0)}(\P;J)$,
there exists a vector bundle $\ti\F^1\T^m\!\lra\!\ti{U}$ such~that
\begin{equation}\label{bundleextdfn_e}
\ti\F^1\T^m\big|_{\ti\U_{\T;1}^m(\P;J)}=\ti\F^1\T\big|_{\ti\U_{\T;1}^m(\P;J)}
\qquad\hbox{and}\qquad
\ti\F^1\T|_{\ti{U}}  \subset \ti\F^1\T^m\subset \ti\F\T.
\end{equation}
The neighborhood $\ti{U}$ and the bundle $\ti\F^1\T^m$
can be chosen so that they are preserved
by the actions of $\Aut(\T)\!\propto\!(S^1)^I$.
We then define the vector bundle ${\cal W}_{1,k;\T}^{d;m}\!\lra\!U$ by
\begin{gather*}
U=\big\{[b]\!\in\!{\cal U}_{\T;1}(\P;J)\!: b\!\in\!\ti{U}\big\}
\qquad\hbox{and}\\
{\cal W}_{1,k;\T}^{d;m}=\big\{ [\xi]\!\in\!\V_{1,k}^d|_b \!:
b\!\in\!U;~ \D_{\cal T}(\xi\!\otimes\!w)\!=\!0
~\forall\,w\!\in\!\ti\F^1\T^m|_b\big\}.
\end{gather*}
By the same argument as at the end of Subsection~\ref{g1conelocalstr_subs2},
${\cal W}_{1,k;{\T}}^{d;m}\!\lra\!U$ is a vector bundle of rank $da\!+\!1\!-\!m$.
By the middle statement of Lemma~\ref{g1conebdstr_lmm} and~\e_ref{bundleextdfn_e},
this vector bundle satisfies the requirements~\e_ref{bundleextcond_e1}
and~\e_ref{bundleextcond_e2}, as~needed.\\

\vspace{.2in}

{\it 
\begin{tabbing}
${}\qquad$
\= Department of Mathematics, Stanford University, Stanford,
CA 94305-2125\\ 
\> azinger@math.stanford.edu
\end{tabbing}}

\end{document}